\documentclass[11pt, letterpaper]{amsart}
\usepackage[margin=2.8cm]{geometry}
\pdfoutput=1
\usepackage[utf8]{inputenc}

\usepackage[english]{babel}
\usepackage[autostyle, english = american]{csquotes}
\MakeOuterQuote{"}

\usepackage{amsmath}
\usepackage{amsfonts}
\usepackage{amssymb}
\usepackage{amsthm}
\usepackage{thmtools}
\usepackage{mathrsfs}
\usepackage{mathabx, bbm}
\usepackage{stmaryrd}
\usepackage[all]{xy}
\usepackage{multirow}
\usepackage{graphicx}
\usepackage[usenames,dvipsnames,svgnames,table]{xcolor}
\usepackage{transparent}         
\usepackage{float}

\usepackage{tikz, comment}
\usetikzlibrary{cd, fit, knots, positioning, arrows.meta, matrix, mindmap, trees, shadows, shapes, calc, fadings, decorations.pathreplacing, intersections, shapes.callouts, shapes.arrows, shapes.geometric, decorations.pathmorphing}

\usepackage[bookmarks=true]{hyperref}
\usepackage{cleveref}

\definecolor{azulmaneiro}{rgb}{0, 0.125, 0.666}
\definecolor{brownmaneiro}{rgb}{0.4, 0.2, 0}
\hypersetup{colorlinks=true, linkcolor=azulmaneiro, citecolor=brownmaneiro, urlcolor=red}

\numberwithin{equation}{section}

\theoremstyle{plain}
\newtheorem{mainthm}{Theorem}

\newtheorem{maincor}[mainthm]{Corollary}

\newtheorem{thm}{Theorem}[section]
\newtheorem{lem}[thm]{Lemma}
\newtheorem{prop}[thm]{Proposition}
\newtheorem{cor}[thm]{Corollary}

\newtheorem{fact}[thm]{Fact}

\theoremstyle{definition}
\newtheorem{defn}[thm]{Definition}
\newtheorem{ex}[thm]{Example}

\newtheorem{notat}[thm]{Notation}

\theoremstyle{remark}
\newtheorem{rem}[thm]{Remark}

\theoremstyle{definition}
\newtheorem*{defn*}{Definition}
\newtheorem*{overw*}{\,\,Overview of the results}

\newcommand{\Z}{\mathbb{Z}}

\newcommand{\Q}{\mathbb{Q}}
\newcommand{\C}{\mathbb{C}}
\newcommand{\R}{\mathbb{R}}
\newcommand{\F}{\mathbb{F}}

\DeclareMathOperator{\hH}{H} 

\DeclareMathOperator{\ct}{\mathrm{ct}} 

\DeclareMathOperator{\Ext}{Ext}

\DeclareMathOperator{\dd}{d}

\newcommand{\FP}{\mathtt{FP}} 

\newcommand{\tF}{\mathtt{F}}

\DeclareMathOperator{\image}{Im}

\DeclareMathOperator{\End}{End}

\DeclareMathOperator{\op}{op}

\DeclareMathOperator{\stab}{stab}
\DeclareMathOperator{\Aut}{Aut}

\DeclareMathOperator{\Lk}{Lk}

\DeclareMathOperator{\opH}{\operatorname{H}}
\DeclareMathOperator{\frC}{\mathfrak{C}}
\DeclareMathOperator{\Spec}{\mathrm{Spec}}

\DeclareMathOperator{\SL}{SL}

\newcommand{\caG}{\mathcal{G}} 

\newcommand{\Sym}{\mathrm{Sym}}
\newcommand{\barG}{\overline{G}}

\newcommand{\Ch}{\mathrm{Ch}} 

\newcommand{\Ce}{\mathrm{Ce}}
\DeclareMathOperator{\type}{\mathrm{type}}

\newcommand{\veps}{\varepsilon} 
\newcommand{\wtil}[1]{\widetilde{#1}} 
\newcommand{\LL}{\langle\!\langle} 
\newcommand{\GG}{\rangle\!\rangle} 

\newcommand{\bbone}{\mathbbm{1}}

\newcommand{\caA}{\mathcal{A}}

\newcommand{\caH}{\mathcal{H}}

\newcommand{\caR}{\mathcal{R}}
\newcommand{\caC}{\mathcal{C}}
\newcommand{\caE}{\mathcal{E}}
\newcommand{\caB}{\mathcal{B}}

\newcommand{\caT}{\mathcal{T}}
\newcommand{\caU}{\mathcal{U}}

\newcommand{\caD}{\mathcal{D}}

\newcommand{\qmin}{q_{\min}}
\newcommand{\qmax}{q_{\max}}
\newcommand{\qminJ}{q_{\min,J}}
\newcommand{\qmaxJ}{q_{\max,J}}
\DeclareMathOperator{\Res}{\mathrm{Res}}
\DeclareMathOperator{\proj}{\mathrm{proj}}

\DeclareMathOperator{\opE}{\mathbf{E}}
\newcommand{\co}{$(\mathrm{co})$}

\newcommand{\Tk}{[T_k]}

\newcommand{\argu}{\hbox to 1.5ex{\hrulefill}}  

\usepackage[style=alphabetic,sorting=nyt,maxbibnames=7,maxcitenames=5,block=space,backend=bibtex]{biblatex}
\title[Higher properties T and finiteness properties]{Higher properties T and finiteness properties\\ for twin building lattices}
\author[Higher properties T and finiteness properties]{Bianca Marchionna} 
\address{Université catholique de Louvain (IRMP), Chemin du Cyclotron 2, 1348 Louvain-la-Neuve, Belgium}
\email{bianca.marchionna@uclouvain.be}
\date{\today}

\subjclass{20G44, 20J05, 20J06, 22D10, 22D55}
\keywords{twin building lattices, Kac--Moody groups, finiteness properties, higher property T}

\begin{document}
\maketitle
\begin{abstract}
Properties~$[T_k]$ and finiteness properties~$\FP_k(\Q)$ and~$\tF_k$ are refinements of the following fundamental notions for groups: Kazhdan property~T in the first case, and finite generation/presentation in the last two cases.
Under suitable hypotheses on the thickness, we prove that every twin building lattice of $(k+1)$-spherical type~$(W,S)$, where~$k\geq 1$, has property~$[T_k]$ and is of type~$\FP_k(\Q)$. If~$(W,S)$ is classical, that group is even of type~$\tF_k$. To the best of the author's knowledge, this produces the first examples of infinite finitely presented simple groups with property~$[T_k]$ for some~$k\geq 2$, and the first examples of finitely presented simple Kazhdan groups that are Frobenius stable and not Frobenius approximable. 
Along the way, we establish new criteria for property~$[T_k]$ for groups acting on CW complexes that extend results of Dymara--Januszkiewicz and Grinbaum-Reizis--Oppenheim. As a key step in our proofs, we determine high acyclicity and, in classical type, high connectivity properties for the opposite complex of a Moufang twin building of $(k+1)$-spherical type and suitable finite thickness.
\end{abstract}

\setcounter{tocdepth}{1}
\tableofcontents

\section{Introduction}\label{s:intro}

\subsection{Overview of the paper}\label{sus:overv}
 Properties~$\Tk$ and~$(T_k)$ (in short, higher properties~T) and
finiteness properties~$\FP_k$ and~$\tF_k$ form families of properties for groups labelled by a positive integral index~$k$ that refine the notion of Kazhdan property~T for the first two families, and the concepts of finite generation and finite presentation for the last two. 

Recall that a locally compact second countable (lcsc) group~$G$ has \emph{property~$\Tk$} if its continuous cohomology~$\hH_{\ct}^\ell(G,\rho)$ over every strongly continuous unitary representation~$\rho\colon G\to \caU(V)$ vanishes in degree~$1\leq \ell\leq k$, and~$G$ has \emph{property~$(T_k)$} if the same vanishing holds only for all those~$\rho$ that have no non-zero invariant vectors~\cite[Definition~1.1]{basa23}. 
By Delorme--Guichardet theorem, both property~$[T_1]$ and~$(T_1)$ are equivalent to property~T. Clearly, property~$[T_k]$ always implies property~$(T_k)$; for~$k\geq 2$, however, the converse does not hold in general~\cite[Remark~4]{BaSau:survey}.
Recent advances~\cite{DGLT:stab, basa23, bbbs:ban} have brought up increased attention to the study of higher properties~T, due to some notable connections to group stability, high dimensional expanders, and fixed point theorems. See~\cite{bader:ICM, BaSau:survey} for an overview.

On the other hand, a group~$G$ is of \emph{type~$\tF_k$} if there is an Eilenberg--MacLane space~$K(G,1)$ with finite $k$-skeleton, and~$G$ is of \emph{type~$\FP_k(R)$} for a non-zero unital ring~$R$, if the trivial $RG$-module~$R$ admits a projective resolution that is finitely generated in dimension~$\leq k$. Both type~$\tF_1$ and type~$\FP_1(R)$ are equivalent to finite generation, and type~$\tF_2$ is equivalent to finite presentation. It is well known that type~$\tF_k$ always implies~$\FP_k(R)$ for every~$R$ although, for~$k\geq 2$, the converse does not hold in general~\cite{StallingsSB2, bieri, BestvinaBrady}.
Finiteness properties are not only a well-established tool to "organise and separate groups" but also make the (co)homology theory of groups having them "better behaved" (e.g., for groups of type~$\FP_k$ for every~$k$ we have a well-defined notion of Euler--Poincaré characteristic~\cite{chis:euler}).

Linear algebraic groups over local fields and their $S$-arithmetic lattices have so far constituted an important "test class" for the study of higher properties~T and finiteness properties, which culminated with recent results by Bader--Sauer~\cite{basa23, BaSau:survey} in the first case, and with the Rank Theorem by Bux--K\"ohl--Witzel~\cite{bkw:rank} in the second. 

It is common to try to extend results known for~$S$-arithmetic lattices to \emph{twin building lattices}. The latter groups are suitable irreducible lattices in the product~$\Aut(\caC_+)\times\Aut(\caC_-)$ of the automorphism groups of the two halves~$\caC_+$ and~$\caC_-$ of a Moufang\footnote{I.e., a thick twin building such that, for every twin root~$\alpha$, the associated root group acts transitively on the set of all twin apartments containing~$\alpha$~\cite[Definition~8.18]{ab:build}.} twin building of type~$(W,S)$ and finite thickness vector.
Twin building lattices comprise some $S$-arithmetic groups but also minimal Kac--Moody groups over finite fields modulo their centre. 
Twin building lattices of Euclidean\footnote{A Coxeter system~$(W,S)$ is of \emph{Euclidean type} if all irreducible subsets~$J\subseteq S$ are spherical or affine.} type~$(W,S)$ share important structural properties with $S$-arithmetic lattices, like being linear and residually finite. In contrast, Caprace--Rémy~\cite{caprem:KMsup} proved that, under mild hypotheses, twin building lattices of irreducible non-Euclidean type are simple and non-residually finite. In particular, being finitely generated, the latter groups are non-linear. See \Cref{thm:TBLsimp} and \Cref{rem:KMsimp} for details.

In this paper, we focus on twin building lattices of \emph{$(k+1)$-spherical type}~$(W,S)$ for some~$k\in \Z_{\geq 1}$, i.e., we require that all subsets~$J\subseteq S$ with~$|J|\leq k+1$ generate a finite subgroup of~$W$. 
Within suitable assumptions on the thickness vector of the associated twin building, every twin building lattice~$\Lambda$ of $(k+1)$-spherical type~$(W,S)$ satisfies the following:
\begin{enumerate}
    \item[(A)] \label{itemA} $\Lambda$ has \emph{property~$\Tk$} and if, in particular~$|S|=k+2$ and~$W$ is infinite, $\Lambda$~does not have property~$(T_{k+1})$ (see \Cref{sus:introhT});
    \item[(B)] \label{itemB} $\Lambda$ is of \emph{type~$\FP_k(\Q)$} and, if~$(W,S)$ is classical\footnote{I.e., all irreducible spherical special subgroup~$(W_J,J)$ is of type~$A_n$, $B_n$, $C_n$ or~$D_n$ for some~$n\geq 0$.}, it is even of~\emph{type~$\tF_k$} (see \Cref{sus:introFP}).
\end{enumerate}
The results in~\hyperref[itemA]{(A)} are, in fact, applications of new general criteria for property~$\Tk$ for groups acting on CW complexes recorded in \Cref{s:mainCoho}. 

As a key step towards these results, we determine high acyclicity and, in classical types, high connectivity properties of the so-called "negative opposite complex" of the associated twin building (see \Cref{sus:introhT}). This addresses a long-standing open problem in building theory: for example, Abramenko~\cite[Conjecture~1]{ab:twin} conjectured that every spherical Moufang building of type~$(W,S)$ and "suitably large thickness" is~$(|S|-2)$-connected. As a consequence of our \Cref{Y-conn}, we prove the $\Q$-homological counterpart of this conjecture.

Thanks to the results listed in~\hyperref[itemA]{(A)} and~\hyperref[itemB]{(B)}, we are able to produce:
\begin{itemize}
    \item[(a)] \label{item-a} the first examples of infinite finitely presented simple groups that have property~$\Tk$ for some~$2\leq k\leq 7$. This generalises an analogous statement for property~T recorded in~\cite[Corollary~21]{caprem:KMsup}. The given upper bound on~$k$ is motivated in \Cref{rem:Sph8}; 
    \item[(b)] \label{item-b} the first examples of finitely presented simple Kazhdan groups that are Frobenius stable and not Frobenius approximable (see \Cref{sus:introFrob}).
\end{itemize} 
In the next subsections of the introduction, we give a detailed statement of our main results and frame them in the pre-existing literature.

\smallskip

We conclude this first introductory part with a disclaimer and a glimpse into the author's approach followed in the preparation of this paper.
Firstly, we find important to mention that this paper does \emph{not} (at least directly) contribute to the solution of the following open problem: determine which finiteness properties are enjoyed by groups with higher properties~T (see~\cite[Example~14, Question~15]{BaSau:survey} for details). 
The "connection" between higher properties~T and finiteness properties that we brought up here is rather at the level of methods: key steps towards the results in~\hyperref[itemA]{(A)} and~\hyperref[itemB]{(B)} have been established by exporting ideas and techniques typically used to detect finiteness properties to prove results on higher properties~T and vice-versa.
For instance, inspired by some well-known criteria for finiteness properties by Brown~\cite{BrownFP}, we establish a sort of "Brown's criterion" for higher properties~T (compare, for instance, \Cref{prop:cohoPlainCC} and \Cref{prop:BrownSuff}), and this constitutes a crucial step towards~\hyperref[itemA]{(A)}.
Moreover, a well-known technique in the study of higher properties~T and high dimensional expanders (see \Cref{thm:OppVan}) plays a key role in determining high acyclicity properties of the opposite complex of a Moufang twin building. This result, in turn, is the crucial new "piece" that we used to get~\hyperref[itemB]{(B)}.

\subsection{The object of study: groups with locally finite root group data}\label{sus:obj}
Rather than just twin building lattices, our results refer more generally to groups~$G$ with a \emph{locally finite root group datum of type~$(W,S)$}, i.e., a family of finite subgroups~$(U_\alpha)_{\alpha\in \Phi}$ of~$G$ indexed by the roots~$\Phi$ of a given Coxeter system~$(W,S)$, which generates\footnote{This requirement is usually relaxed by asking that~$G$ is generated by~$\{U_\alpha\}_{\alpha\in \Phi}$ and a subgroup~$T\leq G$ normalising each~$U_\alpha$. This discrepancy has no impact for the results of the paper, and the version of the definition we use allows us a shorter formulation of Caprace--Rémy's simplicity result; see \Cref{fact:defnRGD} and the line afterwards.}~$G$ and satisfies suitable compatibility axioms (see \Cref{defn:RGD}). These axioms guarantee the existence of a \emph{twin BN-pair}~$(B_+,B_-,N)$ of type~$(W,S)$ for~$G$, which consists of two "suitably interlaced" BN-pairs~$(B_+,N)$ and~$(B_-,N)$ of type~$(W,S)$ (see \Cref{sus:BNpairs}). In our case, both BN-pairs have common finite thickness vector~$(q_s)_{s\in S}$ that satisfies~$q_s\geq 2$ for every~$s\in S$. In other words, the simplicial building~$X_\pm=\Delta(G,B_\pm)$ associated to the BN-pair~$(B_\pm,N)$ is thick. For future reference, we also set
\begin{equation*}
   \qmin:=\min\{q_s\mid s\in S\}\qquad\text{and}\qquad \qmax:=\max\{q_s\mid s\in S\}.
\end{equation*}
The class of twin buildings~$(X_+,X_-)$ described before is fairly broad: every Moufang twin building with finite thickness vector arises from the construction described before (see \Cref{rem:Mouf}).

If~$B_+\cap B_-$ is finite and~$\qmin>|S|$, the group~$G/Z(G)$ embeds diagonally as an irreducible lattice in~$\barG_+\times\barG_-$, where~$\barG_\pm$ is the completion of the action of~$G$ on~$X_\pm$. The \emph{twin building lattice} associated to~$(G,(U_\alpha)_{\alpha\in \Phi})$ is precisely the image of~$G/Z(G)$ under this embedding. Remarkably, since the centre~$Z(G)$ is finite,~$G/Z(G)$ and~$G$ have the same higher properties~T and the same finiteness properties (see \Cref{thm:Tkperm} and \Cref{prop:permFP}). 
We refer to \Cref{sus:TBL} for details.

Another notable group on which we focus is the subgroup of~$G$ given by~$U_+:=\langle U_\alpha\mid \alpha\in\Phi^+\rangle$, where~$\Phi^+$ is the set of all positive roots of~$\Phi$. Sometimes, $U_+$ is called the \emph{positive unipotent subgroup}, given its analogy with the unipotent radical of a (positive) Borel subgroup in a linear algebraic group.
The subgroup~$U_+$ admits a "controlled"\footnote{See \Cref{cor:Y-Uchar} for a more precise statement.} cobounded action on the so-called \emph{negative opposite complex}~$Y_-$, i.e., the subcomplex of~$X_-$ spanned by all chambers of~$X_-$ opposite to the fundamental chamber of~$X_+$; see \eqref{eq:Ypm}. 

We conclude with a concrete well-known example of what presented before; see Examples~\ref{ex:genCart} and~\ref{ex:KMgroups} for details. Every Kac--Moody group~$\caG_A(\F_q)$ over a finite field~$\F_q$ and generalised Cartan matrix~$A=[A_{st}]_{s,t\in S}$ canonically admits a family of subgroups~$(U_\alpha)_{\alpha\in \Phi_A}$, each of them isomorphic to the additive group~$\F_q$, that form a root group datum (in the sense mentioned before) for the finite-index subgroup~$\caG_A(\F_q)^\dagger:=\langle U_\alpha\mid\alpha\in \Phi_A\rangle$. The type of this root group datum is the Weyl group~$(W_A,S)$ of~$A$ and, in this case,~$q_s=q$ for every~$s\in S$.
For~$q>\dim(A)$, we denote by~$\Lambda_A(\F_q):=\caG_A(\F_q)/Z(\caG_A(\F_q)^\dagger)$ the associated \emph{Kac--Moody lattice} of type~$A$ over~$\F_q$. 
If~$(W_A,S)$ is $3$-spherical and~$q\geq 5$, the corresponding group~$U_+$ recovers the \emph{Kac--Moody--Steinberg group}~$\caU_A(\F_q)$, another notable kind of groups that has featured in recent constructions of high-dimensional expanders~\cite{GdP-VP:KMS, OpVal2, OpVal1}. 

\subsection{Acyclicity and connectivity properties of the opposite complex}\label{sus:introY-}
The crucial step towards our results on higher property~T and finiteness properties is the following theorem, which might be of independent interest:
\begin{mainthm}\label{Y-conn}
    Let~$G$ be a group with a locally finite root group datum~$(U_\alpha)_{\alpha\in \Phi}$ of $(k+1)$-spherical type~$(W,S)$ for some~$1\leq k\leq |S|-1$. Assume that
     \begin{equation}\label{Cond3}
        \frac{\qmin}{\sqrt{\qmax}}>2(|S|-1).
    \end{equation}
    Then~$Y_-$ is $(k-1)$-acyclic over~$\Q$. If additionally~$\qmin\geq 2^{2|S|-1}$ and~$(W,S)$ is classical, then~$Y_-$ is $(k-1)$-connected.
\end{mainthm}
We refer to \Cref{thm:Yacyc} and \Cref{thm:Y-homot} for slightly more general statements. Note that \Cref{Y-conn} can be equivalently formulated for the~"$Y_-$"-complex of any Moufang twin building of $(k+1)$-spherical type and finite thickness (see \Cref{rem:Mouf}).

To prove \Cref{Y-conn}, we considered a well-known filtration~$Y_-=X_0\subseteq X_1\subseteq \dots \subseteq X_n\subseteq \dots$ of the negative simplicial building~$X_-$ introduced by Abels~\cite{abels} (see \Cref{defn:Abfiltr}).  
This filtration has a notable property: every~$X_d$,~$d\geq 1$, has the form
\begin{equation*}
    X_d=X_{d-1}\cup \bigcup_{i\in I_d}Y_{id}
\end{equation*}
for some contractible complexes~$Y_{id}$, $i\in I_d$, and~$Y_{id}\cap Y_{jd}\subseteq X_{d-1}$ for all~$i\neq j$ in~$I_d$. For simplicity, we call the various~$X_{d-1}\cap Y_{id}$ the \emph{descending complexes} of the filtration.
By standard arguments (see \Cref{prop:filtHom}), the key part towards \Cref{Y-conn} is to prove that each descending complex is $(k-1)$-acyclic over~$\Q$ (resp.~$(k-1)$-connected). Here is a quick overview of our proof for the latter claim:
\begin{enumerate}
    \item[(1)] Extending a strategy by Abramenko~\cite{ab:twin} in the case when~$(W,S)$ is $(|S|-1)$-spherical, we describe the descending subcomplexes in terms of "$Y_-$-subcomplexes" of smaller rank (see \Cref{prop:descLinks}). This allows us to argue inductively on the difference~$|S|-\kappa(W,S)$, where~$\kappa(W,S)$ is the largest integer~$1\leq \ell\leq |S|$ such that~$(W,S)$ is~$\ell$-spherical. See the first part of the proof of \Cref{thm:Yacyc} for details.
    \item[(2)] The base case of the previously mentioned induction  is when~$(W,S)$ is spherical. In this specific case, the connectivity statement was proved by Abramenko in~\cite[\S II.8, Corollary~19]{ab:twin}, and this motivates the additional restrictions on~$\qmin$ and~$(W,S)$ in our theorem. The acyclicity part uses a completely different strategy: since the spectral properties of the $1$-dimensional links of~$Y_-$ are well understood, the claim follows by invoking Oppenheim's upgraded version of Garland's vanishing theorem for finite simplicial complexes (see \Cref{thm:OppVan}). This allowed us to substantially relax the requirement on the thickness vectors and drop the hypothesis that~$(W,S)$ is classical. We refer to \Cref{prop:YconnSph} for details.
\end{enumerate}

\subsection{Higher properties~T and twin building lattices}\label{sus:introhT}
To motivate our contribution, we find it important to first retrace the main results of higher properties~T established over the years. Before this, two important remarks are required. Firstly,  
every compact group has property~$\Tk$ for all~$k\geq 1$~\cite[Lemma~3.3]{basa23}, thus the focus conveys to non-compact groups. The second remark concerns lattices: while property~T passes from an ambient lcsc~group to all its lattices, at the moment of writing it is not clear whether this is also the case for higher properties~T. What is known is that both properties~$[T_k]$ and~$(T_k)$, $k\geq 1$, are inherited by all cocompact lattices~\cite[Lemma~3.13]{basa23}. Thus, in the current state of knowledge, determining higher properties~T for lcsc groups and for their non-cocompact lattices are \emph{a priori} different problems. 

\smallskip

Although without the current terminology, the first examples of non-compact groups with higher properties~T have been known since the~$1970$s. Garland~\cite{garland} proved that, for every $p$-adic field~$K$ and sufficiently large~$p$, all cocompact lattices in a simple algebraic group~$\caG(K)$ of rank~$k+1\geq 2$ have property~$\Tk$. 

Extending Garland's result to non-cocompact lattices remained open for roughly fifty years until a recent breakthrough by Bader--Sauer~\cite{basa23}. They proved that all simple algebraic groups~$\caG(K)$ of rank~$k+1\geq 2$ on a local field~$K$ of characteristic zero, as well as \emph{all} its lattices, have property~$(T_k)$ if~$K$ is Archimedean, and even property~$[T_k]$ if~$K$ is non-Archimedean (i.e., is a $p$-adic field). They obtained an analogous result for connected semisimple Lie groups of real rank and finite centre, and all their lattices. See~\cite[Theorems~A and~B]{basa23} and~\cite[Theorem~1]{BaSau:survey} for references.
For non-Archimedean local fields of positive characteristic, a similar generalisation is open~\cite[Remark~3]{BaSau:survey}. 

The range of examples is broader when we look at possibly non-discrete lcsc groups. Dymara--Januszkiewicz~\cite[Theorem~E]{DJ:build} proved that an lcsc group~$G$ with a BN-pair~$(B,N)$ of $(k+1)$-spherical type~$(W,S)$ and~$B$ compact open subgroup has property~$\Tk$ provided the associated building has minimal thickness~$\geq\frac{1}{25}(1764)^{|S|}$. 
Their theorem is formulated for groups acting on more general contractible or finite simplicial complexes, although its main known applications concern lcsc BN-groups. Crucially, this result extends beyond BN-pairs of Euclidean type. For example, if~$G=\caG_A(\F_q)$ is a Kac--Moody group with $(k+1)$-spherical generalised Cartan matrix~$A$, its completions~$\barG_\pm$ have property~$\Tk$ whenever~$q\geq\frac{1}{25}(1764)^{\dim(A)}$. 

Dymara--Januszkiewicz's results and methods have inspired several generalisations through the years, resulting in broad criteria for property~T~\cite{ejz:propT, kassabov}, or in more general results on properties~$\Tk$ within lcsc~groups that act cocompactly on simplicial complexes~\cite{Op15:vanish, gro:curv}. 
Cornerning the latter type of results, there are two assumptions that always occur: 
\begin{enumerate} 
\item[(1)] the stabilisers of top-dimensional simplices are compact and open; 
\item[(2)] every infinite\footnote{Here "infinite" means "with infinitely many simplices".} link of the complex is contractible. \end{enumerate}
It remains unclear how to apply them to the case of an arbitrary twin building lattice~$\Lambda$. The main obstructions are the following:
\begin{enumerate}
    \item \label{obstr1} $\Lambda$ naturally comes with two BN-pairs, in fact a twin BN-pair~$(B_+,B_-,N)$ induced from the one of~$G$. However,~$B_+$ and~$B_-$ are infinite countable unless the type~$(W,S)$ is spherical, i.e., unless~$\Lambda$ is finite. Since countable locally compact groups are discrete, there is no locally compact group topology on~$\Lambda$ that makes~$B_\pm$ a compact open subgroup.
    \item \label{obstr2} $\Lambda$ is a lattice in the lcsc~group~$\barG:=\barG_+\times \barG_-$, and the previously mentioned criteria apply to~$\barG$ acting on the reducible building~$X_+\times X_-$. However, if~$(W,S)$ is not spherical,~$\Lambda$ is never cocompact in~$\barG$~\cite{remy:integrab}. Thus, we cannot automatically transfer higher properties~T from~$\barG$ to~$\Lambda$.
    \item \label{obstr3} There is a natural choice of a simplicial complex (built from the opposition complex of~$(X_+,X_-)$) on which~$\Lambda$ acts simplicially, properly, and with finitely many orbits. However, when infinite, this complex is hardly ever contractible, see \Cref{cor:homotRes}.
\end{enumerate}
A non-contractibility obstruction similar to~\eqref{obstr3} appears for the opposite complex~$Y_-$ (see~\Cref{cor:OpXhomot}). Once again, the previously mentioned criteria for property~$\Tk$ remain generally inapplicable to the $U_+$-action on~$Y_-$.

\smallskip

Motivated by these obstructions, building on results by Dymara--Januszkiewicz~\cite[Theorem~E]{DJ:build} and Grinbaum-Reizis--Oppenheim~\cite[Theorem~1.6]{gro:curv}, we obtain general criteria for property~$\Tk$ within the class of lcsc~groups acting on CW (typically simplicial) complexes, weakening the connectivity assumption for the complex and the compactness assumption for the stabilisers of top-dimensional cells.
We have collected them in \Cref{s:mainCoho}: see, in particular, \Cref{thm:CohoGX} and \Cref{thm:CohoGXPlain3} for the main results, and \Cref{cor:cohoPlainBN} and \Cref{cor:vanishBN} for specific criteria for groups with BN-pairs.
The key step in our extension is generalising the core cohomological argument in Dymara--Januszkiewicz's work~\cite[Lemma~3.5]{DJ:build}, which is in turn based on a general result by Borel--Wallach~\cite[Ch.~IX, Theorem~4.1]{bowal:coho}. \Cref{lem:cohoPlain} summarises our generalisation.

\smallskip

We now focus on twin building lattices~$\Lambda$ coming from a locally finite root group datum~$(U_\alpha)_{\alpha\in \Phi}$ of $(k+1)$-spherical type~$(W,S)$. To establish higher properties~T, we proceed as follows: we first determine property~$\Tk$ for~$U_+$, applying our criteria to the $U_+$-action on~$Y_-$.
Since~$U_+$ is a finite-index subgroup of~$B_+$, property~$\Tk$ transfers from~$U_+$ to~$B_+$. In turn, our criterion for property~$\Tk$ for BN-pairs allows us to transfer that property from~$B_+$ to~$G$ and then to the quotient~$\Lambda\cong G/Z(G)$.
In detail, we obtain the following: 
\begin{mainthm}\label{TBLTk}
  Let~$\Lambda$ be a twin building lattice coming from a locally finite root group datum~$(U_\alpha)_{\alpha\in \Phi}$ for a group~$G$ of $(k+1)$-spherical type~$(W,S)$, for some~$k\geq 1$. Assume that
  \begin{equation}\label{Cond}
      \frac{\qmin-1}{\sqrt{\qmax}}>2(|S|-1).
  \end{equation}
  Then:
  \begin{enumerate}
      \item both~$U_+$ and~$G$ have property~$\Tk$. If~$(W,S)$ is not $(k+2)$-spherical then~$U_+$ is not of type~$[T_{k+1}]$. Moreover, if~$|S|=k+2$ and~$(W,S)$ is not spherical, then~$G$ does not have property~$(T_{k+1})$.
      \item  $\Lambda$ has property~$\Tk$. If~$|S|=k+2$ and~$(W,S)$ is not spherical, then~$\Lambda$ does not have property~$(T_{k+1})$.
  \end{enumerate}
\end{mainthm}
See \Cref{thm:RGDpropTk}, \Cref{cor:twinTk} and \Cref{thm:U+nonTk} for slightly more general statements.
In the Kac--Moody case, \Cref{TBLTk} implies the following:
\begin{maincor}\label{KM(S)Tk}
  Let~$A=[A_{st}]_{s,t\in S}$ be a generalised Cartan matrix whose Weyl group~$(W_A,S)$ is $(k+1)$-spherical for some~$k\geq 1$. Let also~$q>3(|S|-1)^2+1$ be a prime power.
  Then: 
  \begin{enumerate}
      \item \label{KMTK} $\Lambda_A(\F_q)$ has property~$\Tk$ and, if~$|S|=k+2$ and~$(W_A,S)$ is not spherical, it does not have property~$(T_{k+1})$;
      \item \label{KMSTK} If~$k\geq 2$ and~$(W_A,S)$ is not $(k+2)$-spherical,~$\caU_A(\F_q)$ has property~$\Tk$ but not~$[T_{k+1}]$. 
  \end{enumerate}
\end{maincor}
\noindent See Corollaries~\ref{cor:KMpropTk} and~\ref{cor:KMSTk} for slightly more general statements.

\subsection{Property~$[T_2]$ and Frobenius stability}\label{sus:introFrob}
One of the first motivations for considering higher properties~T has been certain group stability problems~\cite{DGLT:stab}. Roughly speaking, a sample problem is to determine which countable groups can be approximated by asymptotic homomorphisms to symmetric groups~$\Sym(n)$ or to finite-dimensional unitary groups~$U(n)$, with respect to various norms on the latter groups. This leads, for instance, to the notions of \emph{sofic}, \emph{hyperlinear} and \emph{Frobenious approximable} groups. We refer to~\cite{arzh:asy} or the introduction of~\cite{DGLT:stab} for an overview.

What is particularly relevant in connection with higher properties~T is the method used by De Chiffre--Glebsky--Lubotzky--Thom~\cite{DGLT:stab} to construct the first examples of finitely presented groups that are not Frobenius approximable. They proved, more generally, that every finitely presented group~$G$ that satisfies~$\hH_{\ct}^2(G,\rho)=\{0\}$ for every unitary representation~$(\rho,V)$ is \emph{Frobenius stable}, which roughly means that every asymptotic homomorphism from~$G$ to finite-dimensional unitary groups is closed to an actual homomorphism in the Frobenius norm distance~\cite[Theorem~5.1]{DGLT:stab}. In particular, provided they are also non-residually finite, all those groups are not Frobenius approximable~\cite[Corollary~5.4]{DGLT:stab}. 

Combining~\cite[Theorem~5.1, Corollary~5.4]{DGLT:stab}, the simplicity theorem on twin building lattices by Caprace--Rémy~\cite{caprem:KMsup} (recalled in \Cref{thm:TBLsimp}) and \Cref{KM(S)Tk}, we deduce the following.
\begin{maincor}\label{TBLstab}
    Let~$\Lambda$ be a twin building lattice of irreducible, non-Euclidean, $3$-spherical type~$(W,S)$ whose associated root groups are finite nilpotent and whose thickness vector~$(q_s)_{s\in S}$ satisfies~\eqref{Cond}. Then the finitely presented simple Kazhdan group~$\Lambda$ is Frobenius-stable and not Frobenius approximable.
\end{maincor}
A concrete example of a group satisfying the hypotheses on \Cref{TBLstab} is the Kac--Moody lattice~$\Lambda_A(\F_q)$ with~$q>3(|S|-1)^2+1$ and generalised Cartan matrix~$A=[A_{st}]_{s,t\in S}$ of compact hyperbolic type and dimension~$\geq 4$, e.g., 
\begin{equation*}
   \begin{bmatrix}
         \phantom{-}2 & -2 & \phantom{-}0 & -1\\
        -1 & \phantom{-}2 & -1 & \phantom{-}0\\
        \phantom{-}0 & -1 & \phantom{-}2 & -1\\
        -1 & \phantom{-}0 & -1 & \phantom{-}2
    \end{bmatrix}
   \qquad \text{or}\qquad 
   \begin{bmatrix}
         \phantom{-}2 & -2 & \phantom{-}0 & -1\\
        -1 & \phantom{-}2 & -1 & \phantom{-}0\\
        \phantom{-}0 & -1 & \phantom{-}2 & -1\\
        -1 & \phantom{-}0 & -2 & \phantom{-}2
    \end{bmatrix}.
\end{equation*}

To the best of the author's knowledge, \Cref{TBLstab} provides the first examples of finitely presented simple Kazhdan groups that are Frobenius stable and not Frobenius approximable.
It is worth remarking that there are already known examples of groups that have the latter two properties and have all but one qualification among "finitely presented", "simple", and "Kazhdan". 
In detail:
\begin{itemize}
    \item[(a)] The families constructed in~\cite[Theorem~1.1]{DGLT:stab} and~\cite[Theorem~1.3]{BLSW:stab} consist of suitable non-residually finite central extensions of finitely generated linear groups. The groups produced are \emph{finitely presented, Kazhdan, Frobenius stable} and \emph{not Frobenius approximable}. However, they cannot be simple by design (recall that, by Mal'cev theorem, every finitely generated linear group is residually finite). 
    \item[(b)] Fournier--Facio~\cite{ff:stab} proved that every simple group containing a non-Frobenius approximable subgroup is Frobenius stable and not Frobenius approximable~\cite[Lemma~2.1, Remark~2.2]{ff:stab}. Taking up what observed in~\cite[Remark~2.2]{ff:stab}, from a countable non-Frobenius approximable group~$A$ (possibly with solvable word problem, see~(b2)), thanks to~\cite[Lemma~2.1]{ff:stab} it is possible to construct:
    \begin{itemize}
        \item[(b1)] \emph{simple Kazhdan} groups containing~$A$ that are \emph{Frobenius stable} and \emph{not Frobenius approximable}. In fact, by~\cite[Theorem~4.1, Remark~4.5]{couff}, every countable group embeds in a simple Kazhdan group. However, that result explicitly also proves that the latter simple group is never finitely presented.
        \item[(b2)] \emph{finitely presented simple} groups that are \emph{Frobenius stable} and \emph{not Frobenius approximable}. It suffices to take~$A$ a non-Frobenius approximable finitely generated group with solvable word problem, e.g., finite-degree central extensions of suitable finitely generated linear groups as in~\cite[\S 5.2]{DGLT:stab}. OpenAI~\cite{OA:boohig} has recently claimed that the Boone--Higman conjecture is true. In this is correct, then $A$ embeds in a finitely presented simple group~$\Gamma$ which is Frobenius stable and not Frobenius approximable. 
        It is unclear to the author whether, with this construction, this group~$\Gamma$ can be additionally taken to be Kazhdan.
    \end{itemize}
\end{itemize}

\subsection{Finiteness properties and twin building lattices}\label{sus:introFP}
The study of finiteness properties has played a central role in the research of groups with root group data for at least three decades, motivated by the wealth of results of those properties for the nearby theory of $S$-arithmetic lattices. Recall that every twin building lattice~$\Lambda\cong G/Z(G)$ has the same finiteness properties of~$G$, as the centre~$Z(G)$ finite (see \Cref{prop:permFP}), so the focus naturally converges on general groups~$G$ with a locally finite root group datum~$(U_\alpha)_{\alpha\in \Phi}$ and their subgroup~$U_+$.

The majority of the results on finiteness properties are concentrated in criteria on finite generation and finite presentation for~$G$ and~$U_+$. 
Clearly, $G$ is finitely generated. 
If~$(W,S)$ is $2$-spherical and~$\qmin\geq 4$, say\footnote{In the original theorem, the hypothesis~"$\qmin\geq 4$" is replaced by the weaker condition~\hyperref[eq:co]{\co}. See \Cref{thm:OpConn}.}, the "Curtis--Tits theorem" proved by Abramenko--M\"uhlherr~\cite{AM97} gives an explicit finite presentation for~$G$. Under the same hypotheses, they proved that~$U_+$ is finitely generated~\cite[Théorème~M]{AM97}.  
It is worth mentioning that outside the $2$-spherical case or for small thickness, $G$ might not be even of type~$\FP_2(\Z)$ (see Abramenko--Gates~\cite[\S 5.1--5.3]{AbGates}) or at least not be finitely presented (see Bischof~\cite[Theorem~8.7]{bisch:unc}), and~$U_+$ might not be finitely generated (see Abramenko~\cite[\S I.6, Corollary~8]{ab:twin} and Bischof~\cite[Theorem~8.7]{bisch:unc}).
However, if~$(W,S)$ is $3$-spherical and~$\qmin\geq 16$, Devillers--M\"uhlherr~\cite[Corollary~1.2]{DevMu} proved that~$U_+$ is finitely presented.

Known positive results for higher finiteness properties are restricted to the case when~$|S|=k+2$ and~$(W,S)$ is $(k+1)$-spherical and classical. Under these assumptions and if~$\qmin\geq 2^{2|S|-1}$ and~$B_+\cap B_-$ is finite, Abramenko~\cite[\S I.6, Theorem~A and Corollary~11]{ab:twin} proved that~$U_+$ is of type~$\tF_k$ but not of type~$\FP_{k+1}(\Q)$, and~$G$ is of type~$\tF_k$. 
Abramenko's proof applies Brown's criterion (recalled in~\Cref{thm:BroCrit}) to Abels' filtration mentioned in \Cref{sus:introY-}. 
Thanks to \Cref{Y-conn}, we are able to generalise Abramenko's results as follows:
\begin{mainthm}\label{finProp1}
 Let~$\Lambda$ be a twin building lattice coming from a locally finite root group datum~$(U_\alpha)_{\alpha\in \Phi}$ for a group~$G$ of $(k+1)$-spherical Weyl group~$(W,S)$, for some~$k\geq 1$. Assume that~\eqref{Cond3} holds.
    Then:
    \begin{enumerate}
        \item $U_+$ is of type~$\FP_{k}(\Q)$ and, if~$(W,S)$ is classical and~$\qmin\geq 2^{2|S|-1}$, of type~$\tF_{k}$.
        Moreover, if~$|S|=k+2$ and~$(W,S)$ is non-spherical,~$U_+$ not of type~$\FP_{k+1}(\Q)$. 
        \item ~$G$ and~$\Lambda$ are of type~$\FP_{k}(\Q)$ and, if~$\qmin\geq 2^{2|S|-1}$ and~$(W,S)$ is classical, of type~$\tF_{k}$.
        \item If~$(W,S)$ is non-spherical, neither~$G$ nor~$\Lambda$ have type~$\FP_{2\cdot \sigma(W,S)-2}(\Q)$, where~$\sigma(W,S)$ is the largest cardinality of a spherical subset of~$S$. 
    \end{enumerate}
\end{mainthm}
Notice that, if~$|S|\geq 2$,~\eqref{Cond3} implies automatically that~$\qmin\geq 4$ (and for~$|S|\leq 1$, the relevant group~$G$ is finite).
See \Cref{thm:GFP} and \Cref{cor:TBLFP} for slightly more general statements.

\noindent The only real new contributions of \Cref{finProp1} are the positive parts of the statements. The negative ones were already known by work of Abramenko~\cite[\S I.6, Theorem~A]{ab:twin}, Gandini~\cite{gandini} and Abramenko--Gates~\cite[Proposition~3.2]{AbGates}.

In the Kac--Moody case, \Cref{finProp1} yields the following:
\begin{maincor}\label{U+FPk}
 Let~$A=[A_{st}]_{s,t\in S}$ be a generalised Cartan matrix with associated $(k+1)$-spherical Weyl group~$(W_A,S)$ for some~$k\geq 1$. Let also~$q>4(|S|-1)^2$ be a prime power.
Then:
\begin{enumerate}
    \item $\caG_A(\F_q)$ is of type~$\FP_{k}(\Q)$ and, if~$(W_A,S)$ is classical and~$q\geq 2^{2|S|-1}$, of type~$\tF_k$.
    \item If~$(W_A,S)$ is not spherical, then~$\caG_A(\F_q)$ is not of type~$\FP_{2\cdot\sigma(W_A,S)-2}(\Q)$.
    \item~$\caU_A(\F_q)$ is of type~$\FP_{k}(\Q)$ and, if~$|S|=k+2$ and $(W_A,S)$ is not spherical, it is not of type~$\FP_{k+1}(\Q)$.
\end{enumerate} 
\end{maincor}

\subsection{Acknowledgements}
The author thanks Pierre-Emmanuel Caprace for sharing several questions and points for reflection that shaped much of the material presented, and for valuable feedback on earlier versions of the draft. She also thanks Bharatram Rangarajan for sharing the question regarding the Frobenius stability/non-Frobenius approximability of Kac--Moody groups over~$\F_q$ that eventually led to \Cref{TBLstab}. She also thanks Alex Loué and Timothée Marquis for enlightening discussions that led to the idea of weakening the connectivity hypotheses in the results of~\cite{DJ:build}. The incentive to develop this project came after attending the reading group on higher properties~T that ran at UCLouvain in the a.a.~2025/2026. BM is grateful for the valuable contributions from all the participants.

The author's research is funded by F.R.S.-FNRS via a postdoctoral grant. She also acknowledges membership to the group "GNSAGA" within "Istituto Nazionale di Alta Matematica" (INdAM).

\subsection{AI statement}
The research presented in this paper and the writing of the manuscript have been carried out without any use of~AI at any stage, except for a sporadic use of the spellcheck "Writefull" embedded in Overleaf.

\section{Conventions and preliminaries}\label{s:prels}
\subsection{General conventions and notation}\label{sus:conv}
Throughout this paper, 
\begin{itemize}
    \item $Z(G)$ denotes the centre of a group~$G$.
    \item $\mathbb{D}^n:=\{x\in \R^n\colon |x|\leq 1\}$ denotes the $n$-dimensional disc in~$\R^n$ and~$\mathbb{S}^{n-1}=\partial\mathbb{D}^n$ the $(n-1)$-dimensional sphere in~$\R^n$.
    \item Group actions on sets are always on the left.
    \item Given a $G$-set~$X$, sometimes (when there is no ambiguity) we write~$X/G$ as a shortcut to denote an arbitrary set of representatives of~$X$ mod~$G$.
    \item Quotients of topological groups/modules are always endowed with the quotient topology. Direct products of topological groups/modules are always endowed with the product topology.
    \item An "lcsc~group" is a short way to say "locally compact second countable" group.
    \item When we refer to a "unitary representation"~$\rho\colon G\to \caU(V)$ over~$V$, we assume by default that~$V$ is a complex Hilbert space and~$\rho$ is strongly continuous.
    \item In a CW complex, we refer only to its closed cells and call them "cells".
    \item All cellular $G$-actions on CW complexes are \emph{admissible}, i.e., for every cell~$\sigma$, its pointwise stabiliser coincides with its setwise stabiliser (see \Cref{sus:PolyActions}).
    \item Given a cellular $G$-action on~$X$, we will always choose a cell orientation on~$X$ so that, for every cell~$\sigma$ and every~$g\in G$,~$g\cdot \sigma$ and~$\sigma$ have the same orientation (see \Cref{rem:Gorient}).
\end{itemize}

\subsection{Coxeter groups}\label{sus:cox}
Given a finite non-empty set~$S$, a \emph{Coxeter matrix} $[m_{st}]_{s,t\in S}$ is any symmetric matrix with diagonal entries equal to~$1$ and whose non-diagonal entries belong to~$\Z_{\geq 2}\cup\{+\infty\}$. The {\it Coxeter system}~$(W,S)$ associated to such a matrix~$[m_{st}]_{s,t\in S}$ is determined by defining~$W$ as the group generated by~$S$ subject to the relations $(st)^{m_{st}}=1$ whenever $m_{st}\neq\infty$. By~\cite[Proposition~5.3]{hum:ref},~$m_{st}$ coincides with the order of~$st$ in~$W$. The \emph{cosine matrix}~$\frC_{(W,S)}=[c_{st}]_{s,t\in S}$ of~$(W,S)$ is defined for all~$s,t\in S$ as
\begin{equation}\label{eq:cosmat}
    c_{st}:=-\cos\Big(\frac{\pi}{m_{st}}\Big).
\end{equation}
Throughout,~$\veps_{(W,S)}$ denotes the smallest eigenvalue of~$\frC_{(W,S)}$. The definition of~$(W,S)$ is usually extended to include~$S=\emptyset$ and, in that case,~$W$ is the trivial group.

For every~$J\subseteq S$, the subgroup~$W_J$ of~$W$ generated by~$J$ is called \emph{special subgroup of type~$J$}.
Notice that~$(W_J,J)$ is a Coxeter system for every $J\subseteq S$~\cite[Theorem~5.5(a)]{hum:ref}.  
 
A Coxeter system $(W,S)$ is said to be \emph{spherical} if $W$ is finite. A subset~$J\subseteq S$ is \textit{spherical} if~$(W_J,J)$ is spherical. Given~$k\in \Z_{\geq 1}$, a Coxeter system~$(W,S)$ is \emph{$k$-spherical} if all subsets~$J\subseteq S$ with~$|J|\leq k$ are spherical. Let
\begin{equation}\label{eq:sk(W,S)}
    \begin{split}
        \sigma(W,S) & :=\max\{|J|\colon J\subseteq S \text{ spherical}\};\\
        \kappa(W,S) & :=\max\{1\leq k\leq |S|\colon (W,S)\text{ is $k$-spherical}\}.
    \end{split}
\end{equation}

A subset~$J\subseteq S$ is said to be \emph{irreducible} if the graph~$\Gamma(W_J,J)$ with vertex set~$J$ and edge set~$\{\{s,t\}\mid s,t\in J,\,m_{st}\geq 3\}$ is connected. A Coxeter system~$(W,S)$ is called \emph{irreducible} if the set~$S$ is irreducible.
 A Coxeter system~$(W,S)$ is \emph{Euclidean} if all its irreducible subsets are spherical or affine, and~$(W,S)$ is \emph{classical} (or, following~\cite[\S 1.1]{OpVal1}, \emph{$|S|$-classical}) if, for every irreducible spherical subset~$J\subseteq S$, $(W_J,J)$ is of type~$A_n$, $B_n$, $C_n$ or~$D_n$ for some~$n\geq 1$. 

\begin{rem}\label{rem:Sph8}
    Note that~$\kappa(W,S)\leq 8$ whenever~$(W,S)$ is non-Euclidean. Indeed, every irreducible component~$(W_J,J)$ of~$(W,S)$ is at least~$\kappa(W,S)$-spherical. It is recorded in~\cite[Proposition~5.2]{lopn:betti} and in~\cite[Remark~6.7]{barn:cox} that all irreducible $9$-spherical Coxeter systems are spherical or affine, and this proves the claim.

    On the other hand, examples of non-Euclidean Coxeter systems~$(W,S)$ with~$\kappa(W,S)=8$ are recorded in~\cite[\S 5]{lopn:betti}.
\end{rem}

\subsection{BN-pairs and twin BN-pairs}\label{sus:BNpairs}
Let~$G$ be a group. A \emph{BN-pair~$(B,N)$ of type~$(W,S)$} for~$G$ consists of two subgroups~$B,N\leq G$ such that~$G=\langle B,N\rangle$, $T:=B\cap N\unlhd N$ and there is a group isomorphism~$f\colon W\to N/T$ and~$W$ admits a set of generators~$S$ such that the following hold:
\begin{enumerate}
    \item[(BN1)] $BsBwB\subseteq BwB\cup BswB$ for all~$s\in S$ and~$w\in W$; and,
    \item[(BN2)] $sBs^{-1}\not\subseteq B$ for all~$s\in S$.
\end{enumerate}
It is a standard fact that the pair~$(W,S)$ is a Coxeter system.

Every group~$G$ with a BN-pair~$(B,N)$ of type~$(W,S)$ admits the so-called \emph{Bruhat decomposition}~\cite[p.~320]{ab:build}, i.e., the following decomposition in $B$-double cosets:
\begin{equation*}
    G=\bigsqcup_{w\in W}BwB.
\end{equation*}
For every~$J\subseteq S$, the subset
\begin{equation*}
    P_J=BW_JB:=\bigsqcup_{w\in W_J}BwB
\end{equation*}
is a subgroup of~$G$~\cite[\S 5.2, p.~59]{ronan:book} and is called the \emph{standard parabolic subgroup of type~$J$}. Notice that~$P_\emptyset=B$ and~$P_S=G$. For simplicity, for every~$s\in S$ we write~$P_s$ in place of~$P_{\{s\}}$.
Given~$N_J:=f(W_J)$, it is easily proved that~$(B,N_J)$ is a BN-pair of type~$(W_J,J)$ for~$P_J$~\cite[Chapter~5, Exercise~6]{ronan:book}. Moreover, for all~$I,J\subseteq S$, recall that~\cite[Theorem~II.5.4(ii)]{ronan:book}:
\begin{equation}\label{eq:intersParab}
    P_I\cap P_J=P_{I\cap J}.
\end{equation}

The family of cardinals~$(q_s)_{s\in S}$ defined as~$q_s:=|BsB/B|=|P_s:B|-1$ for every~$s\in S$ is said to be the \emph{thickness vector} of the BN-pair~$(B,N)$. In particular, $(B,N)$ has \emph{finite thickness vector} if~$q_s$ is finite for every~$s\in S$.

\begin{rem}\label{rem:GBspher}
    Assume that~$(B,N)$ has finite thickness vector~$(q_s)_{s\in S}$. Then, for every spherical subset~$J\subseteq S$,
    $$|P_J:B|<\infty.$$ 
    Indeed, by~\cite[Proposition~4]{brw:build}, for every reduced expression~$w=s_1\cdots s_n$ in~$(W,S)$ we have
    \begin{equation*}
        |BwB/B|=|B:B\cap wBw^{-1}|=q_{s_1}\cdots q_{s_n}=:q_w.
    \end{equation*}
    By~\cite[Corollary~1.4.8(ii)]{bb:cox}, the quantity~$q_w$ is well-defined, i.e., it does not depend on the chosen reduced expression of~$w$ in~$(W,S)$.
    Hence,
    \begin{equation*}
        |{P_J}:B|=\sum_{w\in W_J}|BwB/B|=\sum_{w\in W_J}q_w<\infty.
    \end{equation*}
\end{rem}

\smallskip

Let~$G$ be a group with subgroups~$B_+, B_-,N\leq G$ such that~$B_+\cap N=B_-\cap N=:T\unlhd N$, and let~$W:=N/T$. Following~\cite[Definition~6.78]{ab:build}, the triple~$(B_+,B_-,N)$ is said to be a \emph{twin BN-pair of type~$(W,S)$} for~$G$ if there is a finite generating set~$S$ of~$W$ such that the following hold for all~$\veps\in\{\pm\}$:
\begin{enumerate}
    \item[(TBN0)] \label{TBN0} $(B_\veps, N)$ is a BN-pair of type~$(W,S)$ for~$G$;
    \item[(TBN1)] \label{TBN1} $B_\veps wB_\veps sB_{-\veps}=B_\veps wsB_{-\veps}$ for all~$w\in W$ and~$s\in S$ such that~$\ell(sw)<\ell(w)$;
    \item[(TBN2)] \label{TBN2} $B_+s\cap B_-=\emptyset$ for every~$s\in S$.
\end{enumerate}

In addition to the Bruhat decompositions, a group~$G$ with a twin BN-pair~$(B_+,B_-,N)$ of type~$(W,S)$ admits the so-called \emph{Birkhoff decompositions}~\cite[Proposition~6.81]{ab:build}, i.e., the following double-coset decompositions:
\begin{equation}\label{eq:birk}
    G=\bigsqcup_{w\in W}B_+wB_-=\bigsqcup_{w\in W}B_-wB_+.
\end{equation}
For every~$J\subseteq S$, we call~$P_J^\pm:=B_\pm W_JB_\pm$ the \emph{positive/negative standard parabolic subgroup of type~$J$}.
Once again, the family of cardinals~$(q_s^+,q_s^-)_{s\in S}$, where~$q_s^\pm=|B_\pm sB_\pm:B_\pm|$ for every~$s\in S$, is called the \emph{thickness vector} of the twin BN-pair~$(B_+,B_-,N)$. If~$q_s^+=q_s^-$ for every~$s\in S$, we write~$(q_s)_{s\in S}$ instead of~$(q_s^+,q_s^-)_{s\in S}$.

\section{CW complexes and groups acting on them}\label{s:Mkcompl}
\subsection{CW complexes}\label{sus:CWcompl}
We also assume familiarity with the basic definitions of CW complexes and cellular maps (see for instance~\cite[pp.~15-16 and p.~28]{GeoBook} for background). We take the occasion to recall some terminology and set a few conventions. Given a CW complex~$X$,
\begin{itemize}
    \item The term "cell(s)" in~$X$ always refers only to "closed cell(s)", i.e., "closures of open cells".
    \item For every~$k\in \Z_{\geq -1}$, a \emph{$k$-cell} is a $k$-dimensional cell. Denote by~$X(k)$ the collection of all $k$-cells of~$X$. Conventionally, $X(-1)$ contains a unique cell, the \emph{empty cell}, which corresponds to the empty subset in~$X$. Set~$\Ce(X):=\bigcup_{k\geq -1}X(k)$.
    \item Denote by~$X^{(k)}=\bigcup_{\ell\leq k}\bigcup_{\sigma\in X(\ell)}\sigma$ the $k$-th \emph{skeleton} of~$X$.
    \item Let~$\dim(X):=\min\big\{\{k\geq -1\mid X=X^{(k)}\}\cup\{+\infty\}\big\}$ be the \emph{dimension} of~$X$.
\end{itemize}

A CW complex $X$ is said to be:
\begin{itemize}
    \item \emph{pure $n$-dimensional} for an integer~$n$, if~$\dim(X)=n$ and every cell in~$X$ is contained in some $n$-cell;
    \item \emph{countable} (resp.~\emph{finite}) if the set~$\Ce(X)$ is countable (resp.~finite);
    \item \emph{locally finite} if every cell of~$X$ is contained in finitely many cells;
    \item \emph{connected} if the $1$-skeleton of~$X$ is a connected graph.
\end{itemize}
In a CW complex~$X$, its $0$-cells and its $1$-cells are called \emph{vertices} and~\emph{edges}, respectively. Moreover, if~$X$ is finite dimensional, borrowing a term from building theory, we call its top-dimensional cells \emph{chambers}.

\smallskip

In this paper, by default, the term "simplicial complex" refers to \emph{abstract simplicial complex}, that is, any collection~$X$ of subsets of a given set~$V$ (called the \emph{vertex set} of~$X$) closed under taking subsets of its elements and that contains each singleton~$\{v\}$, for all~$v\in V$. As usual, we will simply write~$v$ to denote the $0$-simplex~$\{v\}$.
Every abstract simplicial complex admits a canonical structure of CW complex~\cite[\S 5.2]{GeoBook}. Thus, we will use the previous terminology and notation for CW complexes also for simplicial complexes (e.g., $X(k)$ denotes the set of $k$-simplices of~$X$, $\dim(X)$ is the dimension of~$X$, etc...).
The only exception is that we simply write~$X$ in place of~$\Ce(X)$. For an $n$-dimensional simplicial complex~$X$, where~$1\leq n<\infty$, we also have the following refined notion of connectivity:~$X$ is said to be \emph{gallery connected} if, for all chambers~$\gamma,\eta\in X(n)$, there is a finite sequence in~$X(n)$:
$$\gamma=\gamma_0\sim \gamma_1 \sim \dots \sim \gamma_m=\eta \qquad  m\geq 0$$ 
where $\gamma_{i-1}\sim \gamma_i$ means that~$\gamma_{i-1}\cap \gamma_i$ is a codimension-$1$ face of both~$\gamma_{i-1}$ and~$\gamma_i$.

In this paper, the term \emph{graph} is always used to denote a simplicial complex of dimension~$\leq 1$.
Given a graph~$\Gamma$, a subgraph~$\Lambda$ is said to be \emph{induced} if, whenever~$\{s,t\}\in \Gamma(1)$ and~$s,t\in \Lambda(0)$, the edge~$\{s,t\}$ belongs to~$\Lambda(1)$.

\smallskip

For every simplex~$\sigma$ in a simplicial complex~$X$, the \emph{link} of~$\sigma$ in~$X$ is the subcomplex
\begin{equation}\label{eq:link}
    X_\sigma=\Lk_X(\sigma):=\{\tau\in X\mid \sigma\cap \tau=\emptyset,\,\sigma\cup\tau\in X\}.
\end{equation}
Notice that~$\dim(X_\sigma)=\dim(X)-\dim(\sigma)-1$. Moreover, if~$X$ is pure $n$-dimensional then~$X_\sigma$ is pure $(n-\dim(\sigma)-1)$-dimensional.

\smallskip

 A \emph{cell orientation} on a CW complex~$X$ is the choice of an orientation for each cell of~$X$.
 We refer to~\cite[\S 2.5]{GeoBook} for details. If in particular~$X$ is a simplicial complex, every total ordering~$\leq$ on its vertex set~$X(0)$ canonically induces an orientation on each simplex~$\sigma=\{v_0<\dots<v_n\}$ as described in~\cite[pp.~139-140]{GeoBook}.

\smallskip

In the next two subsections, we will focus on two particular subfamilies of simplicial complexes that will be widely used in this paper.

\subsection{Partite simplicial complexes}\label{sus:partCompl}
Let~$(X_s)_{s\in S}$ be a family of non-empty sets indexed by a set with~$n+1$ elements, $n\in \Z_{\geq 0}$. An $(X_s)_{s\in S}$-\emph{partite} (or, more generically, an \emph{$(n+1)$-partite}) \emph{simplicial complex} is a pure $n$-dimensional simplicial complex~$X$ together with a vertex partition~$X(0)=\bigsqcup_{s\in S}X_s$ with $n+1$~terms (called \emph{sides}) such that every~$\sigma\in X(n)$ contains exactly one vertex per side. Each such partition has an associated \emph{type function}~$\tau\colon X\to 2^S$, which is a map that assigns to each simplex~$\sigma$ of~$X$ the set~$\{s\in S\mid \sigma\cap X_s\neq \emptyset\}$. 

If~$n=1$, we prefer the more standard terminologies \emph{$(A,B)$-partite graph} (where~$A$ and $B$ are non-empty sets) or simply \emph{bipartite graph}.
For~$a,b\in \Z_{\geq 1}$, an~$(A,B)$-partite graph~$\Gamma$ is said to be \emph{$(a,b)$-regular} (or simply \emph{biregular}) if, for all~$v\in A$ and~$w\in B$, 
\begin{equation*}
    \#\{e\in \Gamma(1)\mid e\ni v\}=a\quad \text{and}\quad \#\{e\in \Gamma(1)\mid e\ni w\}=b.
\end{equation*}

\smallskip

In this paper, we mainly deal with a particular instance of partite complexes: the \emph{coset complex}~$\caC\caC(G, (H_i)_{i\in I})$ of a group~$G$ and a finite non-empty family~$(H_i)_{i\in I}$ of subgroups. We briefly recall its definition, following~\cite[\S 2.2]{GdP-VP:KMS} and the references therein.
The simplicial complex~$\caC\caC(G, (H_i)_{i\in I})$ has vertex set
\begin{equation}\label{eq:CComplx}
    \bigsqcup_{i\in I}G/H_i
\end{equation}
and its $k$-simplices, where~$k\geq 1$, are all those sets of the form:
\begin{equation*}
    \sigma(g,J):=\big\{gH_{j}\mid j\in J\big\},
\end{equation*}
for some~$g\in G$ and~$J\subseteq I$ with~$|J|=k+1$.
Putting~$n=|I|$, notice that~$\caC\caC(G, (H_i)_{i\in I})$ is an $(n+1)$-partite simplicial complex with vertex partition as in~\eqref{eq:CComplx} and type function
\begin{equation}\label{eq:typeCC}
    \type\big(\sigma(g,J)\big):=J\qquad\forall\,g\in G,\,J\subseteq I.
\end{equation}
Throughout, for every non-empty subset~$J\subseteq S$ we set
\begin{equation*}
    H_J:=\bigcap_{i\in J}H_i.
\end{equation*}

\begin{lem}[\protect{\cite[Lemmas~3.3-3.4]{HarSap}}]\label{lem:CCprop}
    Given~$X=\caC\caC(G, (H_i)_{i\in I})$ as before, the following hold:
    \begin{enumerate}
        \item \label{lem:CCprop1} $X$ is gallery-connected if and only if~$G=\langle H_i\mid i\in I\rangle$; 
        \item \label{lem:CCprop2} for every simplex~$\sigma$ of~$X$ with~$\type(\sigma)=J\neq \emptyset$, the link~$X_\sigma$ is isomorphic to the coset complex~$\caC\caC(H_J,(H_{J\cup\{i\}})_{i\in I\setminus J})$.
    \end{enumerate}
\end{lem}

\subsection{(Twin) buildings induced by (twin) BN-pairs}\label{sus:build}
Here we follow~\cite[\S 1.1]{ronan:book}. A {\em chamber system} over a set~$S$ is a set $\caC$, whose elements are called {\em chambers}, together with a family of equivalence relations $(\sim_s)_{s\in S}$ on~$\caC$ indexed by~$S$. Two chambers~$c$ and $d$ are said to be {\em $s$-adjacent} if~$c\sim_s d$. A {\em gallery} from~$c$ to~$d$, where~$c,d\in \caC$, is any chain of chambers of the form~$c=c_0\sim_{s_1}c_1\sim_{s_2}\dots \sim_{s_n}c_n=d$. The sequence~$(s_1,\dots, s_n)$ is called the \emph{type} of the gallery. A chamber system is said to be \emph{gallery-connected} if every two chambers can be connected by a gallery. 

\begin{ex}[\protect{\cite[\S 1, Example~1]{ronan:book}, \cite[Remark~2.14]{GdP-VP:KMS}}]\label{ex:GrpChambSys}
  Let~$G$ be a group and~$(P_i)_{i\in I}$ be a finite non-empty family of subgroups with a subgroup~$B\leq \bigcap_{i\in I}P_i$. The coset set~$G/B$ is a chamber system over~$I$ with $i$-adjacency relation ($i\in I$) defined as follows:
   \begin{equation*}
       hB\sim_i gB\,\Longleftrightarrow\, g^{-1}h\in P_i.
   \end{equation*}
   The associated chamber system is denoted by~$\caC(G, (P_i)_{i\in I},B)$. It is immediate to check that~$\caC(G, (P_i)_{i\in I},B)$ is connected if and only if~$G=\langle P_i\mid i\in I\rangle$.

   Assume additionally that there are subgroups~$(H_i)_{i\in I}$ of~$G$ such that, given~$H_J:=\bigcap_{j\in J}H_j$ for every~$J\subseteq I$, the following holds for all~$i\in I$ and~$J\subseteq I$:
   \begin{equation}\label{eq:CCcond}
       P_i=H_{I\setminus\{i\}},\quad B=H_I,\quad\text{and}\quad H_{I\setminus J}=\langle H_{I\setminus\{j\}}\mid j\in J\rangle.
   \end{equation}
   A typical example of the previous situation is when~$B$ is the "$B$"-group of a BN-pair of type~$(W,S)$ for~$G$, $I=S$ and, for every~$s\in S$, $P_s$ and~$H_s$ are the standard parabolic subgroups~$P_s$ and~$P_{S\setminus\{s\}}$, respectively.
   
  Returning to the general framework set by~\eqref{eq:CCcond}, let~$X=\caC\caC(G, (H_i)_{i\in I})$ and fix~$n=|I|-1$. The set of top-dimensional simplices of~$X$:
   \begin{equation*}
       X(n):=\Big\{\sigma(g,I):=\{gH_i\mid i\in I\}\,\Big\vert\, g\in G\Big\}
   \end{equation*}
   is a chamber system over~$I$ with $i$-adjacency ($i\in I$):
   \begin{equation*}
       \sigma(h,I)\sim_i\sigma(g,I)\,\Longleftrightarrow\,g^{-1}h\in H_{I\setminus \{i\}}.
   \end{equation*}
   Moreover, the assignment~$gB\mapsto \sigma(g,I)$, $g\in G$, determines a bijective map~$\caC(G,(P_i)_{i\in I},B)\stackrel{\simeq}{\longrightarrow} X(n)$ preserving the $i$-adjacencies for all~$i\in I$. In this case, we call~$X$ the \emph{coset complex associated to~$\caC(G,(P_i)_{i\in I},B)$}.
\end{ex}

Let~$(B,N)$ be a BN-pair of type~$(W,S)$ for a group~$G$. For every~$J\subseteq S$, recall that~$P_J:=BW_JB$ denotes the standard parabolic subgroup of type~$J$, and that~$P_s:=P_{\{s\}}$ for every~$s\in S$. The chamber system over~$S$
\begin{equation}\label{eq:CBN}
    \caC:=\caC(G, (P_s)_{s\in S},B)
\end{equation}
as in \Cref{ex:GrpChambSys} is, in fact, a building of type~$(W,S)$ and Weyl distance~$\delta\colon \caC\times \caC\to W$ given for all~$w\in W$ as follows:
\begin{equation}\label{eq:twinn}
    \delta(hB,gB)=w\,\Longleftrightarrow\, g^{-1}h\in BwB,
\end{equation}
see~\cite[Theorem~5.3]{ronan:book}. For every~$s\in S$, the $s$-adjacency relation in~$\caC$ is
\begin{equation*}
    hB\sim_s gB\,\Longleftrightarrow\,g^{-1}h\in P_s.
\end{equation*}
We refer to~\cite[Ch.~3, \S 1]{ronan:book} or~\cite[\S 5.1.1]{ab:build} for background on buildings as chamber systems. Here we only recall a few details that we need.
For all~$J\subseteq S$ and~$c\in \caC$, the \emph{residue of type~$J$ centred at~$c$} is the chamber (sub)system over~$J$ defined by
\begin{equation*}
    \Res_J(c):=\{d\in \caC\mid \delta(c,d)\in W_J\}
\end{equation*}
and, for every~$s\in J$, with restricted $s$-adjacency relations from~$\caC$. For every~$J\subseteq S$ and~$g\in G$, notice that
\begin{equation*}
    \Res_J(gB)=g\cdot\Res_J(1B)\qquad\text{and}\qquad \Res_J(1B)=P_J/B.
\end{equation*}

As already mentioned before, in this setting~\eqref{eq:CCcond} is satisfied by putting~$H_s:=P_{S\setminus \{s\}}$ for every~$s\in S$. The coset complex associated to~$\caC$, namely
\begin{equation}\label{eq:CCBN}
    \Delta(G,B)=\caC\caC(G, (P_{S\setminus\{s\}})_{s\in S}),
\end{equation}
is usually called the \emph{simplicial building} associated to the BN-pair~$(B,N)$ of~$G$. 
By \Cref{lem:CCprop}\eqref{lem:CCprop2}, for every simplex~$\sigma=\{gP_{S\setminus\{s\}}\mid s\in J\}$ of type~$J$,  its link in~$\Delta(G,B)$ is isomorphic to the simplicial building associated to the BN-pair for~$P_{S\setminus J}$ induced by~$(B,N)$; in symbols
\begin{equation*}
   \Delta(G,B)_\sigma\cong\caC\caC\big(P_{S\setminus J}, (P_{(S\setminus J)\setminus\{s\}})_{s\in S\setminus J}\big)=\Delta(P_{S\setminus J},B)
\end{equation*}

\medskip

Let now~$(B_+,B_-,N)$ be a twin BN-pair of type~$(W,S)$ for~$G$. 
Recall that~$(B_\pm,N)$ is a BN-pair of~$G$ of type~$(W,S)$ on its own, and thus it comes with an associated building~$\caC_\pm:=\caC(G, (P_s^\pm)_{s\in S},B_\pm)$ of type~$(W,S)$. 
In this case, $\caC_+$ and~$\caC_-$ are inter-related through the so-called \emph{twinning} (or \emph{codistance}), which is the function
$\delta\colon \caC_+\times \caC_-\cup \caC_-\times \caC_+\longrightarrow W$ defined for all~$w\in W$ as follows:
\begin{equation}\label{eq:defndelta}
    \delta(hB_\pm,gB_\mp)=w\,\Longleftrightarrow\, g^{-1}h\in B_\mp wB_\pm.
\end{equation}
The triple~$(\caC_+,\caC_-,\delta)$ is indeed an instance of a \emph{twin building of type~$(W,S)$} (see~\cite[Definition~5.133]{ab:build} for the general definition and~\cite[p.~332]{ab:build} for a proof).

\subsection{Weighted adjacency matrices and cosine matrices}\label{sus:matrices}
Let~$\Gamma$ be a finite graph with vertex set~$V$ and non-empty edge set~$E$. 
The \emph{adjacency matrix}~$A_\Gamma=[G_{\Gamma}(v,w)]_{v,w\in V}$ of~$\Gamma$ is the real square matrix defined as:
\begin{equation}\label{eq:adjGamma}
    A_\Gamma(v,w):=\left\{
        \begin{array}{cl}
          1  & \text{if }\{v,w\}\in E,\\
          0 & \text{otherwise.}
        \end{array}
        \right.
    \end{equation}
Being real and symmetric, $A_\Gamma$ self-adjoint. Thus we can list its eigenvalues in some non-increasing order as
\begin{equation}\label{eq:mui}
    \mu_1(\Gamma)\geq \mu_2(\Gamma)\geq \dots \geq \mu_N(\Gamma),\qquad N=|V|.
\end{equation}

If~$\Gamma$ has no leaves (i.e., $\Gamma$ is pure $1$-dimensional), the graph~$\Gamma$ also admits an associated \emph{normalised weighted adjacency matrix}, i.e., the matrix~$G_\Gamma=[G_{\Gamma}(v,w)]_{v,w\in V}$ defined as follows:
    \begin{equation}\label{eq:GGamma}
        G_\Gamma(v,w):=\left\{
        \begin{array}{cl}
           \frac{1}{\#\{e\in E\colon v\in E\}}  & \text{if }\{v,w\}\in E,\\
           0 & \text{otherwise.}
        \end{array}
        \right.
    \end{equation}
  Once again, being real and symmetric, the square matrix~$G_\Gamma$ is self-adjoint and we can list its eigenvalues in some non-increasing order as
    \begin{equation}\label{eq:lambdai}
        \lambda_1(\Gamma)\geq \lambda_2(\Gamma)\geq \dots\geq \lambda_N(\Gamma),\qquad N=|V|.
    \end{equation}

The following proposition is a generalisation of~\cite[Lemma~5.5]{HarSap} to arbitrary bi-regular bipartite graphs.
\begin{prop}\label{prop:specSubgraph}
    Let~$\Gamma$ be a finite $(a,b)$-regular $(A,B)$-partite, and let~$\Lambda\subseteq \Gamma$ be an induced subgraph which is an $(a',b')$-regular $(A',B')$-partite graph, where~$A':=A\cap \Lambda(0)$ and~$B':=B\cap \Lambda(0)$.
    Then
    \begin{equation*}
        \lambda_2(\Lambda)\leq \sqrt{\frac{ab}{a'b'}}\cdot \lambda_2(\Gamma).
    \end{equation*}
\end{prop}
\begin{proof}
For simplicity, let~$n:=|\Gamma(0)|$ and~$m:=|\Lambda(0)|$. Since~$\Lambda$ is an induced subgraph of~$\Gamma$, $A_\Lambda(v,w)=A_\Gamma(v,w)$ for all~$v,w\in \Lambda(0)$.
Moreover, as both~$\Gamma$ and~$\Lambda$ bipartite graphs, up to reordering the entries of~$A_\Gamma$ we have the following block decomposition:
   \begin{equation}\label{eq:GmatComp}
       A_\Gamma=\begin{bmatrix}
           0 & U\\
           W & 0
       \end{bmatrix}\quad\text{and}\quad A_\Lambda=\begin{bmatrix}
           0 & U'\\
           W' & 0
       \end{bmatrix},
   \end{equation}
   where~$U=[A_\Gamma(v,w)]_{v\in A,w\in B}$, $W=[A_\Gamma(v,w)]_{v\in B,w\in A}$, $U'=[A_\Gamma(v,w)]_{v\in A',w\in B'}$ and~$W'=[A_\Gamma(v,w)]_{v\in B',w\in A'}$. Similarly, we have
\begin{equation}\label{eq:AmatComp}
       G_\Gamma=\begin{bmatrix}
           0 & \frac{1}{a}U\\
           \frac{1}{b}W & 0
       \end{bmatrix}\quad\text{and}\quad G_\Lambda=\begin{bmatrix}
           0 & \frac{1}{a'}U'\\
           \frac{1}{b'}W' & 0
       \end{bmatrix}.
   \end{equation}
   
   For this paragraph, we now argue as in the proof of~\cite[Lemma~5.5]{HarSap}: 
   by the Courant--Fischer theorem (recorded for instance in~\cite[Theorem~A.1]{HarSap}) we have
   \begin{equation*}
       \mu_2(\Lambda)=\max_{v\in \R^{m}\colon v\perp \bbone_{m},\,v\neq 0}\frac{v^tA_\Lambda v}{v^tv}\quad\text{and}\quad\mu_2(\Gamma)=\max_{v\in \R^{n}\colon v\perp \bbone_{n},\,v\neq 0}\frac{v^tA_\Gamma v}{v^tv},
   \end{equation*}
   where~$\bbone_m\in \R^m$ and~$\bbone_n\in \R^n$ are the vectors with all entries equal to~$1$. Every~$v\in \R^{m}$ orthogonal to~$\bbone_{m}$ can be padded with zero entries to obtain a vector~$\overline{v}\in \R^{n}$. Observe that~$\overline{v}\perp \bbone_n$,~$v^tv=\overline{v}^t\overline{v}$ and~$v^tA_\Lambda v=\overline{v}^tA_{\Gamma} \overline{v}$. Hence,
    \begin{equation}\label{eq:lambda2compar}
        \mu_2(\Lambda)=\max_{v\in \R^{m}\colon v\perp \bbone_{m},\,v\neq 0}\frac{v^tA_\Lambda v}{v^tv}\leq\max_{v\in \R^{n}\colon v\perp \bbone_{n},\,v\neq 0}\frac{v^tA_\Gamma v}{v^tv}=\mu_2(\Gamma).
    \end{equation}

    To conclude the argument, we claim that
    \begin{equation}\label{eq:lamba2compar2}
\lambda_2(\Gamma)=\frac{1}{\sqrt{ab}}\mu_2(\Gamma)\quad\text{and}\quad \lambda_2(\Lambda)=\frac{1}{\sqrt{a'b'}} \mu_2(\Lambda).
    \end{equation}
    We prove~\eqref{eq:lamba2compar2} for~$\Gamma$, as for~$\Lambda$ one can proceed analogously.
  In fact, we will prove the more general equality of spectra:
   \begin{equation}\label{eq:specCompar}
       \Spec(G_\Gamma)=\Bigg\{\frac{1}{\sqrt{ab}} \lambda\,\Bigg\vert\,\lambda\in \Spec(A_\Gamma)\Bigg\}.
   \end{equation}
   Since the map~$\R\to \R$, $x\mapsto \frac{1}{\sqrt{ab}}x$, is strictly increasing, \eqref{eq:specCompar} implies~\eqref{eq:lambda2compar}.
   
   Firstly, by~\eqref{eq:GmatComp} and~\eqref{eq:AmatComp} we have
   \begin{equation*}
       G_\Gamma=\begin{bmatrix}
           \frac{1}{a}I_{|A|} & 0 \\
           0 & \frac{1}{b}I_{|B|}
       \end{bmatrix}\cdot A_\Gamma,
   \end{equation*}
   where~$I_{|A|}$ and~$I_{|B|}$ are the real identity matrices of dimension~$|A|$ and~$|B|$, respectively.
   Hence~$\det(G_\Gamma)=\frac{1}{ab}\det(A_\Gamma)$ and
   \begin{equation}\label{eq:specAG1}
       0\in \Spec(A_\Gamma)\,\Longleftrightarrow\,0\in \Spec(A_\Gamma).
   \end{equation}
   Let now~$\lambda\in \R$ satisfy~$\lambda\neq 0$. We first recall \emph{Schur's formula}: given a square real matrix with block decomposition
   $$M=\begin{bmatrix}
       A & B\\
       C & D
   \end{bmatrix}$$ 
   made up of square blocks~$A, D$ with~$A$ invertible, we have
   \begin{equation}\label{eq:schur}
       \det(M)=\det(A)\cdot \det(D-CA^{-1}B).
   \end{equation}
  Hence, by Schur's formula, we obtain
   \begin{equation*}
   \begin{split}
       \det(A_\Gamma-\lambda I_n) & =\det\left(\begin{bmatrix}
           -\lambda I_{|A|} & U\\
           V & -\lambda I_{|B|}
       \end{bmatrix}\right)\\
       & =\det(-\lambda I_{|A|})\cdot \det\big(-\lambda I_{|B|}-\lambda^{-1}WU\big)\\
       & =(-1)^{|A|}\lambda^{|A|-|B|}\det(WU-\lambda^2I_{|B|}).
    \end{split}
   \end{equation*}
  Proceeding analogously for~$G_\Gamma$, we obtain
  \begin{equation*}
      \det(G_\Gamma-\lambda I_n)=\frac{(-1)^{|A|}}{(ab)^{|B|}}\lambda^{|A|-|B|}\det(WU-ab\lambda^{2}I_{|B|}).
  \end{equation*}
  Since~$a,b>0$ and~$\Spec(A_\Gamma)$ is symmetric around~$0$, we deduce that
  \begin{equation}\label{eq:specAG2}
      \lambda\in \Spec(G_\Gamma)\setminus\{0\}\Longleftrightarrow\pm \sqrt{ab}\cdot |\lambda|\in \Spec(A_\Gamma)\setminus\{0\}\Longleftrightarrow\sqrt{ab}\cdot \lambda\in \Spec(A_\Gamma)\setminus\{0\}.
  \end{equation}
  Combining~\eqref{eq:specAG1} and~\eqref{eq:specAG2} we obtain~\eqref{eq:specCompar}. This concludes the argument.
\end{proof}
For future reference, we give the following explicit computation.
\begin{ex}\label{ex:Bip}
    Let~$\Gamma$ be a finite graph with a vertex bipartition~$\Gamma(0)=A\sqcup B$ such that~$|A|=|B|=q+1$ for some integer~$q\geq 1$, and that induces the following square-block decomposition of~$A_\Gamma$:
    \begin{equation*}
        A_\Gamma=\begin{bmatrix}
            0 & \bbone-I\\
        \bbone-I & 0
        \end{bmatrix},
    \end{equation*}
where~$\bbone$ is the $(q+1)$-dimensional matrix with all entries equal to~$1$. In particular, $\Gamma$ is~$q$-regular.
We claim that the spectrum~$\Spec(G_\Gamma)$ is
    \begin{equation*}
        \Spec(G_\Gamma)=\textstyle{\frac{1}{q}}\cdot\Spec(A_\Gamma)=\Big\{\pm \textstyle{\frac{1}{q}},\pm 1\Big\}.
    \end{equation*}
The fist equality follows from~\eqref{eq:specCompar}. It remains to determine~$\Spec(A_\Gamma)$.
    Notice that~$0\not\in \Spec(A_\Gamma)$ because~$A_\Gamma$ is invertible. By Schur's formula~\eqref{eq:schur}, for every~$\lambda\in \R_{\neq 0}$ we get
\begin{equation*}
\begin{split}
  \det(A_\Gamma-\lambda I)=0 &\Longleftrightarrow\,\det(-\lambda I+\lambda^{-1}(\bbone-I)^2)=0\,\Longleftrightarrow\,\det((\bbone-I)^2-\lambda^2I)=0 \\
  & \Longleftrightarrow\, \det(\bbone-(1+\lambda)I)=0\text{ or }\det(\bbone-(1-\lambda)I)=0\\
    & \Longleftrightarrow\, 1+\lambda\in \Spec(\bbone)\text{ or }1-\lambda\in \Spec(\bbone).
\end{split}
\end{equation*}    
Notice that~$\Spec(\bbone)=\{0,q+1\}$: indeed, since~$\bbone$ has rank~$1$, $0$ is an eigenvalue of~$\bbone$ of geometric multiplicity~$q$. Moreover, it is easily checked that~$q+1$ is an eigenvalue of~$\bbone$. Since~$\bbone$ has~$q+1$ eigenvalues counted with their algebraic multiplicity, this yields the claim on~$\Spec(\bbone)$. By what observed before, we conclude that~$\Spec(A_\Gamma)=\{\pm 1, \pm q\}$.
\end{ex}

Another matrix that will be useful throughout is the following.
\begin{defn}[\protect{\cite[Definition~1.3]{gro:curv}}]\label{defn:CosMatX}
    For~$n\geq 1$, let~$X$ be a $(n+1)$-partite simplicial complex with type function~$\tau\colon X\to 2^{S}$, with finite $1$-dimensional links and such that the group~$\Aut_\tau(X)$ of all type-preserving automorphisms~$X\to X$ is transitive on~$X(n)$.

    The \emph{cosine matrix} of~$X$ is the real $(n+1)$-dimensional matrix~$A(X)=[A(X)_{st}]_{s,t\in S}$ defined as:
    \begin{equation*}
        A(X)_{st}:=\left\{
        \begin{array}{cl}
            1 & \text{if }s=t, \\
            \lambda_{st} & \text{otherwise}, 
        \end{array}\right.
    \end{equation*}
    where~$\lambda_{st}$ is the second largest eigenvalue of the link in~$X$ of \emph{any}\footnote{Since~$\Aut_\tau(X)$ acts type-preservingly on~$X$ and transitively on~$X(n)$, the choice of the simplex of a given type does not change the isomorphism class of the link~\cite[Remark~1.4]{gro:curv}.} simplex of~$X$ of type~$S\setminus\{s,t\}$.
\end{defn}

A typical exemplification of \Cref{defn:CosMatX} (which also motivates the name) is the following.
\begin{ex}\label{ex:cosMatBuil}
    Let~$X=\Delta(G,B)$ be a simplicial building coming from a BN-pair~$(B,N)$ of $2$-spherical type~$(W,S)$ for a group~$G$. Denote by~$M=[m_{st}]_{s,t\in S}$ the Coxeter matrix of~$(W,S)$. Assume the given BN-pair has finite thickness vector~$(q_s)_{s\in S}$, and let
    \begin{equation*}
        \qmin:=\min\{q_s\mid s\in S\}\quad\text{and}\quad \qmax:=\max\{q_s\mid s\in S\}.
    \end{equation*}

    Recall that~$X$ is a $(n+1)$-partite simplicial complex, where~$n=|S|-1\geq 1$ and with type function~$\tau\colon X\to 2^S$ given in~\eqref{eq:typeCC}. By hypothesis, all $1$-dimensional links (i.e., the rank-$2$ residues) of~$X$ are finite.

    Assume now~$\qmin\geq 2$. Since~$(W,S)$ is $2$-spherical, a theorem by Feit--Higman~\cite{FeHig} implies that~$m_{st}\in\{2,3,4,6,8\}$ for all~$s\neq t$. In this case we have an explicit expression of~$\lambda_{st}$ in terms of~$m_{st}$, $q_s$ and~$q_t$, as recorded in~\cite[\S 3.1]{BallSwiat} (see also~\cite[Proposition~7.10]{garland} for the case~$m_{st}\leq 6$).
    Based on this, for all~$s\neq t$ in~$S$ we deduce that
    \begin{equation}\label{eq:AXineq}
        -A(X)_{st}=\lambda_{st}\leq 2\cdot\cos\Big(\frac{\pi}{m_{st}}\Big)\cdot\sqrt{\frac{\max\{q_s,q_t\}}{(q_s+1)(q_t+1)}}\leq 2\cdot\cos\Big(\frac{\pi}{m_{st}}\Big)\cdot\frac{\sqrt{\qmax}}{\qmin+1}.
    \end{equation}
   The chain of inequalities in~\eqref{eq:AXineq} is inspired from the proof of~\cite[Lemma~4.5]{gro:curv}. (However, to the best of the author's understanding, their proof seems to contain a typo: the authors claims that~$\lambda_{st}\leq 2\cdot\cos(\frac{\pi}{m_{st}})\cdot\frac{\sqrt{\qmin}}{\qmin+1}$).
    
    In a compact way, we have the following entry-wise inequality:
    \begin{equation}\label{eq:AXestim}
        A(X)\geq \frac{2\sqrt{\qmax}}{\qmin+1}\frC_{(W,S)}+\Bigg(1-\frac{2\sqrt{\qmax}}{\qmin+1}\Bigg)I,
    \end{equation}
    where~$\frC_{(W,S)}$ is the cosine matrix of~$(W,S)$. 
    Remarkably, both matrices in~\eqref{eq:AXestim} have diagonal entries equal to~$1$ and all off-diagonal entries in~$\R_{\leq 0}$. Thus, by~\cite[Proposition~2.8]{gro:curv}, the smallest eigenvalue of~$A(X)$ is at least the smallest eigenvalue of the matrix on the right-hand side of~\eqref{eq:AXestim}. Recall that~$\veps_{(W,S)}$ denotes the smallest eigenvalue of~$\frC_{(W,S)}$. Therefore, if
    \begin{equation}\label{eq:ineqPosDef}
        \frac{\qmin+1}{2\sqrt{\qmax}}>1-\veps_{(W,S)}
    \end{equation}
    from the previous observations we have that~$A(X)$ is positive definite.
\end{ex}

\subsection{Groups actions on CW complexes}\label{sus:PolyActions}
Let~$G$ be a group and~$X$ be a CW complex. A \emph{cellular $G$-action~$G\curvearrowright X$ on~$X$} is any $G$-action on~$X$ by cellular homeomorphisms that is \emph{admissible}, i.e., for every cell~$\sigma$, its setwise stabiliser coincides with its pointwise stabiliser. We denote by~$G_\sigma$ the pointwise stabiliser of~$\sigma$ in~$G$. Note that every cellular $G$-action on~$X$ induces a $G$-action on~$X(k)$ for every~$k\geq -1$. If~$X$ has a given cell orientation, a cellular $G$-action on~$X$ is said to be \emph{orientation-preserving} if, for all~$g\in G$ and~$\sigma\in \Ce(X)$, the cell~$g\cdot \sigma$ has the same orientation as~$\sigma$.

\begin{rem}\label{rem:Gorient}
    Given a $G$-action~$G\curvearrowright X$ on~$X$, there is always a cell orientation on~$X$ that makes~$G\curvearrowright X$ orientation-preserving.

    For instance, for every~$k\geq 0$, let~$R_k$ be a set of representatives for the $G$-orbits on~$X(k)$ and fix an orientation~$\veps_\sigma\in\{\pm 1\}$ for each~$\sigma\in R_k$. Then, for every~$g\in G$ and~$\sigma\in R_k$, choose$\veps_{g\cdot \sigma}:=\veps_\sigma$ as the orientation on~$g\cdot\sigma$. 
\end{rem}

If~$X$ is a simplicial complex, a \emph{simplicial $G$-action on~$X$} is any $G$-action on~$X$ by simplicial isomorphisms. Endowing~$X$ with its canonical CW structure (see \Cref{sus:CWcompl}), every simplicial action can be seen as a particular cellular action. 

\smallskip

A cellular action~$G\curvearrowright X$ is said to be:
\begin{itemize}
    \item \emph{discrete} if all pointwise stabilisers of cells are open;
    \item \emph{proper} if all pointwise stabilisers of cells are compact;
    \item \emph{with finitely many orbits} if the induced $G$-action on~$\Ce(X)$ has finitely many orbits.
\end{itemize}

\smallskip

In case of group actions on simplicial complexes, we further expand the terminology as follows. Firstly, given a simplicial $G$-action on a simplicial complex~$X$, a \emph{simplicial strict fundamental domain} is any simplicial subcomplex~$Y$ of~$X$ such that~$X=G\cdot Y$ and~$G\cdot \sigma\cap Y=\{\sigma\}$ for every~$\sigma\in Y$ (recall that we are regarding~$X$ and~$Y$ as collections of simplices).
For convenience, for every integer~$n\geq 1$, we introduce a family~$\caB_n$ of group actions on simplicial complexes as follows. 
\begin{defn}\label{Bn}
   A tuple~$(G,X,\Delta,\gamma)$ \emph{belongs to~$\caB_n$} if:
\begin{itemize}
    \item $X$ is a $(n+1)$-partite simplicial complex with associated type function~$\tau\colon X\to 2^{S}$, for some set~$S$ with~$|S|=n+1$;
    \item all positive dimensional links of~$X$ (included~$X$) are gallery-connected;
    \item $G$ is an lcsc group with a discrete, type-preserving, simplicial action on~$X$ that is transitive on~$X(n)$;
    \item $\Delta$ is the subcomplex of all faces of a given chamber~$\gamma$ of~$X$;
\end{itemize} 
\end{defn}
\begin{rem}\label{rem:defnBn}
   In \Cref{Bn}, the type function~$\tau$ restricted to~$\Delta$ is a bijection onto~$2^{S}$. Thus, for every subset~$J\subseteq S$,~$\Delta$ contains exactly one simplex~$\gamma_J$ of type~$J$.
\end{rem}
The family~$\hyperref[Bn]{\caB_n}$ should be regarded as a generalisation of the family of group actions satisfying conditions~$\caB1$--$\caB4$ considered in~\cite[\S 1]{DJ:build} and in~\cite[p.~845]{gro:curv}. In particular, the family~$\hyperref[Bn]{\caB_n}$ includes the following standard example:
\begin{ex}\label{ex:BNBn}
    Let~$G$ be an lcsc~group with a BN-pair~$(B,N)$ of type~$(W,S)$ such that~$B$ is an open subgroup of~$G$. 
    Let~$X=\Delta(G,B)$ be the simplicial building associated to the BN-pair~$(B,N)$. Denote by~$\gamma=\{1P_{S\setminus\{s\}}\mid s\in S\}$ its fundamental chamber, and let~$\Delta$ be the subcomplex of~$X$ of all faces of~$\gamma$. For every~$J\subseteq S$, the simplex~$\gamma_J=\{1P_{S\setminus \{j\}}\mid j\in J\}$ has stabiliser~$P_{S\setminus J}$ and its link in~$\Delta(G,B)$ is $P_{S\setminus J}$-equivariantly isomorphic to~$\Delta(P_{S\setminus J},S\setminus J)$.
    In particular,~$(G,X,\Delta,\gamma)\in \caB_{|S|-1}$.
\end{ex}

\begin{lem}\label{lem:Fams}
    Let~$(G,X, \Delta,\gamma)$ be in~$\hyperref[Bn]{\caB_n}$ for some~$n\geq 1$, and denote by~$\tau\colon X\to 2^S$ the type function of~$X$. Then: 
    \begin{enumerate}
        \item \label{lem:Fams0} The subcomplex~$\Delta$ is a strict fundamental domain for the $G$-action on~$X$.
        \item \label{lem:Fams1} For every~$J\subseteq S$, denote by~$G_J$ the stabiliser of~$\gamma_J$ in~$G$. Then
        \begin{equation*}
            X\longrightarrow \caC\caC(G, (G_{\{s\}})_{s\in S}),\qquad g\cdot\gamma_J\longmapsto \{gG_{\{s\}}\mid s\in J\}\quad\forall\,g\in G,\,J\subseteq S
        \end{equation*}
        is a $G$-equivariant simplicial isomorphism.
        \item \label{lem:Fams2} For every simplex~$\sigma\subsetneq \gamma$ we have~$G_\sigma=\langle G_\eta\mid \gamma\supseteq\eta\supseteq \sigma,\,\dim(\eta)=n-1\rangle$.
       \item \label{lem:Fams3} For every~$J\subseteq S$, the link $X_{\gamma_J}$ is a pure $|J|$-partite simplicial complex with type function~$\tau_J\colon X_\sigma(0)\to 2^{J}$ given by~$\tau_J(\sigma):=\tau(\sigma\cup \gamma_J)$.
       Moreover, $(G_{\gamma_J},X_{\gamma_J},\Delta_{\gamma_J},\gamma_J)$ belongs to~$\caB_{n-|J|}$. 
    \end{enumerate}
\end{lem}
\begin{proof}
\eqref{lem:Fams0} We first claim that~$X(k)=G\cdot \Delta(k)$ for every~$0\leq k\leq n$. Clearly,~$X(n)=G\cdot \gamma=G\cdot \Delta(n)$. For~$k\leq n-1$, since~$X$ is pure and $G$ acts transitively on~$X(n)$, for every~$\sigma\in X(k)$ there is~$g\in G$ such that~$\sigma\subseteq g\cdot \gamma$. As~$G$ acts type-preserving on~$X$,~$g^{-1}\cdot \sigma$ is the only face of~$\gamma$ of type~$J:=\tau(\sigma)$. We conclude that~$\sigma=g\cdot \gamma_J$.
This completes the proof that~$X=G\cdot \Delta$. 

The last part of the claim follows from the fact that~$\tau{\vert_\Delta}\colon \Delta\to 2^S$ is a bijection and that every~$g\in G$ acts as a type-preserving self-map on~$X$.
 
\eqref{lem:Fams1} It is a straightforward application of part~\eqref{lem:Fams0}.

\eqref{lem:Fams2} We proceed by induction on~$n-\dim(\sigma)\geq 1$. For the base case~$\dim(\sigma)=n-1$ there is nothing to prove.
For every~$J\subseteq S$, we keep denoting by~$G_J$ the stabiliser of~$\gamma_J$ in~$G$. The inductive step is a straightforward consequence of the following general observation: for every~$J\subseteq S$ with~$|J|\leq n-1$, we have
\begin{equation}\label{eq:Gsigmagen}
    G_J=\langle G_{J\cup \{s\}}\mid s\in S\setminus J\rangle.
\end{equation}
We now prove~\eqref{eq:Gsigmagen}.
By part~\eqref{lem:Fams1} and \Cref{lem:CCprop}\eqref{lem:CCprop2}, the link of~$X_\sigma$ is $G_J$-equivariant isomorphic to the coset complex~$\caC\caC(G_J,(G_{J\cup\{s\}})_{s\in S\setminus J})$. Hence, the latter coset complex is connected. Now \Cref{lem:CCprop}\eqref{lem:CCprop1} implies~\eqref{eq:Gsigmagen} and this concludes the argument.

\eqref{lem:Fams3} For every face~$\sigma\subseteq \gamma$ and every simplex~$\eta$ in~$X_\sigma$, notice that~$(X_\sigma)_\eta=X_{\sigma\cup \eta}$. Moreover, by part~\eqref{lem:Fams1} and \Cref{lem:CCprop}, every link~$X_\sigma$ of~$X$ can be written as an appropriate coset complex of~$G_\sigma$. The claim now follows.
\end{proof}

\section{Fréchet modules and cohomology theories}\label{s:FreCoho}

\subsection{Fréchet $G$-modules}\label{sus:freGmod}
A \emph{Fréchet space} is a topological $\C$-vector space that is metrisable, locally convex, and complete.
Hilbert spaces and, more generally, Banach spaces are examples of Fréchet spaces.  
Remarkably, the Open Mapping theorem holds for Fréchet spaces~\cite[Corollaries~2.12]{rudin:funAn}: every continuous surjective linear map between Fréchet spaces is open.

From now on in this section,~$G$ denotes a lcsc~group.
A \emph{Fréchet $G$-module} is a Fréchet space~$V$ with a continuous $G$-action~$G\times V\to V$ by continuous $\C$-linear automorphisms.
We typically denote a Fréchet $G$-module as a pair~$(\rho,V)$ where~$V$ is the underlying Fréchet space and~$\rho\colon G\to \End(V)$ is the group homomorphism induced by the continuous $G$-action on~$V$. Given~$v\in V$ and~$g\in G$, we sometimes write~$g\cdot v$ or~$\rho(g)v$ to denote the image of~$v$ via~$\rho(g)$.
For instance, if~$V$ is a Hilbert space and~$\rho\colon G\to \caU(V)$ is a strongly continuous unitary representation, then~$V$ has an induced structure of Fréchet $G$-module~\cite[Propositions~1.1-1.2]{cass:rep}.
\begin{ex}[\protect{\cite[Ch.~IX, \S 1.2]{bowal:coho}}]\label{ex:CXV}
    Let~$X$ be a locally compact $\sigma$-compact space with a continuous $G$-action by homeomorphisms. For every Fréchet $G$-module~$V$, the space of continuous functions~$C(X,V)$ from~$X$ to~$V$ is a Fréchet $G$-module with the compact-open topology and with $G$-action given by:
       \begin{equation*}
           (g\cdot f)(x):=g\cdot f(g^{-1}\cdot x) \qquad\forall\,g\in G,\,f\in C(X,V).
       \end{equation*}
\end{ex}

A \emph{Fréchet $G$-morphism} is a continuous $G$-equivariant $\C$-linear map between Fréchet $G$-modules. 
A \emph{Fréchet $G$-isomorphism} is a bijective Fréchet $G$-morphism whose inverse is a Fréchet $G$-morphism. It might be useful to recall the following standard characterisation, which follows from the Open Mapping theorem and~\cite[III, \S 6.2, Propositions~3 and~4]{bou:top}.
\begin{fact}\label{fact:FreIso}
Let~$V$ and~$W$ be Fréchet $G$-modules. Then a map~$f\colon V\to W$ is a Fréchet $G$-isomorphism if and only if~$f$ is continuous, bijective, $\C$-linear, and $G$-equivariant. 
\end{fact}

Below are some general properties of Fréchet $G$-modules that we use in this paper. 
\begin{fact}\label{fact:FreG}
    For every Fréchet $G$-module~$(\rho,V)$ the following hold:
\begin{enumerate}
    \item \label{fact:freG1} A closed $G$-invariant subspace of~$V$ is a Fréchet $G$-module.
    \item \label{fact:freG2} Given a $G$-invariant subspace~$C\subseteq V$, 
     \begin{equation*}
           V/C\text{ is a Fréchet $G$-module}\,\Longleftrightarrow\, V/C\text{ is Hausdorff}\,\Longleftrightarrow\,C\subseteq V\text{ is closed.}
       \end{equation*}
    \item \label{fact:freG3} Countable direct products of Fréchet $G$-modules are again Fréchet $G$-modules with respect to the entry-wise $G$-action.
    \item \label{fact:freG4} Let~$X$ be as in \Cref{ex:CXV}. Then every short exact sequence in the category of Fréchet $G$-modules
       \begin{equation*}
           \xymatrix{0\ar[r] & V\ar[r]^-{u} & W\ar[r]^-v & Z\ar[r] & 0}
       \end{equation*}
       induces an exact sequence in the category of Fréchet $G$-modules
       \begin{equation*}
           \xymatrix{0\ar[r] & C(X,V)\ar[r]^-{u^\ast} & C(X,W)\ar[r]^-{v^\ast} & C(X,Z)\ar[r] & 0,}
       \end{equation*}
      where~$u^\ast=u\circ\argu$ and~$v^\ast=v\circ \argu$.
       In particular, $\image(u^\ast)=\ker(v^\ast)$ is a closed subspace of~$C(X,W)$ that is isomorphic to~$C(X,V)$ as a Fréchet $G$-modules, and~$v^\ast$ induces an isomorphism of Fréchet $G$-modules:
       \begin{equation*}
          C(X,W)/\ker(v^\ast)\stackrel{\cong}{\longrightarrow} C(X,Z).
       \end{equation*}
\end{enumerate}
\end{fact}
\begin{proof}
    For~\eqref{fact:freG1} see~\cite[Theorem~1.2]{vogt:fre}. For~\eqref{fact:freG2} see~\cite[Theorem~1.41]{rudin:funAn} or~\cite[Theorem~1.2]{vogt:fre}. For~\eqref{fact:freG3} see~\cite[\S C, Remark~13]{nagy}. For~\eqref{fact:freG4} see~\cite[Ch.~IX, \S 1.2]{bowal:coho} and \Cref{fact:FreIso}.
\end{proof}

\subsection{Continuous cohomology of lcsc groups}\label{sus:contcoho}
Let~$G$ be an lcsc group and~$(\rho,V)$ be a Fréchet $G$-module. 
For every~$k\in \Z_{\geq 0}$, let
\begin{equation}\label{eq:contCochG}
    C^k(G,V):=C(G^{k+1},V)
\end{equation}
be the space of continuous functions from the $(k+1)$-fold direct product of~$G$ to~$V$. Every~$C^k(G,V)$ is endowed with the structure of Fréchet $G$-module as in~\Cref{ex:CXV}. 
The \emph{$k$-th coboundary map} $\partial^k\colon C^k(G,V)\to C^{k+1}(G,V)$ is defined as
\begin{equation}\label{eq:coboundG}
    (\partial^kf)(x_0,\dots, x_{k+1}):=\sum_{i=0}^{k+1}(-1)^if(x_0,\dots, \check{x}_i,\dots, x_{k+1}).
\end{equation}
By standard arguments,~$\partial^k$ is continuous and $G$-equivariant. The assignment in~\eqref{eq:coboundG} induces a cochain complex of Fréchet spaces and morphisms between them:
\begin{equation*}
    \xymatrix{\{0\} \ar[r]^-{\partial^{-1}} & C^0(G,V)^G\ar[r]^-{\partial^0} & \dots \ar[r]^-{\partial^{k-1}} & C^k(G,V)^G\ar[r]^-{\partial^k} & C^{k+1}(G,V)^G\ar[r]^-{\partial^{k+1}} & \dots}
\end{equation*}
For every~$k\geq 0$, the quotient space
$$\hH_{\ct}^k(G,\rho)= \hH_{\ct}^k(G,V):=\ker(\partial^k\vert_{C^k(G,V)^G})/\image(\partial^{k-1}\vert_{C^{k-1}(G,V)^G})$$ is called the \emph{$k$-th continuous cohomology} of~$G$ with coefficients in~$(\rho,V)$. Notice that~$\hH^k(G,\rho)$ is generally not Hausdorff (see~\Cref{fact:FreG}\eqref{fact:freG2}).

\smallskip

In fact, up to topological isomorphisms, one can prove (see~\cite[Ch.~IX, \S 1.5]{bowal:coho}) that~$\hH_{\ct}^k(G,\rho)$ is the $k$-th cohomology space of the cochain complex
\begin{equation*}
    \xymatrix{\{0\}\ar[r]^-{\partial^{-1}} & (E_0)^G\ar[r]^-{\partial^0} & \dots \ar[r]^-{\partial^{k-1}} & (E_k)^G\ar[r]^-{\partial^k} & (E_{k+1})^G\ar[r]^-{\partial^{k+1}} & \dots}
\end{equation*}
obtained by applying the fixed point functor~$(\argu)^G$ to any \emph{augmented strongly injective resolution} of~$V$ in the category of Fréchet $G$-modules:
\begin{equation*}
    E_\bullet:\qquad \xymatrix{\{0\} \ar[r]& V\ar[r]^-{\veps} & E_0\ar[r]^-{\partial^0} & \dots \ar[r]^-{\partial^{k-1}} & E_k\ar[r]^-{\partial^k} & E_{k+1}\ar[r]^-{\partial^{k+1}} & \dots}.
\end{equation*}
Recall that an exact sequence of Fréchet $G$-modules and Fréchet $G$-morphism~$E_\bullet$ of the form above is an \emph{augmented strongly injective resolution} of~$V$ (or determines a \emph{strongly injective resolution}~$E_0\to E_1\to \dots$ of~$V$) if every~$E_i$ is \emph{strong injective} as a Fréchet $G$-module, i.e., for every injective Fréchet $G$-morphism~$\varphi\colon A\hookrightarrow B$ such that the inclusion maps~$\ker(\varphi)\hookrightarrow A$ and~$\image(\varphi)\hookrightarrow B$ admit continuous $G$-invariant left inverses, every Fréchet $G$-morphism~$f\colon A\to E_i$ extends to a Fréchet $G$-morphism~$\wtil{f}\colon B\to E_i$. For details, see~\cite[Ch.~IX, \S 1.5]{bowal:coho} or~\cite[\S 2]{blanc}.
For instance, the sequence
\begin{equation}\label{eq:augmCoch}
   \xymatrix{\{0\} \ar[r]& V\ar[r]^-{\veps} & C^0(G,V)\ar[r]^-{\partial^0} & \dots \ar[r]^-{\partial^{k-1}} & C^k(G,V)\ar[r]^-{\partial^k} & C^{k+1}(G,V)\ar[r]^-{\partial^{k+1}} & \dots},
\end{equation}
where~$\veps\colon V\hookrightarrow C^0(G,V)$ is the constant map embedding, is an augmented strongly injective resolution of~$V$~\cite[Ch.~IX, Lemma~1.6]{bowal:coho}.

It follows from standard arguments that countable products of strongly injective Fréchet $G$-modules (which are Fréchet $G$-modules by \Cref{fact:FreG}\eqref{fact:freG3}) is again strongly injective. Using this, one may prove the following.

\begin{prop}\label{prop:HctProd}
    Let~$(V_i)_{i\in \Z}$ be a countable family of Fréchet $G$-modules, and let~$V:=\prod_{i\in \Z}V_i$ be endowed with the structure of product Fréchet $G$-module as in \Cref{fact:FreG}\eqref{fact:freG3}. Then, for every~$k\geq 0$, 
    $\hH_{\ct}^q(G,V)\cong_\C \prod_{i\in \Z}\hH_{\ct}^q(G,V_i)$ as topological $\C$-vector spaces.
\end{prop}
\begin{proof}[Proof sketch]
    Both~\eqref{eq:augmCoch} and the following sequence
\begin{equation*}
   \xymatrix{\{0\} \ar[r]& V\ar[r]^-{\veps} & \prod_{i\in \Z}C^0(G,V_i)\ar[r]^-{\prod_i\partial_i^0} & \dots \ar[r]^-{\prod_i\partial_i^{k-1}} & \prod_{i\in \Z}C^k(G,V_i)\ar[r]^-{\prod_i\partial_i^k} & \prod_{i\in \Z}C^{k+1}(G,V_i)\ar[r]^-{\prod_i\partial_i^{k+1}} & \dots}
\end{equation*}
  determine two strongly injective resolutions of~$V$. One now observes that the canonical map~$\varphi_k\colon C^k(G,V)\longrightarrow \prod_{i\in \Z}C^k(G,V_i)$, $k\geq 0$, determine an isomorphisms of Fréchet $G$-modules (recall \Cref{fact:FreIso}), and
  $\varphi_{k+1}\circ\partial^k=\prod_i\partial_i^k\circ \varphi_k$ for every~$k\geq 0$. 
\end{proof}

Apart from a few exceptions, such as for compact metrisable groups, the continuous cohomology of~$G$ is usually difficult to compute simply by applying its definition. As proved in various ways in~\cite{DJ:build, BallSwiat, Op15:vanish, gro:curv}, if~$G$ has a suitable action on a countable simplicial complex~$X$, it is possible to describe the vanishing of~$\hH_{\ct}^\bullet(G,\rho)$ in terms of the vanishing of so-called \emph{twisted simplicial cohomology} of $X$ with coefficients in~$(\rho,V)$ (see \Cref{sus:twistCoho} for the definition). In \Cref{s:mainCoho} we will generalise some results of the previously mentioned citations, producing new criteria for higher properties~T.

\smallskip

Speaking of higher properties~T, for future reference we recall the following permanence properties.
\begin{thm}[\protect{\cite{basa23}, \cite{bbbs:ban}}]\label{thm:Tkperm}
    Let~$G$ be an lcsc~group and let~$k\geq 1$.
    \begin{enumerate}
        \item \label{thm:Tkperm1} If~$G$ is compact then~$G$ has property~$\Tk$.
        \item \label{thm:Tkperm2} Let~$K\unlhd G$ be a compact normal subgroup. Then~$G$ has property~$\Tk$ (resp.~$(T_k)$) if and only if~$G/N$ has property~$\Tk$ (resp.~$(T_k)$).
        \item \label{thm:Tkperm3} Let~$H\leq G$ be a closed cocompact subgroup. Then~$G$ has property~$\Tk$ if and only if~$H$ has property~$\Tk$. Moreover, if~$G$ has property~$(T_k)$ then~$H$ has property~$(T_k)$.
    \end{enumerate}
\end{thm}
\begin{rem}
    In \Cref{thm:Tkperm}, part~\eqref{thm:Tkperm1} is  recorded~\cite[Lemma~3.3]{basa23}. Part~\eqref{thm:Tkperm2} is~\cite[Lemma~3.4]{basa23}. For part~\eqref{thm:Tkperm3} apply~\cite[Theorem~3.6, Corollary~3.7]{bbbs:ban}.
\end{rem}

\subsection{Cellular (co)homology of CW complexes}\label{sus:CWcoho}
Let~$X$ be a CW complex with a given cell orientation, and let~$V$ be an $R$-module. We denote by
\begin{equation*}
   \xymatrix{\cdots \ar[r]^-{d_k} & C_k(X,V)\ar[r]^-{d_{k-1}} & \cdots\ar[r]^-{d_1} & C_1(X,V) \ar[r]^-{d_0}& C_0(X,V)\ar[r]^-{d_{-1}} & \{0\}}
\end{equation*}
and
\begin{equation}\label{eq:cochainComplx}
   \xymatrix{\{0\} \ar[r]^-{d^{-1}} &C^0(X,V) \ar[r]^-{d^0} & C^1(X,V)\ar[r]^-{d^1} & \cdots \ar[r]^-{d^{k-1}} & C^k(X,V)\ar[r]^-{d^k} & \cdots }
\end{equation}
the \emph{cellular chain complex} and the \emph{cellular cochain complex} of~$X$ over~$V$, respectively. For every integer~$k\geq 0$, we denote by
$$\hH_k(X,V):=\ker(d_k)/\image(d_{k+1})\quad\text{and}\quad \hH^k(X,V):=\ker(d^k)/\image(d^{k-1})$$ 
the \emph{$k$-th cellular homology} and the \emph{$k$-th cellular cohomology} of~$X$ over~$V$, respectively.
One may refer to~\cite[\S 2.3, \S 12.1]{GeoBook} for a detailed background. In the following, we recall only those details that will be used in some arguments throughout. That is, for every~$k\geq 0$, recall that
\begin{equation}\label{eq:cochainNice}
    C^k(X,V):=\prod_{\sigma\in X(k)}V.
\end{equation}
For every~$f=\big(f(\sigma)\big)_{\sigma\in X(k)}\in C^{k}(X,V)$, the map~$d^kf\in C^{k+1}(X,V)$ is defined as
\begin{equation}\label{eq:coboundNice}
    (d^kf)(\tau):=\sum_{\sigma\subset \tau,\,\sigma\in X(k)}[\tau:\sigma]f(\sigma)\qquad \forall\,\tau\in X(k+1),
\end{equation}
where~$[\tau:\sigma]$ is the \emph{incidence number} of~$\sigma$ in~$\tau$~\cite[\S IV.5]{massey}. Recall that~$[\tau:\sigma]=0$ whenever~$\sigma$ is not contained in~$\tau$~\cite[Proposition~2.5.8]{GeoBook}. 
For instance, if~$X$ is a simplicial complex with orientation induced by a total ordering~$\leq$ on~$X(0)$, for all simplices~$\sigma,\tau$ with~$\tau=\{v_0<\dots < v_k\}$ and~$\dim(\sigma)=k-1$ we have
\begin{equation*}
    [\tau:\sigma]=\left\{
    \begin{array}{cl}
        (-1)^i & \text{if }\sigma=\tau\setminus\{v_i\}, \\
        0 & \text{otherwise;}
    \end{array}
    \right.
\end{equation*}
see~\cite[Proposition~5.4.2]{GeoBook}.
In particular, the cellular (co)chain complex of~$X$ over~$V$ is isomorphic to the corresponding simplicial (co)chain complex over~$V$.

\begin{thm}[Universal coefficient theorem; \protect{\cite[Theorem~12.5.1]{GeoBook}}]\label{thm:UCT}
Let~$R$ be a PID. Then, for every integer~$k\geq 0$, there is a short exact sequence of $R$-modules:
    \begin{equation*}
      \xymatrix{\{0\} \ar[r] & \Ext_R\big(\hH_{k-1}(X,R),V\big)\ar[r] & \hH^k(X,V)\ar[r] & \hom_R\big(\hH_k(X,R),V\big)\ar[r] & \{0\}}
    \end{equation*}
    with the convention that~$\hH_{-1}(X,R)=\{0\}$. If~$R$ is a field, this yields an $R$-linear isomorphism
    \begin{equation*}
        \hH^k(X,V)\cong_R \hom_R\big(\hH_k(X,R),V\big).
    \end{equation*}
\end{thm}

Given an integer~$k\geq 0$, we say that~$X$ is \emph{$k$-acyclic over~$V$} if~$\hH_0(X,V)=V$ and~$\hH_\ell(X,V)=\{0\}$ for all~$1\leq \ell\leq k$. This is the homological counterpart of asking that~$X$ is \emph{$k$-connected}, i.e., $X$ path-connected and~$\pi_\ell(X)=\{0\}$ for every~$1\leq \ell\leq k$. In particular, $X$ is said to be \emph{acyclic over~$V$} if it is $k$-acyclic over~$V$ for every~$k\geq 0$. 

By Hurewicz theorem~\cite[Theorem~4.32]{hatcher} and the Universal Coefficient Theorem recalled before, for every $R$-module~$V$ one easily deduces that
  \begin{equation}\label{eq:remcoho}
       X\text{ is $k$-connected}\,\Longrightarrow\,X\text{ is $k$-acyclic over~$\Z$}\,\Longrightarrow\, X\text{ is $k$-acyclic over~$V$.}
  \end{equation}

\begin{rem}\label{rem:cochFré}
   Let~$G$ be an lcsc group with a continuous cellular action on a countable CW complex~$X$. By \Cref{rem:Gorient}, we may set a cell orientation on~$X$ that is preserved by the $G$-action on~$X$. For every Fréchet $G$-module~$V$, each space~$C^k(X,V)$ is a Fréchet $G$-module as given in \Cref{ex:CXV}. Moreover, it is easily checked that each coboundary map~$d^k$ is continuous.
   
Since the $G$-action on~$X$ is orientation-preserving, by~\cite[Propositions~2.5.6-2.5.7]{GeoBook} we also have
   \begin{equation*}
       [g\cdot\sigma:g\cdot\tau]=[\sigma:\tau]
   \end{equation*}
   for all~$g\in G$, all cells~$\sigma$ and every codimension-$1$ face~$\tau\subset \sigma$.
Thus each~$d^k$ is also~$G$-invariant. This latter property will be crucial in \Cref{sus:twistCoho}.
\end{rem}

\subsection{Twisted cellular cohomology of CW complexes}\label{sus:twistCoho}
For this subsection, let~$X$ be a countable CW complex, and let~$G$ be an lcsc~group with a continuous (e.g., discrete) cellular action on~$X$. Set a cell orientation on~$X$ that is preserved by the $G$-action on~$X$ (see \Cref{rem:Gorient}).
Let also~$(\rho,V)$ be a Fréchet $G$-module.

\smallskip

For every~$k\geq 0$, regard~$X(k)$ as a discrete $G$-set and consider the Fréchet $G$-module structure on~$C^k(X,V)$ as given in \Cref{ex:CXV}. Since~$X(k)$ is discrete, the compact-open topology on~$C^k(X,V)=C(X(k),V)=\prod_{\sigma\in X(k)}V$ coincides with the product topology. Thus, we may write~$C^k(X,V)=\prod_{\sigma\in X(k)}V$ without ambiguities. Recall also that, by~\Cref{rem:cochFré}, each map~$d^k$ is continuous and $G$-equivariant. For every~$k\geq 0$ let
\begin{equation}\label{eq:coboundGNice}
  C^k(X,V)^G:=\Bigg\{\big(f(\sigma)\big)_\sigma\in \prod_{\sigma\in X(k)}\!V\,\Bigg\vert\,\forall\,g\in G,\,\sigma\in X(k),\,f(g\cdot \sigma)=g\cdot f(\sigma)\Bigg\}.
\end{equation}
Hence, the cochain complex in~\eqref{eq:cochainComplx} induces the following cochain complex of Fréchet spaces and continuous linear maps:
\begin{equation}\label{eq:twistXcochain}
   \xymatrix{\{0\}\ar[r]^-{d^{-1}} & C^0(X,V)^G \ar[r]^-{d^0} & C^1(X,V)^G\ar[r]^-{d^1} & \cdots \ar[r]^-{d^{k-1}} & C^k(X,V)^G\ar[r]^-{d^k} & \cdots }
\end{equation}
that we call the \emph{twisted cochain complex} of~$X$ with coefficients in~$(\rho,V)$. For every~$k\geq 0$, the \emph{$k$-th twisted cohomology}~$\hH^k(X,\rho)$ of~$X$ with coefficients in~$(\rho,V)$ is the $k$-th cohomology group of the cochain complex in~\eqref{eq:twistXcochain}.

\smallskip

In this paper (especially towards \Cref{s:mainCoho} and \Cref{sus:TkRGD}), we will be mainly interested in the vanishing of~$\hH^\ell(X,\rho)$, rather than a complete description of them. To this end, in the remaining part of the section, we recall some vanishing criteria for~$\hH^\ell(X,\rho)$ present in the literature, taking the change to weaken the connectivity hypotheses on~$X$ assumed in previous versions. 

We first recall the following vanishing criterion due to Oppenheim in~\cite[\S 3.2]{Op15:vanish}, which in turn generalises~\cite[Theorem~2.5]{BallSwiat}\footnote{In turn, \cite[Theorem~2.5]{BallSwiat} is a generalisation of what is usually called "Garland's vanishing theorem" (see~\cite{garland} and~\cite[Théorème~1.5]{borel}).}, and here presented in a slightly different form (see \Cref{rem:OppVan}).
\begin{thm}[\protect{\cite[Corollary~3.6]{Op15:vanish}}]\label{thm:OppVan}
    Let~$n\geq 1$ be an integer and let~$X$ be a locally finite, pure $n$-dimensional simplicial complex such that all links of~$X$ of positive dimension ($X$ included) are connected. Assume that
    \begin{equation}\label{eq:OppVan}
        \sup_{\sigma\in X(n-2)}\lambda_2(X_\sigma)<\frac{1}{n},
    \end{equation}
    where~$\lambda(X_\sigma)$ is the second largest eigenvalue of the normalised weighed adjacency matrix of~$X_\sigma$ (see \Cref{sus:matrices}). Then the following hold.
    \begin{enumerate}
        \item \label{thm:OppVan1} Let~$G$ be a unimodular lcsc group acting on~$X$ with a proper, discrete, simplicial action with finitely many orbits. Then, for every unitary representation~$(\rho,V)$ of~$G$ we have~$\hH^\ell(X,\rho)=\{0\}$ for all~$1\leq \ell\leq n-1$.
        \item \label{thm:OppVan2} If~$X$ is finite, then~$X$ is $(n-1)$-acyclic over any field of characteristic zero.
    \end{enumerate}
\end{thm}
\begin{rem}\label{rem:OppVan}
The statement of \Cref{thm:OppVan} merges~\cite[Corollary~3.6]{Op15:vanish} and~\cite[Theorem~2.5]{BallSwiat} with a few adaptations to our notation. To hopefully streamline the comparison, we clarify a few points below.
\begin{enumerate}
    \item The spectral condition in~\eqref{eq:OppVan} is equivalent to the one given in~\cite[Corollary~3.6]{Op15:vanish} for~$k=n-1$ and~$r=0$. In~\cite{Op15:vanish}, the equivalent condition is given in terms of the smallest positive eigenvalue of the weighted Laplacian~$\lambda(X_\sigma)$ of the link~$X_\sigma$. The initial weight set on~$X$ that we consider is~$w\colon X\to \Z_{\geq 0}$ given by~$w(\sigma):=\#\{\gamma\in X(n)\mid\gamma\supseteq\sigma\}$ for all~$\sigma\in X$. This is also the weight considered in~\cite{BallSwiat}. To see the previously mentioned equivalence, recall that for all~$\sigma\in X(n-2)$ we have~$\lambda_2(X_\sigma)=1-\lambda(X_\sigma)$ (see~\cite[pp.~629-630]{BallSwiat} for details). Moreover, as already observed in~\cite[Example~2.3]{Op15:vanish}, our definition of~$\lambda_2(X_\sigma)=1-\lambda(X_\sigma)$ and our choice of the weight~$w$ on~$X$ force the constant~$C_{n-1}$ to be equal to~$1$. In summary,
    \begin{equation*}
       \sup_{\sigma\in X(n-2)}\lambda_2(X_\sigma)<\frac{1}{n}\,\Longleftrightarrow\,\inf_{\sigma\in X(n-2)}\lambda(X_\sigma)>C_{n-1}\frac{n-1}{n}. 
    \end{equation*}
    \item As already partly observed in~\cite[Remark~3.3]{Op15:vanish}, the background hypothesis of contractibility of~$X$ assumed in~\cite[Corollary~3.6]{Op15:vanish} (see~\cite[p.~81]{Op15:vanish}) is unnecessary for the proof of the part of the statement we cite. Indeed, the only place where the contractibility of~$X$ is used in~\cite{Op15:vanish} is to have~$\hH_{\ct}^k(G,\rho)\cong \hH^k(X,\rho)$ (see the last part of the proof of~\cite[Theorem~2.6]{Op15:vanish}). 

   \item The statement in~\eqref{thm:OppVan2} is a straightforward  application of part~\eqref{thm:OppVan1} to the case when~$G$ is the trivial group and~$(\rho,V)$ is the trivial representation of~$\{1\}$ over~$V=\C$. This yields the $(n-1)$-acyclicity over~$\C$. The Universal Coefficient theorem (see \Cref{thm:UCT}) yields the $(n-1)$-acyclicity over~$\Q$ and hence over any field of characteristic zero.
\end{enumerate}
\end{rem}

Another vanishing criterion for~$\hH^\ell(X,\rho)$ is~\cite[Theorem~3.3(i)]{gro:curv} which is, in turn, based on~\cite[Theorem~5.1]{DJ:build}. We recall it below, and take the occasion to minimise the connectivity hypotheses on~$X$ assumed in both the cited results. This weakening will be crucial in view of one of the key results of the paper: \Cref{thm:U+propTk}.
\begin{thm}[\protect{\cite[Theorem~5.1]{DJ:build}, \cite[Theorem~3.3(i)]{gro:curv}}]\label{thm:GRO}
   Let~$(G,X,\Delta,\gamma)$ be in~\hyperref[Bn]{$\mathcal{B}_n$} for some~$n\geq 1$ and let~$(\rho,V)$ be a unitary representation of~$G$ with the following property: there is a family of closed subspaces~$(V^\eta)_{\eta\subseteq \gamma}$ such that, for every~$\sigma\subset \gamma$ we have
   \begin{equation}\label{eq:VsigmaDec}
       V^{G_{\sigma}}=\bigoplus_{\eta\subseteq\sigma}V^\eta.
   \end{equation}
   Then 
   \begin{equation}\label{eq:HXrho}
        \hH^\ell(X,\rho)=\left\{\begin{array}{ll}
           V^G  & \text{ if}~\ell=0, \\
           \{0\}  & \text{ if}~1\leq\ell\leq n-1,\\
           V^{\Delta} &\text{ if}~\ell=n,
        \end{array}\right.
    \end{equation}
    where~$V^\Delta$ is the orthogonal complement of~$\sum_{\sigma\subsetneq \gamma}V^{G_\sigma}$ in~$V^{G_\gamma}$.
\end{thm}

As proved in~\cite[Theorem~2.6]{gro:curv}, the condition~\eqref{eq:VsigmaDec} holds for every~$\sigma$ if a suitable "cosine matrix" associated to~$(\rho,V)$ and~$\Delta$.
This will be an instance of the following more general definition:

\begin{defn}\label{defn:CosMatRho}
Let~$V_0,\dots, V_n$, $n\geq 0$, be closed subspaces of a Hilbert space~$V$. The \emph{cosine matrix} associated to the tuple~$(V_0,\dots, V_n)$ is the real $(n+1)$-dimensional matrix~$A(V_0,\dots, V_n)=[A_{ij}]_{0\leq i,j\leq n}$ defined as
\begin{equation*}
        A_{ij}:=\left\{\begin{array}{cl}
            0 & \text{if }i=j, \\
            -\cos\angle (V_i,V_j) & \text{otherwise}, 
        \end{array}\right.
    \end{equation*}
    where~$\cos\angle (V_i,V_j)$ is the \emph{cosine of the angle}\footnote{This is a specialisation of the general notion by of cosine of the angle~$\cos\angle (W_1,W_2)$ of any pair of closed subspaces~$W_1,W_2$ of an arbitrary Hilbert space~$V$. See for instance~\cite[Definition~3.2]{kassabov}.} between~$V_i$ and~$V_j$, i.e., $\cos\angle (V_i,V_j):=1$ if~$V_i\subseteq V_j$ or~$V_j\subseteq V_i$ and, otherwise, 
    \begin{equation*}
    \cos\angle (V_i,V_j):=\sup\{|\langle v,w\rangle|\colon v\in V_i,\,w\in V_j,\,|v|=|w|=1,\,v,w\in(V_i\cap V_j)^\perp\}.
\end{equation*}
\end{defn}

\begin{cor}\label{cor:GRO1}
    Let~$(G,X,\Delta,\gamma)$ be in~\hyperref[Bn]{$\mathcal{B}_n$} for some~$n\geq 1$. Let~$(\rho,V)$ be a unitary representation of~$G$ and, for every~$0\leq i\leq n$,
    let~$V_i:=V^{G_{\gamma_{S\setminus\{i\}}}}$ (see \Cref{rem:defnBn}).
    If~$A(V_0,\dots, V_n)$ is positive definite, then~\eqref{eq:HXrho} holds for~$\hH^\ell(X,\rho)$ for all~$0\leq\ell\leq n$.
\end{cor}

The cosine matrix~$A(V_0,\dots, V_n)$ associated to~$(\rho,V)$ and~$(G,X,\Delta,\gamma)$ as in \Cref{cor:GRO1} should not be confused with the cosine matrix~$A(X)$ of~$X$ as in \Cref{defn:CosMatX}. However, as already partly recorded in~\cite[Corollary~3.8]{gro:curv}, under suitable hypotheses on the~$G$-action on~$X$, the positive definiteness of~$A(X)$ implies the positive definiteness of~$A(V_0,\dots, V_n)$ for every unitary representation~$(\rho,V)$ of~$G$. Here we present the generalisation of~\cite[Corollary~3.8]{gro:curv} corresponding to our weaker connectivity hypotheses on~$X$ we assumed in \Cref{thm:GRO}.

\begin{prop}\label{prop:GRO2}
Let~$(G,X,\Delta,\gamma)$ be in~\hyperref[Bn]{$\mathcal{B}_n$} for some~$n\geq 1$. Assume also that the $G$-action on~$X(n)$ is proper and that all $1$-dimensional links of~$X$ are finite. Denote by~$A(X)$ the cosine matrix of~$X$ with respect to the given type function~$\tau\colon X\to 2^{\{0,\dots, n\}}$ of~$X$, and let~$V_i=V^{G_{\gamma_{S\setminus\{i\}}}}$ for every~$0\leq i\leq n$.
Then:
\begin{enumerate}
    \item \label{prop:GRO2a} Let~$(\rho,V)$ be a unitary representation~$(\rho,V)$ of~$G$. For all~$i,j\in\{0,\dots, n\}$ with~$i\neq j$,
    \begin{equation*}
        \cos\angle(V_i,V_j)\leq \lambda_2\big({X_{\gamma_{S\setminus\{i,j\}}}}\big),
    \end{equation*}
   where~$\lambda_2\big(G_{X_{\gamma_{S\setminus\{i,j\}}}}\big)$ is the second largest eigenvalue of the normalised weighted adjacency matrix of the link of the codimension-$2$ face~$\gamma_{S\setminus\{i,j\}}$ (see \Cref{sus:matrices}). In other words, we have the following entry-wise inequality:
   $$A(V_0,\dots, V_n)\geq A(X).$$ 
   Moreover, if~$A(X)$ is positive definite then so is~$A(V_0,\dots, V_n)$.
    \item \label{prop:GRO2b} If~$A(X)$ is positive definite, then~\eqref{eq:HXrho} holds for all unitary representations~$(\rho,V)$ of~$G$.
\end{enumerate}
\end{prop}
\begin{proof}
 \eqref{prop:GRO2a}  The first two parts of the claim are already recorded in~\cite[Lemma~3.7]{gro:curv} and the references therein. The only change to make is to observe that the connectivity hypotheses on~$X$ assumed here are sufficient to obtain the conclusion. To make this passage neat, we briefly outline the argument again. The crux is to observe that for all~$i\neq j$,
  \begin{equation}\label{eq:cosIneq}
        \cos\angle (V_i, V_j)\leq\lambda_2(X_{\gamma_{S\setminus\{i,j\}}}\big).
    \end{equation}
For this inequality, one can proceed for instance as in~\cite[Theorem~1.7, Corollary~5.4]{kaufop18}.
Finally,~\cite[Proposition~2.8]{gro:curv} implies the statement on the positive definiteness.
  
  
\eqref{prop:GRO2b} Combine~\eqref{prop:GRO2a} and \Cref{cor:GRO1}.
\end{proof}
\begin{rem}\label{rem:GRO2imp}
    In \Cref{prop:GRO2} notice that for all~$i\neq j$, the pointwise stabiliser of~$\gamma_{S\setminus\{i,j\}}$ is compact. Indeed, since the link of that simplex is finite and contains~$\gamma$, we have~$|G_{\gamma_{S\setminus\{i,j\}}}:G_\gamma|<\infty$ and~$G_\gamma$ is compact. To the author's knowledge, all the known proofs of \Cref{prop:GRO2}\eqref{prop:GRO2a} use in an essential way peculiar properties of the representation theory of compact groups (which, in our case, are the various~$G_{\gamma_{S\setminus\{i,j\}}}$).
\end{rem}

\Cref{thm:GRO} and its corollaries give criteria for the vanishing of~$\hH^\ell(X,\rho)$ in degree~$1\leq\ell\leq n$, and provide a fairly explicit description of these spaces in degrees~$\ell=0$ and~$\ell=n$. A priori, $V^\Delta$ might vanish: this is the case when~$(\rho,V)$ is the trivial representation of~$G$ on~$V=\C$. In the following, slightly generalising~\cite[\S 6]{DJ:build}, we give a standard class of examples for which~$V^\Delta$ is always non-zero.
\begin{ex}\label{ex:H2perp}
  Let~$G$ be an lcsc group with a BN-pair~$(B,N)$ of type~$(W,S)$ and such that~$B$ is an open subgroup of~$G$. Assume also that the associated thickness vector~$(q_s)_{s\in S}$ is finite, and that and~$q_s>|S|$ for every~$s\in S$.

  Consider the unitary representation~$(\rho,V)$ of~$G$ given by the Hilbert space
  \begin{equation}
      V=L^2(G/B,\C):=\Big\{f\colon G/B\to \C\,\Big\vert\,\textstyle{\sum_{xB\in G/B}}|f(xB)|^2<\infty\Big\}
  \end{equation}
  with inner product~$\langle f_1,f_2\rangle:=\sum_{xB\in G/B}f_1(xB)\overline{f_2(xB)}$, and by
  $$\rho\colon G\to \caU(V), \qquad \rho(g)f(xB):=f(g^{-1}xB)$$ for every~$x,g\in G$, $f\in H$.
Clearly,
  \begin{equation*}
      V^B=\{f\in V\mid f\text{ is constant of each $(B,B)$-double coset}\}
  \end{equation*}
and, for every~$s\in S$,
  \begin{equation*}
      V^{P_s}=\{f\in V\mid f\text{ is constant of each $(P_s,B)$-double coset}\}.
  \end{equation*}
By~\cite[Proposition~2.4.4]{bb:cox}, there is a set of representatives~${}^sW$ of the right cosets~$W_s\backslash W$ with the following property: for every~$w\in W$, there are unique~$u\in W_s$ and~$v\in {}^sW$ such that~$w=uv$. As a consequence,
\begin{equation*}
    G=\bigsqcup_{w\in W}BwB=\bigsqcup_{w\in {}^sW}P_swB
\end{equation*}
and thus the set~$\{\bbone_{P_swB}\mid w\in {}^sW\}$ of indicator functions over~$G$ is a orthogonal basis of~$V^{P_s}$.

Let~$V^\Delta$ be the orthogonal complement of~$\sum_{s\in S}V^{P_s}$ in~$V^B$. Recalling \Cref{ex:BNBn}, notice that~$(G,\Delta(G,B),\Delta,\gamma)\in \caB_{|S|-1}$,~$V^B=V^{G_\gamma}$ and~$\sum_{s\in S}V^{P_s}=\sum_{J\subseteq S,\,J\neq\emptyset}V^{P_J}=\sum_{\eta\subsetneq \gamma}V^{G_\gamma}$.

We claim that
\begin{equation}\label{eq:HDeltanon0}
    V^{\Delta}\neq \{0\}.
\end{equation}
Our proof strategy is a variation on theme of~\cite[\S 6]{DJ:build}. For every~$w\in W$, let~$\ell(w)$ be its word length in~$(W,S)$ and let~$q_w:=q_{s_1}\cdots q_{s_n}$, where~$w=s_1\cdots s_n$ is any reduced expression of~$w$ in~$(W,S)$. Clearly,
$$q_w\geq (q_{\min})^{\ell(w)},$$ 
where~$q_{\min}:=\min\{q_s\mid s\in S\}$.
Moreover, arguing as in \Cref{rem:GBspher} we have
$$|B:B\cap wBw^{-1}|=|BwB/B|=q_w.$$

Consider the function
\begin{equation*}
    \phi\colon G/B\longrightarrow \C\qquad \phi(bwB):=(-1)^{\ell(w)}q_w^{-1}\quad\forall\,b\in B.
\end{equation*}
Since~$\qmin>|S|$, we have
\begin{equation*}
    \sum_{xB\in G/B}\hspace{-0.2cm}|\phi(xB)|^2\leq \sum_{w\in W}\sum_{b\in B/B\cap wBw^{-1}} \hspace{-0.3cm}q_w^{-2}\leq \sum_{w\in W}(q_{\min})^{-\ell(w)}\leq \sum_{n=0}^\infty \Bigg(\frac{|S|}{q_{\min}}\Bigg)^n\hspace{-0.2cm}<\infty
\end{equation*}
and hence~$\phi\in V^B$. We prove more precisely that~$\phi$ belongs to~$V^\Delta$. To this end, it suffices to prove that
\begin{equation}\label{eq:phiOrth}
    \langle \phi,\bbone_{P_swB}\rangle=0
\end{equation}
for all~$s\in S$ and~$w\in{}^sW$. Since~$\ell(sw)>\ell(w)$, we have~$(BsB)(BwB)=BswB$~\cite[p.~320, part~(2)]{ab:build}. As~$P_s=B\sqcup BsB$, we deduce that
\begin{equation*}
    P_swB=BwB\cup (BsB)(BwB)=BwB\sqcup BswB.
\end{equation*}
Therefore, since~$l(sw)\in\{\ell(w)\pm 1\}$~\cite[Proposition~1.4.2(iii)]{bb:cox},
\begin{equation*}
    \langle \phi,\bbone_{P_swB}\rangle=\langle \phi,\bbone_{BwB}\rangle+\langle \phi,\bbone_{BswB}\rangle=(-1)^{\ell(w)}+(-1)^{\ell(sw)}=0.
\end{equation*}
This proves~\eqref{eq:phiOrth} and completes the proof of~\eqref{eq:HDeltanon0}. 
\end{ex}

\subsection{Cohomological spectral sequences}\label{sus:specSeq}
In this paper, we introduce the necessary background on spectral sequences towards \Cref{s:mainCoho}. 
We mainly follow~\cite{greenberg} and~\cite[Appendix~A]{ribzal}.
For general background, see also~\cite[\S VII.3]{brown:coho} and~\cite[\S 5.6]{weibel}.

Let~$R$ be a unital ring.
A \emph{double cochain complex}~$\caC=(C^{pq},\dd_I^{pq}, \dd_{I\!I}^{pq})_{(p,q)\in \Z^2}$ of $R$-modules consists of $R$-modules~$C^{pq}$ and $R$-linear morphisms~$\dd_I^{pq}\colon C^{pq}\to C^{p,q+1}$ and~$\dd_{I\!I}^{pq}\colon C^{pq}\to C^{p+1,q}$ such that the following holds for all integers~$p,q$:
    \begin{equation}\label{eq:doubleC}
        \dd_I^{p,q+1}\circ \dd_I^{pq}=0,\quad \dd_{I\!I}^{p+1,q}\circ \dd_{I\!I}^{pq}=0 \quad\text{and}\quad \dd_I^{p+1,q}\circ \dd_{I\!I}^{pq}=-\dd_{I\!I}^{p,q+1}\circ \dd_I^{pq}.
    \end{equation}
Whenever there is no ambiguity, we conventionally omit the superscripts on~$\dd_I$ and~$\dd_{I\!I}$.    
\begin{equation*}
             \xymatrix{
             & \quad\ar@{<..}[d] & \quad\ar@{<..}[d] & \quad\\
            \quad & C^{p,q+1}\ar@{<..}[l]\ar@{->}[r]^-{\dd_{I}}\ar@{<-}[d]_-{\dd_{I\!I}} & C^{p+1,q+1}\ar@{<-}[d]^-{\dd_{I\!I}}\ar@{..>}[r] & \quad\\
             \quad & C^{pq}\ar@{<..}[l]\ar@{->}[r]_-{\dd_{I}}\ar@{<..}[d] & C^{p+1,q}\ar@{..>}[r]\ar@{<..}[d] &\quad\\
             & \quad\ar@{..>}[u] & \quad\ar@{..>}[u] &\quad
             }
         \end{equation*}

In particular,~$\caC$ is said to be \emph{first-quadrant} if~$C^{pq}=\{0\}$ if~$p<0$ or~$q<0$.

\noindent Every double cochain complex~$\caC$ has an associated \emph{total complex}~$\caT=\caT(\caC)=(T^n,\dd^n)_{n\in \Z}$ defined by
\begin{equation}\label{eq:Tcomplx1}
    T^n:=\bigoplus_{(p,q)\in \Z^{2}\colon p+q=n}C^{pq}
\end{equation}
and
\begin{equation}\label{eq:Tcomplx2}
    \dd^n\colon T^n\to T^{n+1},\quad \dd^n\vert_{C^{pq}}:=\dd_I^{pq}+\dd_{I\!I}^{pq}\colon C^{pq}\to C^{p,q+1}\oplus C^{p+1,q}.
\end{equation}
The anti-commutativity condition in~\eqref{eq:doubleC} guarantees that~$\dd^{n+1}\circ \dd^n=0$ for every~$n\in \Z$. Hence~$\caT$ is a cochain complex, and we denote by~$\hH^n(\caT)$ its $n$-th cohomology group.

As recalled in the following \Cref{fact:SPSdouble}, the cohomology of the total complex~$\caT=\caT(\caC)$ can be described in terms of two spectral sequences~$({}^I\!\opE_{r}^{pq}, {}^I\delta_r^{pq})_{r\in\Z_{\geq 0},\,(p,q)\in \Z^2}$ and~$({}^{I\!I}\!\opE_{r}^{pq}, {}^{I\!I}\delta_r^{pq})_{r\in \Z_{\geq 0},\,(p,q)\in \Z^2}$. Before stating the fact, we briefly recall the necessary background on spectral sequences. 

A \emph{spectral sequence}\footnote{The complete terminology should be "bigraded cohomological spectral sequence". Indeed, there is a more general notion of spectral sequence, which does not necessarily stem from a double cochain complex and does not necessarily require a bigrading. Moreover, the adjective "cohomological" distinguishes these spectral sequences from their dual counterparts: the homological spectral sequences.} of $R$-modules~$\caE=(\opE_r^{pq}, \delta_r^{pq})_{r\in \Z_{\geq 0},\,(p,q)\in \Z^2}$ consists of~$R$-modules~$\opE_r^{pq}$ and $R$-linear homomorphisms~$\dd_r^{pq}\colon \opE_r^{pq}\to \opE_r^{p+r,q+1-r}$ such that, for all~$r\geq 0$ and~$(p,q)\in \Z^2$,
    \begin{equation}\label{eq:spec}
        \dd_r^{p+r,q+1-r}\circ \dd_r^{p,q}=0\quad \text{and}\quad \opE_{r+1}^{pq}=\ker(\dd_r^{pq})\Big/\image(\dd_r^{p-r,q+r-1}).
    \end{equation}
Notice that~$\opE_{r+1}^{pq}$ is a subquotient of~$\opE_r^{pq}$, and thus
\begin{equation}\label{eq:subquots}
    \opE_r^{pq}=\{0\}\,\Longrightarrow\,\opE_{r+\ell}^{pq}=\{0\}\quad\forall\,\ell\geq 0.
\end{equation}

Like for double cochain complexes, a spectral sequence~$\caE$ is said to be \emph{first-quadrant} if, for every~$r\geq 0$, $\opE_r^{pq}=\{0\}$ if~$p<0$ or~$q<0$. If~$\caE$ is first-quadrant, by~\eqref{eq:subquots} one easily deduces (see~\cite[Proposition~A.2.1(a)]{ribzal}) the following for all integers~$p,q\geq 0$:
\begin{equation}\label{eq:preEinf}
    \opE_r^{pq}=\opE_{\max\{p+1,q+2\}}^{pq}\qquad\forall\,r\geq\max\{p+1,q+2\}.
\end{equation}
As a consequence of~\eqref{eq:preEinf}, for all integers~$p,q\geq 0$ we may define
\begin{equation*}
    \opE_\infty^{pq}:=\opE_{\max\{p+1,q+2\}}^{pq}.
\end{equation*}
A first-quadrant spectral sequence~$\caE$ is said to \emph{converge} to a given sequence of $R$-modules~$\caH=(\opH^n)_{n\in \Z_{\geq 0}}$, written
\begin{equation*}
    \opE_r^{pq}\Rightarrow \opH^{p+q},
\end{equation*}
if, for all~$n\in \Z$, $\opH^n$ admits a filtration of $R$-submodules
\begin{equation}\label{eq:Hnsubquot}
    \{0\}=(\opH^n)^{n+1}\subseteq (\opH^n)^n\subseteq \dots \subseteq (\opH^n)^1\subseteq (\opH^n)^0=\opH^n
\end{equation}
such that
\begin{equation}\label{eq:EinftConv}
    \opE_\infty^{pq}\cong (\opH^{p+q})^p\big/(\opH^{p+q})^{p+1} \qquad\forall\,(p,q)\in \Z_{\geq 0}^2.
\end{equation}
As a consequence, if~$\opE_r^{pq}\Rightarrow \opH^{p+q}$ then the following holds for all~$n\geq 0$:
\begin{equation}\label{eq:vanishHnEinfty}
    \opH^n=\{0\}\,\Longleftrightarrow\, \opE_\infty^{p,n-p}=\{0\}\qquad\forall\,p\in\{0,\ldots, n\}.
\end{equation}

We are now ready to state \Cref{fact:SPSdouble}.
\begin{fact}[\protect{\cite[\S A.4]{ribzal}, \cite[Theorem~8]{greenberg}}]\label{fact:SPSdouble}
    Let~$\caC=(C^{pq},\dd_I^{pq}, \dd_{I\!I}^{pq})_{(p,q)\in \Z^2}$ be a double cochain complex of~$R$-modules with associated total complex~$\caT=(T^n,d^n)_{n\in \Z}$. Then there are two explicitly defined spectral sequences~${}^I\!\opE_{r}^{pq}$ and~${}^{I\!I}\!\opE_r^{pq}$, where~$r\geq 0$, such that
    \begin{equation}\label{eq:SPSdouble1}
        \begin{array}{lllll}
        {}^I\!\opE_{0}^{pq}=C^{pq},  & \quad & {}^I\!\opE_{1}^{pq}=\hH_{I\!I}^q(C^{p\bullet}); & \quad & {}^I\!\opE_{2}^{pq}=\hH_{I}^p(\hH_{I\!I}^q(C^{\bullet\bullet}));\\
          {}^{I\!I}\!\opE_{0}^{pq}=C^{qp},  & \quad & {}^{I\!I}\!\opE_{1}^{pq}=\hH_{I}^q(C^{\bullet p}); & \quad & {}^{I\!I}\!\opE_{2}^{pq}=\hH_{I\!I}^p(\hH_{I}^q(C^{\bullet\bullet})).
        \end{array}
    \end{equation}
If~$\caC$ is first-quadrant then so are both~${}^I\!\opE_r^{pq}$ and~${}^{I\!I}\!\opE_r^{pq}$, and for all integers~$p,q\geq 0$ we have
     \begin{equation*}
         {}^I\!\opE_\infty^{pq}\Rightarrow \hH^{p+q}(\caT)\qquad\text{and}\qquad{}^{I\!I}\!\opE_\infty^{pq}\Rightarrow\hH^{p+q}(\caT).
     \end{equation*}
\end{fact}
Before we further specialise our discussion, we recall a standard property on spectral sequences that will be frequently used throughout.
\begin{lem}\label{lem:convSP}
     Let~$\caE=(\opE_r^{pq},\delta_r^{pq})_{r\in \Z_{\geq 0},\, (p,q)\in \Z^2}$ be a first-quadrant spectral sequence and let~$k\geq 1$ be an integer.
    Assume\footnote{If~$k=1$, this assumption results in an empty condition.} that~$\opE_1^{pq}=\{0\}$ for all~$0\leq p\leq k-2$ and~$q\geq 1$ such that~$p+q\leq k-1$. Then
        \begin{equation*}
            \opE_\infty^{p0}\cong\opE_2^{p0}\qquad\forall\,p\in\{1,\dots, k\}.
        \end{equation*}
\end{lem}
\begin{proof}
    Clearly~$\opE_\infty^{10}\cong\opE_2^{10}$. In particular, this proves the statement for~$k=1$. Let now~$2\leq p\leq k$ and~$2\leq r\leq p$. We claim that the extrema of the following sequence both vanish:
    \begin{equation}\label{eq:convSPS2A}
        \xymatrix{
        \opE_r^{p-r,r-1}\ar[rr]^-{\delta_r^{p-r,r-1}} && \opE_r^{p0} \ar[rr]^-{\delta_r^{p0}} && \opE_r^{p+r,1-r}}
    \end{equation}
   Notice that~$\opE_r^{p+r,1-r}=\{0\}$ because~$\caE$ is first-quadrant and~$r\geq 2$. By hypothesis we also have~$\opE_1^{p-r,r-1}=\{0\}$. By passage to subquotients,~$\opE_r^{p-r,r-1}=\{0\}$.
    Since the two extrema of~\eqref{eq:convSPS2A} vanish,~$\opE_{r+1}^{p0}=\ker(\delta_r^{p0})/\image(\delta_r^{p-r,r-1})\cong\opE_r^{p0}$ for all~$2\leq p\leq k$ and~$2\leq r\leq p$. We conclude that~$\opE_\infty^{p0}\cong\opE_{p+1}^{p0}\cong\opE_2^{p0}$ for all~$1\leq p\leq k$.
\end{proof}

\smallskip

In this paper, we are primarily interested in double cochain complexes of the following form.

\begin{defn}[\emph{The double cochain complex~$\caC(G,\caA)$}]\label{defn:CGA}
Let~$G$ be an lcsc group and let
\begin{equation*}
    \caA\colon \quad  \xymatrix{\{0\}\ar[r]^-{e^{-1}} & A^0 \ar[r]^-{e^0} & A^1\ar[r]^-{e^1} & \cdots \ar[r]^-{e^{p-1}} & A^p\ar[r]^-{e^p} & A^{p+1} \ar[r]^-{e^{p+1}} & \cdots }
\end{equation*}
be a cochain complex consisting of Fréchet $G$-modules and continuous $G$-equivariant coboundary maps, i.e., $e^{p+1}\circ e^p=0$ for every~$p\geq 0$. 
For every~$p\geq 0$ denote by
$$\hH^p(\caA):=\ker(e^p)\Big/\image(e^{p-1})$$ 
be the $p$-th cohomology group of~$\caA$ (as usual, endowed with the quotient topology). 

Define as the following first-quadrant double cochain complex~$\caC(G,\caA)=(C^{pq},\dd_I^{pq}, \dd_{I\!I}^{pq})$: for all~$p,q\geq 0$ let
\begin{equation}\label{eq:CpqA}
    C^{pq}:=C^q(G,A^p)^G,
\end{equation}
let
\begin{equation}\label{eq:dII}
    \dd_{I}^{pq}\colon C^q(G,A^p)^G\longrightarrow C^q(G,A^{p+1})^G,\quad (\dd_{I\!I}^{pq})(f):=e^p\circ f.
\end{equation}
and notice that the $G$-equivariance of~$e^p$ makes~\eqref{eq:dII} well-defined. Let also
\begin{equation}\label{eq:dI}
    \dd_{I\!I}^{pq}\colon C^q(G,A^p)^G\longrightarrow C^{q+1}(G,A^p)^G
\end{equation}
be the coboundary map as in~\eqref{eq:coboundG}.
\end{defn}

It can be easily checked that~$\caC(G,\caA)$ is a double cochain complex of $\C$-vector spaces, in fact, of Fréchet $G$-modules. Denote by~${}^I\!\opE_r^{pq}$ and~${}^{I\!I}\!\opE_r^{pq}$ the two spectral sequences associated to~$\caC(G,\caA)$ as in~\Cref{fact:SPSdouble}. Apart when~$r=0$, the spaces~${}^I\!\opE_r^{pq}$ and~${}^{I\!I}\!\opE_r^{pq}$ might \emph{not} be Fréchet $G$-modules (with the quotient topology); see \Cref{fact:FreG}\eqref{fact:freG2}.

\smallskip

The following proposition is extracted from the proof of~\cite[\S IX, Theorem~4.1]{bowal:coho}. Note that the theorem stated \emph{in loc.~cit.}~assumes stronger hypotheses on~$\caA$ than the ones needed for the result below (the cited result also gets a stronger conclusion than the one stated here). 

\begin{prop}\label{prop:niceE1}
   Let~$\caC(G,\caA)$ be as in~\Cref{defn:CGA} and let~$p,q\geq 0$ be integers. Then 
   \begin{equation}\label{eq:niceE11a}
       \ker(\dd_{I}^{qp})=C^p(G,\ker(e^q))^G.
   \end{equation}
Assume now there is an integer~$k\geq 1$ such that~$\hH^q(\caA)$ is Hausdorff for every~$0\leq q\leq k$. Then, for all~$p\geq 0$ and~$0\leq q\leq k$ we have
   \begin{equation}\label{eq:niceE11b}
       \image(\dd_{I}^{q-1,p})=C^p(G,\image(e^{q-1}))^G
   \end{equation}
   and there is a canonical isomorphism of Fréchet spaces 
   \begin{equation}\label{eq:niceE12}
{}^{I\!I}\!\opE_1^{pq}\cong C^p(G,\hH^q(\caA))^G.
   \end{equation}
\end{prop}
\begin{proof}
    The equality in~\eqref{eq:niceE11a} is straightforward.
    
    Let~$\hH^q(\caA)$ be Hausdorff for every~$0\leq q\leq k$, and fix arbitrary integers~$p\geq 0$ and~$0\leq q\leq k$.
    Both~$\ker(e^{q-1})$ and~$\image(e^{q-1})$ are closed in~$A^{q-1}$ because~$e^{q-1}$ is continuous and~$\hH^q(\caA)=\ker(e^q)/\image(e^{q-1})$ is Hausdorff (see \Cref{fact:FreG}\eqref{fact:freG2}). 
    Hence, the canonical inclusion~$\ker(e^{q-1})\hookrightarrow A^{q-1}$ and the map~$e^{q-1}$ induce the following exact sequence of Fréchet $G$-modules and Fréchet $G$-morphisms:
    \begin{equation}\label{eq:niceE14}
        \xymatrix{\{0\} \ar[r] & \ker(e^{q-1})\ar@{->}[r] & A^{q-1} \ar@{->}[r]^-{e^{q-1}} & \image(e^{q-1}) \ar[r] & \{0\}.}
    \end{equation}
  By~\cite[Chapter~IX, \S 1.4]{bowal:coho}, for every~$p\geq 0$ and every Fréchet $G$-module~$V$, the following is an isomorphism of Fréchet spaces:
    \begin{equation}\label{eq:Cuniv}
        \Psi_V\colon C^p(G,V)^G\longrightarrow C^{p-1}(G,V)\qquad(\Psi_Vf)(x_1,\dots, x_p):=f(1,x_1,x_1x_2,\dots, x_1x_2\cdots x_p),
    \end{equation}
    with the convention that~$C(G^0,V)=V$. Moreover, every Fréchet $G$-morphism~$e\colon V\to W$ induces a commuting square with invertible rows as follows:
   \begin{equation*}
       \xymatrix{C^p(G,V)^G\ar[r]^-{\Psi_V}\ar[d]_-{e\circ\argu} & C^{p-1}(G,V)\ar[d]^-{e\circ\argu}\\
       C^p(G,W)^G\ar[r]^-{\Psi_W} & C^{p-1}(G,W)}
   \end{equation*}
By the previous observation and by \Cref{fact:FreG}\eqref{fact:freG4}, the exact sequence in \eqref{eq:niceE14} yields an exact sequence of Fréchet spaces
    \begin{equation*}
        \xymatrix{\{0\} \ar[r] & C^p(G,\ker(e^{q-1}))^G\ar@{->}[r] & C^p(G,A^{q-1})^G \ar@{->}[r]^-{d_{I}^{q-1,p}} & C^p(G,\image(e^{q-1}))^G \ar[r] & \{0\}.}
    \end{equation*}
    and hence
    \begin{equation}\label{eq:niceE15}
    \dd_{I}^{q-1,p}\big(C^p(G,A^{q-1})^G\big)=C^p(G,\image(e^{q-1}))^G.
    \end{equation}
   This proves~\eqref{eq:niceE11b}.
    
    Finally, by \Cref{fact:FreG}\eqref{fact:freG2} the quotient~$\hH^{q}(\caA)=\ker(e^{q})/\image(e^{q-1})$ is a Fréchet $G$-module. Thus,
    \begin{equation*}
    \begin{split}
{}^{I\!I}\!\opE_1^{pq}& =\hH_I^q(C^{\bullet p}) =\ker(\dd_{I}^{qp})\Big/\image(\dd_{I}^{q-1,p})\stackrel{\clubsuit}{\cong} C^{p}\big(G,\ker(e^{q})\big)^G/C^{p}\big(G,\image(e^{q-1})\big)^G\\
& \stackrel{\ast}{\cong} C^{p-1}(G,\ker(e^q))\big/C^{p-1}(G,\image(e^{q-1}))\stackrel{\diamondsuit}{\cong} C^{p-1}(G,\hH^q(\caA))\stackrel{\ast}{\cong} C^p(G,\hH^q(\caA))^G.
\end{split}
\end{equation*}
Before, the passage~$\clubsuit$ is implied by~\eqref{eq:niceE11a} and \eqref{eq:niceE11b}, in the passages~$\ast$ we have used the isomorphism in~\eqref{eq:Cuniv}, and~$\diamondsuit$ follows from \Cref{fact:FreG}\eqref{fact:freG4} applied to the canonical exact sequence of Fréchet $G$-modules~$\{0\}\to \image(e^{q-1})\hookrightarrow \ker(e^q)\twoheadrightarrow \hH^p(\caA)\to \{0\}$.
\end{proof}

Combining \Cref{fact:SPSdouble} and \Cref{prop:niceE1} we deduce the following generalisation of~\cite[Chapter~IX, Corollary~4.2]{bowal:coho}.
\begin{cor}\label{cor:niceDGA}
    Let~$\caC(G,\caA)$ be as in~\Cref{defn:CGA}. Then for all integers~$p,q\geq 0$ we have
    \begin{equation}\label{eq:niceDGA1}
        {}^I\!\opE_1^{pq}\cong\hH_{\ct}^q(G,A^p)
    \end{equation}
    and
    \begin{equation}\label{eq:niceDGA2}
        {}^I\!\opE_2^{p0}\cong\frac{\ker\big(e^p\colon (A^p)^G\to (A^{p+1})^G\big)}{\image\big(e^{p-1}\colon (A^{p-1})^G\to (A^p)^G\big)}.
    \end{equation}
    
Furthermore, assume that there is an integer~$k\geq 1$ such that~$\hH^q(\caA)$ is Hausdorff for every~$0\leq q\leq k$. Then for all integers~$p\geq 0$ and~$0\leq q\leq k$ we have
    \begin{equation}\label{eq:niceDGA3}
{}^{I\!I}\!\opE_1^{pq}\cong C^p(G,\hH^q(\caA))^G\qquad\text{and}\qquad {}^{I\!I}\!\opE_2^{pq}\cong\hH_{\ct}^p(G,\hH^q(\caA)).
    \end{equation}
\end{cor}

\section{Groups actions on CW complexes and higher properties~T}\label{s:mainCoho}
Let~$X$ be a countable CW complex and let~$G$ be an lcsc group with a discrete (by default, admissible) cellular action on~$X$. Set a cell orientation of~$G$ so that the $G$-action on~$X$ is orientation-preserving.
Let also~$(\rho,V)$ be a Fréchet $G$-module. 
Throughout, we will freely use the notation introduced in Sections~\ref{sus:CWcoho} and~\ref{sus:twistCoho}.

\smallskip

In this section, we specialise \Cref{defn:CGA} to the case where~$\caA$ is the cellular cochain complex of~$X$ over~$V$; in symbols,
\begin{equation}\label{eq:AXrho}
    \caA=\caA_{(X,\rho)}=(C^p(X,V),d^p)_{p\geq 0}.
\end{equation}
The following result generalises~\cite[Lemma~3.5]{DJ:build} to discrete cellular $G$-actions on countable CW complexes. 
Note that the CW complex is allowed to be infinite and non-contractible, and the $G$-action on it may have infinitely many orbits.
\begin{lem}\label{lem:cohoPlain}
Let~$G$, $X$ and~$(\rho,V)$ be as at the beginning of the current section, and let~${}^I\!\opE_r^{pq}$ and~${}^{I\!I}\!\opE_r^{pq}$ be the canonical spectral sequences induced by~$\caC(G,\caA_{(X,\rho)})$ (see \Cref{fact:SPSdouble}). Then, for all integers~$p,q\geq 0$,
     \begin{equation}\label{eq:cohoP2a}
        {}^{I}\!\opE_2^{p0}\cong\hH^p(X,\rho)
    \end{equation}
    and there is a canonical $\C$-linear isomorphism 
    \begin{equation}\label{eq:cohoP2b}
         {}^{I}\!\opE_1^{pq}\cong \prod_{\sigma\in X(p)/G}\hH_{\ct}^q(G_\sigma,\rho\vert_{G_\sigma}).
    \end{equation}
     
Assume now there is an integer~$k\geq 1$ such that~$\hH^q(X,V)$ is Hausdorff for every~$0\leq q\leq k$. Then we have $\C$-linear isomorphisms for all~$p\geq 0$ and~$1\leq q\leq k$ as follows:
    \begin{equation}\label{eq:cohoP1}
        {}^{I\!I}\!\opE_1^{pq}\cong C^p(G,\hH^q(X,V))^G\quad \text{and}\quad {}^{I\!I}\!\opE_2^{pq}\cong\hH_{\ct}^p(G,\hH^q(X,V));
    \end{equation}
and if~$X$ is additionally connected, we also have~${}^{I\!I}\!\opE_2^{p0}\cong\hH_{\ct}^p(G,\rho)$ for all~$p\geq 0$.
\end{lem}
\begin{proof}
This is an application of \Cref{cor:niceDGA}.
The equality in~\eqref{eq:cohoP2a} is clear. Moreover, \eqref{eq:cohoP1} and the statement just after are applications of~\eqref{eq:niceDGA3} and the fact that~$\hH^0(X,V)=V$ if~$X$ is connected.

It remains to prove~\eqref{eq:cohoP2b}. We proceed essentially as in~\cite[p.~596]{DJ:build} except that we drop the hypothesis of finitely many orbits on the $G$-action on the cells of~$X$.

Let~$\sigma\in X(p)$. Since~$G_\sigma$ is open and~$G$ is second countable, the quotient space~$G/G_\sigma$ is a discrete countable $G$-set~\cite[Corollary~2.C.6(1)]{cdlh:metric}. Hence~$I_{G_\sigma}^G(V):=C(G/G_\sigma,V)$ is a Fréchet $G$-module as in \Cref{ex:CXV}. Since~$G/G_\sigma$ is discrete, note that~$I_{G_\sigma}^G(V)=\prod_{g\in G/G_\sigma}V$ as a set, and the compact-open topology on~$I_{G_\sigma}^G(V)$ coincides with the product topology on~$\prod_{g\in G/G_\sigma}V$ induced from~$V$. 
We can therefore write~$I_{G_\sigma}^G(V)=\prod_{g\in G/G_\sigma}V$ unambiguously.

By \Cref{fact:FreG}\eqref{fact:freG3}, the product space~$I:=\prod_{\sigma\in X(p)/G}I_{G_\sigma}^G(V)$ is a Fréchet $G$-module with $G$-action~$G\times I\to I$ determined, for every~$\sigma\in X(p)/G$, by: 
\begin{equation*}
    G\times I\longrightarrow I_{G_\sigma}^G(V)\qquad \big(g,(f_\sigma)_{\sigma\in X(p)/G}\big)\longmapsto g\cdot f_\sigma.
    \end{equation*}
Using the decomposition of~$X(p)$ in~$G$-orbits and the orbit-stabiliser theorem, we obtain the following canonical isomorphisms of Fréchet $G$-modules:
\begin{equation}\label{eq:cohoP5}
    C^p(X,V):=\prod_{\sigma\in X(p)}V\cong \prod_{\sigma\in X(p)/G}\,\prod_{\tau\in G\cdot\sigma}V\cong\prod_{\sigma\in X(p)/G}\,\prod_{g\in G/G_\sigma}V=\prod_{\sigma\in X(p)/G}I_{G_\sigma}^G(V),
\end{equation}
where~$X(p)/G$ denotes a given set of representatives for the $G$-orbits in~$X(p)$.

We now claim to have the following chain of isomorphisms of $\C$-vector spaces:
    \begin{equation}\label{eq:cohoP6}
        {}^{I}\!\opE_1^{pq}\cong\hH_{\ct}^q(G,C^p(X,V))\cong \prod_{\sigma\in X(p)}\hH_{\ct}^q(G,I_{G_\sigma}^G(V))\cong\prod_{\sigma\in X(p)/G}\hH_{\ct}^q(G_\sigma,\rho\vert_{G_\sigma}).
    \end{equation}
    Starting from the left-hand side in~\eqref{eq:cohoP6}, the first isomorphism is due to~\eqref{eq:niceDGA1} and the second one to \Cref{prop:HctProd}. Finally, the third isomorphism in~\eqref{eq:cohoP6} is due to the Eckmann--Shapiro lemma~\cite[\S IX, Proposition~2.3(ii)]{bowal:coho}, which is applicable because the canonical quotient map~$G\twoheadrightarrow G/G_\sigma$ is onto a discrete space and thus admits a continuous section.
\end{proof}

\begin{thm}\label{thm:CohoGX}
Let~$G$ be an lcsc~group with a discrete cellular action on a connected countable CW complex~$X$. Let also~$(\rho,V)$ be a Fréchet $G$-module and assume there are integers~$h\geq 0$ and~$k\geq 1$ such that:
\begin{enumerate}
    \item \label{thm:CohoGX1} $X$ is $h$-acyclic over~$\C$;
    \item \label{thm:CohoGX2} $\hH_{\ct}^{m-\ell}(G_\sigma,\rho\vert_{G_\sigma})=\{0\}$ for all~$\sigma\in X(\ell)$,~$0\leq \ell\leq m-1$ and~$1\leq m\leq k$.
\end{enumerate}
 Then for all~$0\leq m\leq \min\{h,k\}$ we have a $\C$-linear isomorphism
 \begin{equation}\label{eq:CohoGXmain}
     \hH_{\ct}^m(G,\rho)\cong \hH^m(X,\rho).
 \end{equation}
 Moreover, 
 \begin{enumerate}
     \item[(a)] \label{thm:CohoGXa} if~$h=k+1$ and~$\hH_{\ct}^h(G,\rho)=\{0\}$, then~$\hH^h(X,\rho)=\{0\}$; and
     \item[(b)] \label{thm:CohoGXb} if~$k=h+1$ and~$\hH^k(X,\rho)=\{0\}$, then~$\hH_{\ct}^k(G,\rho)=\{0\}$.
 \end{enumerate}
\end{thm}

\begin{proof}
Clearly, $\hH_{\ct}^0(G,\rho)=V^G=\hH^0(X,\rho)$. Thus, we may assume~$h\geq 1$ and prove~\eqref{eq:CohoGXmain} for all~$1\leq m\leq \min\{h,k\}$.

Let~$1\leq m\leq k$. By~\eqref{eq:cohoP2b}, for all~$0\leq \ell\leq m-1$ we have
    \begin{equation*}
        {}^I\!\opE_1^{\ell,m-\ell}\cong \prod_{\sigma\in X(\ell)/G}\hH_{\ct}^{m-\ell}(G_\sigma,\rho\vert_{G_\sigma})=\{0\}
    \end{equation*}
    and hence, passing to subquotients, ${}^I\!\opE_\infty^{\ell,m-\ell}=\{0\}$. Moreover, by \Cref{lem:convSP} and~\eqref{eq:cohoP2a}, 
    \begin{equation}\label{eq:cohoGX1}
{}^I\!\opE_{\infty}^{p0}\cong{}^I\!\opE_{2}^{p0}\cong\hH^p(X,\rho)\qquad\forall\,p\in\{1,\dots, m+1\}.
    \end{equation}
    By \Cref{fact:SPSdouble},
    $${}^I\!\opE_r^{pq}\Rightarrow \hH^{p+q}(\caT),$$
    where~$\caT$ be the total complex associated to the double cochain complex~$\caC(G,\caA_{(X,\rho)})$ (see~\eqref{eq:AXrho}).
    By~\eqref{eq:Hnsubquot}-\eqref{eq:EinftConv} and~\eqref{eq:cohoGX1}, we have
    \begin{equation}\label{eq:cohoGX1'}
\hH^m(\caT)\cong{}^I\!\opE_\infty^{m0}\cong\hH^m(X,\rho)\qquad\forall\,m\in\{1,\dots, k\}
    \end{equation}
    as well as
    \begin{equation}\label{eq:cohoGX2}
        \hH^{k+1}(\caT)\supseteq (\hH^{k+1}(\caT))^{k+1}\cong{}^I\!\opE_\infty^{k+1,0}\cong \hH^{k+1}(X,\rho).
    \end{equation}

   We now argue similarly for~${}^{I\!I}\!\opE_r^{pq}$. Let~$1\leq m\leq h$. Since~$X$ is $h$-acyclic over~$\C$, by \Cref{lem:cohoPlain}, for all~$p\geq 0$ and~$1\leq q\leq m$ we have
   \begin{equation*}
       {}^{I\!I}\!\opE_1^{pq}\cong C^p(G,\hH^q(X,V))^G=\{0\}.
   \end{equation*}
   Therefore, \Cref{lem:convSP} and \Cref{lem:cohoPlain} (recall that~$X$ is connected) imply that 
   \begin{equation}\label{eq:cohoGX3}
       {}^{I\!I}\!\opE_\infty^{p0}\cong {}^{I\!I}\!\opE_2^{p0}\cong\hH_{\ct}^p(G,\rho)\qquad\forall\,p\in\{1,\dots, m+1\}.
   \end{equation}
  Moreover, for every~$0\leq p\leq m-1$ we have
    \begin{equation}\label{eq:vanEinf}
       {}^{I\!I}\!\opE_{\infty}^{p,m-p}=\{0\}
   \end{equation}
   as already~${}^{I\!I}\!\opE_{1}^{p,m-p}$ vanishes. By \Cref{fact:SPSdouble},
   $${}^{I\!I}\!\opE_r^{pq}\Rightarrow \hH^{p+q}(\caT).$$
  Thus, by~\eqref{eq:Hnsubquot}-\eqref{eq:EinftConv} and~\eqref{eq:vanEinf} we conclude that
  \begin{equation}\label{eq:cohoGX4}
\hH^m(\caT)\cong{}^{I\!I}\!\opE_\infty^{m0}\cong\hH_{\ct}^m(G,\rho)\qquad\forall\,m\in\{1,\dots, h\}
  \end{equation}
   and 
   \begin{equation}\label{eq:cohoGX5}
       \hH^{h+1}(\caT)\supseteq (\hH^{h+1}(\caT))^{h+1}\cong{}^{I\!I}\!\opE_\infty^{h+1,0}\cong \hH_{\ct}^{h+1}(G,\rho).
   \end{equation}

   By~\eqref{eq:cohoGX1'} and~\eqref{eq:cohoGX4} we obtain~\eqref{eq:CohoGXmain}. Moreover, from~\eqref{eq:cohoGX4} and~\eqref{eq:cohoGX2} we have part~\hyperref[thm:CohoGXa]{(a)}, and from~\eqref{eq:cohoGX1'} and~\eqref{eq:cohoGX5} we deduce part~\hyperref[thm:CohoGXb]{(b)}.
\end{proof}

The following corollary is a straightforward consequence of \Cref{thm:CohoGX}.
\begin{cor}\label{cor:CohoGX}
Let~$G$ be an lcsc~group with a discrete cellular action on a connected countable CW complex~$X$. Let~$h\geq 0$ and~$k\geq 1$ be integers such that:
    \begin{enumerate}
        \item \label{cor:CohoGX1a} $X$ is $h$-acyclic over~$\C$;
        \item \label{cor:CohoGX1b} $G_\sigma$ has property~$[T_{k-\dim(\sigma)}]$ for all~$\sigma\in \Ce(X)$ such that~$0\leq \dim(\sigma)\leq k-1$.
    \end{enumerate}
    Then for every unitary representation~$(\rho,V)$ of~$G$ we have
    \begin{equation*}
        \hH_{\ct}^m(G,\rho)=\hH^m(X,\rho)\qquad\forall\,m\in\big\{0,\dots, \min\{h,k\}\big\}. 
    \end{equation*}
    Moreover, if $X$ is $(k+1)$-acyclic over~$\C$ then 
    \begin{equation*}
        \hH_{\ct}^{k+1}(G,\rho)=\{0\}\,\Longrightarrow\,\hH^{k+1}(X,\rho)=\{0\},
    \end{equation*}
   and if~$X$ is $(k-1)$-acyclic over~$\C$ then
   \begin{equation*}
       \hH^k(X,\rho)=\{0\}\,\Longrightarrow\,\hH_{\ct}^k(G,\rho)=\{0\}.
   \end{equation*}
\end{cor}

We now combine \Cref{cor:CohoGX} and the vanishing results for~$\hH^\bullet(X,\rho)$ recalled in \Cref{sus:twistCoho} to obtain criteria for properties~$\Tk$.

Firstly, combining \Cref{cor:CohoGX} and \Cref{thm:OppVan}, we deduce the following generalisation of~\cite[Corollary~1.2]{Op15:vanish}.
\begin{cor}
Let~$X$ be a locally finite pure $n$-dimensional $(n\geq 1)$ simplicial complex which is $h$-acyclic over~$\C$ $(h\geq 0)$, all its positive dimensional links are connected, and such that~\eqref{eq:OppVan} holds. 
Let also~$G$ be a unimodular lcsc~group with a proper discrete simplicial action on~$X$ with finitely many orbits.

Then~$G$ has property~$[T_{\min\{h,n-1\}}]$. If additionally~$h=n-2$, then~$G$ has property~$[T_{n-1}]$.
\end{cor}
\begin{proof}
Let~$(\rho,V)$ be any unitary representation of~$G$.
By \Cref{thm:OppVan}, $\hH^\ell(X,\rho)=\{0\}$ for all~$1\leq \ell\leq n-1$. Moreover, for every non-empty simplex~$\sigma$, $G_\sigma$ is compact and hence has property~$[T_{k}]$ for every~$k\geq 1$ (see \Cref{thm:Tkperm}\eqref{thm:Tkperm1}). Now \Cref{cor:CohoGX} applies with~$k=n-1$. 
\end{proof}

\begin{prop}\label{prop:cohoPlainCC}
    Let~$(G,X,\Delta,\gamma)$ be in~\hyperref[Bn]{$\caB_n$} for some~$n\geq 1$, and let~$h\geq 0$ and~$k\geq 1$ be integers such that:
    \begin{enumerate}
        \item \label{prop:cohoPlainCC1} $X$ is $h$-acyclic over~$\C$;
        \item \label{prop:cohoPlainCC2} $G_\sigma$ has property~$[T_{k-\dim(\sigma)}]$ for all faces~$\sigma\subseteq \gamma$ with~$0\leq \dim(\sigma)\leq k-1$;
        \item \label{prop:cohoPlainCC3} for every unitary representation~$(\rho,V)$, the cosine matrix~$A(V_0,\dots, V_n)$ associated to~$(\rho,V)$ and~$\Delta$ is positive definite (see \Cref{defn:CosMatRho}).
    \end{enumerate}
    Then~$G$ has property~$[T_{\min\{h,k, n-1\}}]$.

   Assume in particular that~$k=n-1$ and that~$X$ is acyclic\footnote{Or equivalently, since~$X$ is $n$-dimensional, $X$ is $n$-acyclic over~$\C$.} over~$\C$. Then~$\hH_{\ct}^n(G,\rho)\neq\{0\}$ for every unitary representation~$(\rho,V)$ of~$G$ such that~$\sum_{\tau\subsetneq \gamma}V^{G_\tau}$ has non-zero orthogonal complement in~$V^{G_\gamma}$.
\end{prop}
\begin{proof}
    \Cref{cor:CohoGX} and \Cref{cor:GRO1} imply the first part of the statement.

    Assume additionally that~$k=n-1$ and that~$X$ acyclic over~$\C$. Let~$(\rho,V)$ be any unitary representation of~$G$ such that~$\sum_{\tau\subsetneq \gamma}V^{G_\tau}$ has non-zero orthogonal complement in~$V^{G_\gamma}$. By \Cref{cor:GRO1}, that orthogonal complement is exactly~$\hH^n(X,\rho)$. Thus, by \Cref{cor:CohoGX} we have~$\hH^n_{\ct}(G,\rho)\neq\{0\}$.
\end{proof}

By \Cref{prop:GRO2}\eqref{prop:GRO2b} and \Cref{prop:cohoPlainCC} we deduce the following.
\begin{cor}\label{cor:cohoPlainCCa}
    Let~$(G,X,\Delta,\gamma)$ be in~\hyperref[Bn]{$\caB_n$} for some~$n\geq 1$, and assume that~$G$ acts properly on~$X(n)$ and that all $1$-dimensional links of~$X$ are finite. Let~$h\geq 0$ and~$k\geq 1$ be integers such that conditions~\eqref{prop:cohoPlainCC1}-\eqref{prop:cohoPlainCC2} in \Cref{prop:cohoPlainCC} hold. Assume furthermore that:
    \begin{itemize}
        \item[(iii')]\label{cor:cohoPlainCC3'} the cosine matrix~$A(X)$ of~$X$ is positive definite (see \Cref{defn:CosMatX}).
    \end{itemize}
    Then~$G$ has property~$[T_{\min\{h,k,n-1\}}]$.

   Let in particular~$n=k+1$ and let~$X$ be acyclic over~$\C$. Then~$\hH_{\ct}^n(G,\rho)\neq\{0\}$ for every unitary representation~$(\rho,V)$ of~$G$ such that~$\sum_{\tau\subsetneq \gamma}V^{G_\tau}$ has non-zero orthogonal complement in~$V^{G_\gamma}$.
\end{cor}

For groups having BN-pairs, \Cref{prop:cohoPlainCC} specialises as follows:
\begin{cor}\label{cor:cohoPlainBN}
    Let~$G$ be an lcsc~group with a BN-pair~$(B,N)$ of type~$(W,S)$, where~$S=\{0,\dots, n\}$ has~$n+1$ elements and~$B\leq G$ is open. Assume that, for every unitary representation~$(\rho,V)$ of~$G$, the cosine matrix~$A(V^{P_0},\dots, V^{P_n})$ is positive definite (here~$P_i$ is the standard parabolic subgroup of type~$\{i\}$). Let~$k\geq 1$ be an integer such that, for every~$J\subseteq S$ with~$|J|>n-k$, $P_{S\setminus J}$ has property~$[T_{k+|J|-n}]$. 
   Then:
   \begin{enumerate}
       \item $G$ has property~$[T_{\min\{k,n-1\}}]$;
       \item Let~$n=k+1$,~$(W,S)$ be non-spherical, and denote by~$\gamma$ the fundamental chamber of~$\Delta(G,B)$. Then~$\hH_{\ct}^n(G,\rho)\neq\{0\}$ for every unitary representation~$(\rho,V)$ of~$G$ such that~$\sum_{\tau\subsetneq \gamma}V^{G_\tau}$ has non-zero orthogonal complement in~$V^{G_\gamma}$.
   \end{enumerate} 
\end{cor}
\begin{proof}
    Let~$X=\Delta(G,B)$ be the simplicial building associated to the BN-pair~$(B,N)$. Denote by~$\gamma=\{1P_{S\setminus\{s\}}\mid s\in S\}$ its fundamental chamber, and let~$\Delta$ be the subcomplex of~$X$ of all faces of~$\gamma$. Then~$(G,X,\Delta,\gamma)\in \caB_n$ where~$n=|S|-1$ (see \Cref{ex:BNBn}). By Solomon--Tits theorem, $X$ is $(n-1)$-spherical and hence $(n-1)$-acyclic over~$\C$. More precisely, if~$(W,S)$ is non-spherical then~$X$ is contractible and hence acyclic over~$\C$. Moreover, for every~$0\leq \ell\leq n$, every $\ell$-dimensional simplex has the form~$\{gP_{S\setminus\{s\}}\mid s\in J\}$ for some~$g\in G$ and~$J\subseteq S$ with~$|J|=\ell+1$, and thus its stabiliser is conjugated to~$P_{S\setminus J}$ (see~\eqref{eq:intersParab}). Now \Cref{prop:cohoPlainCC} applies.
\end{proof}

Another variation of \Cref{thm:CohoGX} and its corollaries is the following generalisation of~\cite[Theorem~E]{DJ:build} and~\cite[Theorem~3.9]{gro:curv}. 

\begin{thm}\label{thm:CohoGXPlain3}
    Let~$G$ be an lcsc~group with a discrete simplicial action on a countable pure $n$-dimensional simplicial complex~$X$ for some integer~$n\geq 1$. Let~$k\geq 1$ be an integer such that:
    \begin{enumerate}
        \item for every~$\gamma\in X(n)$, $G_\gamma$ has property~$[T_k]$;
        \item for every (possibly empty) simplex~$\sigma$ of~$X$, given~$d=\dim(X_\sigma)$, the following hold:
        \begin{enumerate}
            \item if~$d\leq k$ then~$X_\sigma$ is finite; 
            \item if~$d>k$ then~$X_\sigma$ is finite or $(k-1)$-acyclic over~$\C$, and $\hH^\ell(X_\sigma,\rho_\sigma)=\{0\}$ for every unitary representation~$(\rho_\sigma,V_\sigma)$ of~$G_\sigma$ and every~$1\leq \ell\leq k$.
        \end{enumerate}
    \end{enumerate}
    Then~$G$ has property~$\Tk$.
\end{thm}
In \Cref{thm:CohoGXPlain3}, the restriction to simplicial complexes (rather than arbitrary CW complexes) guarantees a well-defined notion of links at every simplex.

\begin{proof}
Let~$k\geq 1$. We firstly prove the claim whenever~$n=\dim(X)\leq k$. In this case~$X$ is finite. Pick~$\gamma\in X(n)$ and, by the orbit-stabiliser theorem, observe that
    \begin{equation*}
        |G:G_\gamma|=|G\cdot \gamma|\leq |X(n)|<\infty.
    \end{equation*}
 Since~$G_\gamma$ has property~$\Tk$, by \Cref{thm:Tkperm}\eqref{thm:Tkperm1} so has~$G$.

We now prove the claim if~$\dim(X)\geq k$ proceeding by induction on~$\dim(X)$. The base case~$\dim(X)=k$ is already proved. Let now~$n>k$, assume that the claim is true in dimension~$<n$ and let~$\dim(X)=n$.

Let~$\sigma$ be a non-empty simplex of~$X$. Notice that~$G_\sigma$ is an open subgroup of~$G$ and hence an lcsc group, and clearly the link~$X_\sigma$ is a countable pure~$(n-1-\dim(\sigma))$-simplicial complex. Note also that~$(X_\sigma)_{\tau}=X_{\sigma\cup\tau}$ for every simplex~$\tau$ of~$X_\sigma$, and the stabiliser of every top dimensional simplex of~$X_\sigma$ is equal to~$G_\gamma$ for some~$\gamma\in X(n)$ with~$\gamma\supseteq \sigma$. By induction,~$G_\sigma$ has property~$\Tk$. By \Cref{cor:CohoGX}, also~$G$ has property~$\Tk$.
\end{proof}

\begin{prop}\label{prop:vanSpecCC}
Let~$(G,X,\Delta,\gamma)\in \hyperref[Bn]{\caB_n}$ for some~$n\geq 1$ and denote by~$\tau\colon X\to 2^S$ the type function of~$X$. Assume there is an integer~$k\geq 1$ such that:
\begin{enumerate}
    \item \label{prop:vanSpecCC1} $G_\gamma$ has property~$\Tk$; 
    \item \label{prop:vanSpecCC2} For all~$J\subseteq S$, the following hold:
    \begin{enumerate}
        \item \label{prop:vanSpecCC2a} if~$|J|\geq n-k$ then~$X_{\gamma_J}$ is finite;
        \item \label{prop:vanSpecCC2b} if~$|J|<n-k$ then~$X_{\gamma_J}$ is finite or $(k-1)$-acyclic over~$\C$. Moreover, if~$S\setminus J=\{j_0,\dots, j_\ell\}$, for every unitary representation~$(\rho,V)$ of~$G_{\gamma_{J}}$ the cosine matrix~$A(V_{j_0},\dots, V_{j_\ell})$ is positive definite, where~$V_{j_i}:=V^{G_{\gamma_{S\setminus\{j_i\}}}}$.
    \end{enumerate}
\end{enumerate}
Then~$G$ has property~$\Tk$. 

If in particular~$k=n-1$ and~$X$ is acyclic over~$\C$, then~$\hH_{\ct}^n(G,\rho)\neq\{0\}$ for every unitary representation~$(\rho,V)$ of~$G$ such that~$\sum_{\tau\subsetneq \gamma}V^{G_\tau}$ has non-zero orthogonal complement in~$V^{G_\gamma}$.
\end{prop}
\begin{proof}
    By \Cref{lem:Fams}\eqref{lem:Fams3}, for every face~$\sigma\subseteq \gamma$ the tuple~$(G_\sigma,X_\sigma, \Delta_\sigma,\sigma)$ belongs to~$\caB_{n-1-\dim(\sigma)}$. More precisely,~$\sigma=\gamma_J$ where~$J=\tau(\sigma)$. Hence~$X_\sigma=X_J$ and~$\dim(X_\sigma)=\dim(X_J)=n-1-\dim(\sigma)=n-|J|$.
    
    By \Cref{cor:GRO1} and the hypothesis~\eqref{prop:vanSpecCC2b}, whenever~$|J|<n-k$ (i.e, $\dim(X_\sigma)>k)$ we have~$\hH^\ell(X_\sigma,\rho_\sigma)=\{0\}$ for all~$1\leq \ell\leq k$ (in fact, for all~$1\leq \ell\leq \dim(X_\sigma)-1$). Now \Cref{thm:CohoGXPlain3} applies and~$G$ has property~$\Tk$.
    
    In fact, for every~$J\subseteq S$, applying the first part of the statement to the tuple~$(G_{\gamma_J}, X_{\gamma_J}, \Delta_{\gamma_J}, \gamma_J)$ we have that~$G_{\gamma_J}$ has property~$\Tk$. If~$S=\{0,\dots, n\}$, the hypothesis~\eqref{prop:vanSpecCC2b} with~$J=\emptyset$ implies that~$A(V_0,\dots, V_n)$ is positive definite for every unitary representation~$(\rho,V)$ of~$G$. The second part of the statement now follows from the second part of \Cref{prop:cohoPlainCC}.
\end{proof}

\begin{cor}\label{cor:vanishBN}
    Let~$G$ be an lcsc~group with a BN-pair~$(B,N)$ of type~$(W,S)$, of finite thickness~$(q_s)_{s\in S}$ and with~$B\leq G$ open. 
    {For every subset~$J\subseteq S$ with~$|J|<|S|-1-k$, given~$J=\{j_0,\ldots, j_\ell\}$, assume that for every unitary representation~$(\rho,V)$ of~$P_J$ the cosine matrix~$A(V^{P_{j_0}},\dots, V^{P_{j_\ell}})$ is positive definite.}
 Let also~$1\leq k\leq |S|-1$ be an integer such that~$B$ has property~$\Tk$. Then:
    \begin{enumerate}
        \item \label{cor:vanishBN1} If~$(W,S)$ is $(k+1)$-spherical, then~$G$ has property~$[T_k]$.
        \item \label{cor:vanishBN2} Let~$|S|=n+1$, assume~$(W,S)$ is $n$-spherical but not spherical, and let~$\min\{q_s\mid s\in S\}>|S|$. Then~$G$ has property~$[T_{n-1}]$ but not property~$(T_{n})$.
    \end{enumerate}
\end{cor}
\begin{proof}
Our goal is to apply \Cref{prop:vanSpecCC} to the tuple~$(G,\Delta(G,B),\Delta,\gamma)$, which belongs to~$\caB_{|S|-1}$ by \Cref{ex:BNBn}. 
For every~$J\subseteq S$, recall that the simplex~$\gamma_J=\{1P_{S\setminus \{j\}}\mid j\in J\}$ has stabiliser~$P_{S\setminus J}$ and its link in~$\Delta(G,B)$ is $P_{S\setminus J}$-equivariantly isomorphic to~$\Delta(P_{S\setminus J}, S\setminus J)$ and has dimension~$|S|-1-|J|$. In particular, $\Delta(P_{S\setminus J},J)$ is finite whenever~$|S\setminus J|\leq k+1$ because it is a spherical building of finite thickness. In any case, for every~$I\subseteq S$, $\Delta(P_{I},I)$ is a building of type~$(W_I,I)$ of finite thickness, and it is finite if~$(W_I,I)$ is spherical or, by Solomon--Tits theorem, contractible if~$(W_I,I)$ is not spherical.
Now \Cref{prop:vanSpecCC} and \Cref{ex:H2perp} imply the statement.
\end{proof}

\begin{cor}\label{cor:vanishBNco}
    Let~$G$ be an lcsc~group with a BN-pair~$(B,N)$ of $2$-spherical type~$(W,S)$, finite thickness vector~$(q_s)_{s\in S}$ and with~$B\leq G$ compact open. Assume also that
    \begin{equation}\label{eq:niceqineq}
        \frac{\qmin+1}{2\sqrt{\qmax}}>1-\veps_{(W,S)},
    \end{equation}
    where~$\qmin:=\min\{q_s\mid s\in S\}$, $\qmax:=\max\{q_s\mid s\in S\}$ and~$\veps_{(W,S)}$ is the smallest eigenvalue of the cosine matrix of~$(W,S)$. Then the following holds for every integer~$1\leq k\leq |S|-1$:
    \begin{enumerate}
        \item If~$(W,S)$ is $(k+1)$-spherical then~$G$ has property~$\Tk$;
        \item Let~$|S|=n+1$, assume that~$(W,S)$ is~$n$-spherical and not spherical and let~$\qmin>|S|$. Then~$G$ has property~$[T_{n-1}]$ but not property~$(T_n)$.
    \end{enumerate}
\end{cor}
\begin{proof}
    We aim to apply \Cref{cor:vanishBN}. By \Cref{thm:Tkperm}\eqref{thm:Tkperm1}, $B$ has property~$\Tk$. Moreover, for every non-empty subset~$J=\{j_0,\dots, j_\ell\}$ of~$S$, the BN-pair~$(B,N_J)$ for~$P_J$ induced by~$(B,N)$ has $2$-spherical type~$(W_J,J)$. This means that the $P_J$-action on~$X_J:=\Delta(P_J,B)$ satisfies the hypotheses of \Cref{prop:GRO2}. Hence, by that corollary, if the cosine matrix~$A(X_J)$ of~$X_J$ is positive definite then~$A(V_{j_0},\dots, V_{j_\ell})$ positive definite for every unitary representation~$(\rho,V)$ of~$P_J$. It remains to prove that~$A(X_J)$ is positive definite for every non-empty~$J\subseteq S$. Note that~$A(X_J)$ is a principal submatrix of the cosine matrix~$A(X)$, where~$X=\Delta(G,B)$. By Sylvester theorem, if~$A(X)$ is positive definite so is every~$A(X_J)$. As argued in \Cref{ex:cosMatBuil}, the condition in~\eqref{eq:niceqineq} guarantees that~$A(X)$ is indeed positive definite. 

    The claim now follows from \Cref{cor:vanishBN}.
\end{proof}
\begin{rem}
    \begin{enumerate}
        \item \Cref{cor:vanishBNco} recovers~\cite[Theorem~4.6]{gro:curv} and corrects a minor typo coming from~\cite[Lemma~4.5]{gro:curv} (see \Cref{ex:cosMatBuil}). In~\cite[Theorem~4.6]{gro:curv}, given~$\wtil{\mu}=\veps_{(W,S)}$, the condition $\wtil{\mu}>1-\frac{\qmin+1}{2\sqrt{\qmin}}$ should be replaced by~\eqref{eq:niceqineq}. 
        \item If~$\qmin=\qmax\geq 2$, the condition~\eqref{eq:niceqineq} is empty if~$(W,S)$ is Euclidean. Indeed, in this case, the cosine matrix of~$(W,S)$ is positive semidefinite~\cite[\S 6.4--6.5]{hum:ref} and hence~$\frac{\qmin+1}{2\sqrt{\qmin}}>1\geq 1-\veps_{(W,S)}$.
    \end{enumerate}
\end{rem}

\section{Root group data and twin building lattices}\label{s:RGD}

\subsection{Roots systems}\label{sus:rootSys}

In this section, we give a brief introduction to root systems following the combinatorial approach of~\cite[\S 5.1]{tits:KM}. 

Let~$(W,S)$ be a Coxeter system and~$\ell=\ell_{(W,S)}$ be its Coxeter length. For every~$s\in S$, define
\begin{equation*}
    \alpha_s:=\{w\in W\mid \ell(sw)>\ell(w)\}.
\end{equation*}
The \emph{root system} associated to~$(W,S)$ is the collection of all roots of~$(W,S)$, that is, 
\begin{equation*}
    \Phi=\Phi(W,S):=\{w\alpha_s\mid w\in W,\,s\in S\}.
\end{equation*}
Let also
\begin{equation*}
  \Phi^+:=\{\alpha\in \Phi\mid 1\in \alpha\}\quad\text{and}\quad \Phi^-:=\{\alpha\in \Phi\mid 1\not\in\alpha\}  
\end{equation*}
be the sets of all \emph{positive} and \emph{negative roots} of~$(W,S)$. Clearly, $\Phi=\Phi^+\sqcup \Phi^-$.

For all~$w\in W$ and~$s\in S$, notice that
\begin{equation*}
    W\setminus w\alpha_s=ws\alpha_s\quad\text{and}\quad 1\in w\alpha_s\,\Longleftrightarrow\,1\not\in ws\alpha_s.
\end{equation*}
The complement operation in~$W$ defines an involution~$-\colon \Phi\to \Phi$, $\alpha\mapsto -\alpha:=W\setminus \alpha$ that satisfies~$-\Phi^+=\Phi^-$.

A pair of roots~$\{\alpha,\beta\}$ is said to be \emph{prenilpotent} if~$\alpha\cap \beta\neq\emptyset$ and~$(-\alpha)\cap (-\beta)\neq \emptyset$.
For every such pair of roots, the following set is finite:
\begin{equation*}
    (\alpha,\beta):=\big\{\gamma\in \Phi\,\big\vert\, \gamma\not\in\{\alpha,\beta\},\,\alpha\cap \beta\subseteq \gamma\text{ and }(-\alpha)\cap (-\beta)\subseteq -\gamma\big\}.
\end{equation*}
\begin{ex}\label{ex:openInt}
    For every~$s,t\in S$ with~$s\neq t$, note that~$\{\alpha_s,\alpha_t\}$ is a prenilpotent pair because~$t\in \alpha_s\cap (-\alpha_t)$ and~$s\in (-\alpha_s)\cap \alpha_t$. In this case, it is immediate to check that~$(\alpha_s,-\alpha_t)=\emptyset$.
\end{ex}

For every~$J\subseteq S$, consider the following subsets of~$\Phi$:
\begin{equation}\label{eq:PhiJ}
\Phi_J:=\{w\alpha_s\mid w\in W_J,\,s\in S\}\quad\text{and}\quad \Phi_J^\pm:=\Phi_J\cap \Phi^\pm.
\end{equation}
\begin{rem}[Identification of~$\Phi_J$ and~$\Phi(W_J,J)$]\label{rem:PhiJ}
    For every~$J\subseteq S$, the subset~$\Phi_J$ has a canonical bijection onto~$\Phi(W_J,J):=\{w\alpha_s^J\mid w\in W_J,\,s\in J\}$, where~$\alpha_s^J:=\{u\in W_J\mid \ell_{(W_J,J)}(su)>\ell_{(W_J,J)}(u)\}$ for every~$s\in J$ and~$\ell_{(W_J,J)}$ is the Coxeter length on~$(W_J,J)$. Indeed, since~$\ell\vert_{W_J}=\ell_{(W_J,J)}$ by~\cite[Theorem~5.5(b)]{hum:ref}, we have~$\alpha_s^J=\alpha_s\cap W_J$ for every~$s\in S$. Therefore, 
    \begin{equation*}
        f\colon \Phi_J\longrightarrow \Phi(W_J,J)\qquad w\alpha_s,\longmapsto w\alpha_s^J\quad\forall\,w\in W_J,\,s\in J
    \end{equation*}
    is a bijection. Notice also that~$f(\Phi_J^\pm)=\Phi(W_J,J)^\pm$. This allows us to identify the root system~$\Phi(W_J,J)$ with~$\Phi_J$.
\end{rem}

Examples of Coxeter groups (and thus of root systems) can be produced from a generalised Cartan matrix, as recalled in the following example.
\begin{ex}\label{ex:genCart}
   Let~$S$ be a finite non-empty set. A \emph{generalised Cartan matrix} over~$S$ is any matrix~$A=[A_{st}]_{s,t\in S}$ with integral entries satisfying the following properties:
   \begin{enumerate}
    \item $A_{ss}=2$ for every~$s\in S$;
    \item $A_{st}\leq 0$ for all~$s\neq t$ in~$S$; and
    \item $A_{st}=0$ if and only if~$A_{ts}=0$.
\end{enumerate}  
As the name suggests, the above notion extends the one of \emph{Cartan matrix}, which is precisely a generalised Cartan matrices that factors as a product~$DX$ of a real diagonal matrix~$D$ and a real square matrix~$X$ which is positive definite and symmetric.

Every generalised Cartan matrix~$A$ has an associated Coxeter matrix~$M_A=[m_{st}]_{s,t\in S}$ whose entries~$m_{st}$,~$s\neq t$, are determined by~$A_{st}A_{ts}$ as follows:
\begin{equation*}
    \begin{tabular}{c|ccccc}
        $A_{st}A_{st}$ & $0$ & $1$ & $2$ & $3$ & $\geq 4$  \\\hline
        $s_{st}$ & $2$ & $3$ & $4$ & $6$ & $+\infty$
    \end{tabular}
\end{equation*}
see~\cite[Proposition~4.22]{marq:km}. We call~$(W_A,S)$ the Coxeter group associated to~$M_A$. For future reference, if~$(W_A,S)$ is $2$-spherical we notice that
\begin{equation}\label{eq:GCMcos}
    A_{st}A_{ts}=4\cos^2\Big(\frac{\pi}{m_{st}}\Big)\leq 3\qquad\forall\,s,t\in S,\,s\neq t.
\end{equation}
We also denote by~$\Phi_A:=\Phi(W_A,S)$ the root system associated to~$(W_A,S)$.
\end{ex}

\subsection{Root group data}\label{sus:RGD}
Let~$(W,S)$ be a Coxeter system and let~$\Phi:=\Phi(W,S)$ be its associated root system. Let also~$\Phi^\pm$ be the set of all positive/negative roots.

\begin{defn}[\protect{\cite[\S 2.1]{CapRem:root}}]\label{defn:RGD}
    A \emph{root group datum} of type~$(W,S)$ for a group~$G$ is a family of subgroups~$(U_\alpha)_{\alpha\in \Phi}$, called \emph{root groups}, that satisfies the following axioms:
\begin{itemize}
    \item[(RGD0)] \label{RGD0} $U_\alpha\neq\{1\}$ for every~$\alpha\in \Phi$, and~$G=\langle U_\alpha\mid\alpha\in \Phi\rangle$; 
     \item[(RGD1)] \label{RGD1} $U_{-\alpha_s}\not\subseteq \langle U_\alpha\mid \alpha\in \Phi^+\rangle$ and~$U_{\alpha_s}\not\subseteq \langle U_\alpha\mid \alpha\in \Phi^-\rangle$ for every~$s\in S$;
    \item[(RGD2)] \label{RGD2} For every~$s\in S$, there is a function~$\mu\colon U_{\alpha_s}\setminus\{1\}\to G$ such that, for all~$u\in U_{\alpha_s}\setminus\{1\}$ and~$\alpha\in \Phi$,
    \begin{equation*}
        \mu(u)\in U_{-\alpha_s}uU_{-\alpha_s}\quad\text{and}\quad \mu(u)U_\alpha \mu(u)^{-1}=U_{s\alpha};
    \end{equation*}
     \item[(RGD3)] \label{RGD3} For every prenilpotent pair~$\{\alpha,\beta\}\subseteq \Phi$, the commutator subgroup~$[U_\alpha,U_\beta]$ is contained in~$U_{(\alpha,\beta)}:=\langle U_\gamma\mid \gamma\in (\alpha,\beta)\rangle$.
\end{itemize}
For future reference, we set the following subgroups of~$G$:
\begin{equation*}
    U_\pm:=\langle U_\alpha\mid \alpha\in \Phi^\pm\rangle.
\end{equation*}
A root group datum is said to be \emph{locally finite} if~$U_\alpha$ is finite for every~$\alpha\in \Phi$.
\end{defn}
\begin{rem}\label{rem:RGDrem}
    In contrast to~\hyperref[RGD1]{(RGD1)}, for every~$s\in S$ we have~$U_{\alpha_s}\subseteq U_+$ and~$U_{-\alpha_s}=U_{s\alpha_s}\subseteq U_-$. Moreover, since~$-\alpha_s=s\alpha_s$, by~\hyperref[RGD2]{(RGD2)} we have~$|U_{\alpha_s}|=|U_{-\alpha_s}|$ for every~$s\in S$.
\end{rem}

In some definitions (see~\cite[\S 3.3]{Tits:twinKM} or~\cite[\S 8.6.1]{ab:build}), the axiom~\hyperref[RGD0]{(RGD0)} is replaced by the following:
\begin{itemize}
    \item[(RGD0')] \label{RGD0'} $U_{\alpha}\neq\{1\}$ for all~$\alpha\in \Phi$, and there is a subgroup~$T\leq G$ that normalises each~$U_\alpha$ such that~$G=T\langle U_+,U_-\rangle$. 
\end{itemize}
The discrepancy between the two definitions is negligible for the purpose of this paper, as observed in the following lemma:
\begin{fact}[\protect{\cite[Exercise~7.130]{ab:build}, \cite[Remarks~2.1(5)]{CapRem:root}}]\label{fact:defnRGD}
    Let~$\wtil{G}$ be a group and~$(U_\alpha)_{\alpha\in \Phi}$ be a family of subgroups satisfying~\hyperref[RGD0']{(RGD0')} and the axioms from~\hyperref[RGD1]{(RGD1)} to~\hyperref[RGD3]{(RGD3)}. Then~$(U_\alpha)_{\alpha\in \Phi}$ is a root group datum of type~$(W,S)$ as in \Cref{defn:RGD} for the subgroup~$G:=\langle U_\alpha\mid\alpha\in \Phi\rangle$.

    Conversely, let~$(U_\alpha)_{\alpha\in \Phi}$ be a root group datum of type~$(W,S)$ for a group~$G$ in the sense of \Cref{defn:RGD}. Then, for every group of the form~$\wtil{G}:=TG$ for some subgroup~$T\leq \wtil{G}$ normalising each~$U_\alpha$,~$(\wtil{G}, (U_\alpha)_{\alpha\in \Phi})$ satisfies~\hyperref[RGD0']{(RGD0')} and the axioms from~\hyperref[RGD1]{(RGD1)} to~\hyperref[RGD3]{(RGD3)}.
\end{fact}

The definition of root group datum that we chose follows the one given in~\cite[\S 2.1]{CapRem:root}, and it will allow us to slightly simplify the formulation of certain results (e.g., \Cref{thm:TBLsimp}).

\begin{rem}\label{rem:defnRGDsym}
 As observed in~\cite[Remarks~2.1(3)]{CapRem:root} and in~\cite[Remarks~8.49(b)]{ab:build}, $(U_\alpha)_{\alpha\in \Phi}$ is a root group datum of type~$(W,S)$ for~$G$ if and only if so is~$(U_{-\alpha})_{\alpha\in \Phi}$. This determines a \textbf{complete symmetry between the "$+$"-case and the "$-$"-case}. 
To ease the notation, in this paper we often formulate the results only for one sign.
\end{rem}

Typical examples of groups with root group data arise from the group~$\caG(K)$ of $K$-points of a linear algebraic $K$-group~$\caG$, where~$K$ is a field (e.g., $\caG=\SL_n$). These root group data have spherical type. We refer the reader to~\cite[\S 7.9.1-7.9.2-7.9.3]{ab:build} for details or to~\cite[\S 2.6]{CapRem:root} for a brief overview. It is also worth mentioning that certain $S$-arithmetic groups, for instance~$\SL_n(K[t,t^{-1}])$, admit a canonical choice of root group datum of affine type; see~\cite[\S 2.7]{CapRem:root}. 
The main source of examples of groups with root group datum of non-spherical and possibly non-Euclidean type, arises from Kac--Moody groups over fields, as briefly recalled in the following example.
\begin{ex}[Kac--Moody and Kac--Moody--Steinberg groups]\label{ex:KMgroups}
   Let~$A=[A_{st}]_{s,t\in S}$ be a generalised Cartan matrix with associated Coxeter group~$(W_A,S)$, and let~$\Phi_A$ be its associated root system (see \Cref{ex:genCart}). Tits~\cite{tits:KM} introduced a functor~$\caG_A$ from the category of commutative rings to the category of groups with the following property: for every field~$K$, the group~$\caG_A(K)$ is naturally endowed with a family of subgroups~$(U_\alpha)_{\alpha\in \Phi_A}$, each of them isomorphic to~$(K,+)$, that forms a root group datum of type~$(W_A,S)$ (in the sense of \Cref{defn:RGD}) for the subgroup~$\caG_A(K)^\dagger:=\langle U_\alpha\mid \alpha\in \Phi_A\rangle\leq \caG_A(K)$. 
   Since~$\caG_A(K)=T(K)\cdot\caG_A(K)^\dagger$ and~$T(K)\cong (K^\times)^{|S|}$, if~$K$ is finite then~$\caG_A(K)^\dagger$ has finite index in~$\caG_A(K).$
   We call the group~$\caG_A(K)$ the \emph{Kac--Moody group of type~$A$ over~$K$}.

   Alongside Kac--Moody groups are {Kac--Moody--Steinberg groups} over a field~$K$. We follow~\cite[\S 3.1]{GdP-VP:KMS} and refer to that article for further background. For every~$J\subseteq S$ spherical, let~$U_J$ be the subgroup of~$\caG_A(K)$ generated by~$\{U_{\alpha_s}\}_{s\in J}$. The groups~$U_J$, where~$J\subseteq S$ is spherical, form a direct system with respect to the inclusion maps~$U_I\hookrightarrow U_J$,~$I\subseteq J$. The direct limit of that direct system, i.e.,
   \begin{equation*}
       \caU_A(K):=\bigast_{J\subseteq S\text{ sph.}}U_J\Big/(U_I\hookrightarrow U_J,\,I\subseteq J,\,J\text{ sph.}),
   \end{equation*}
   is called the \emph{Kac--Moody--Steinberg group of type~$A$ over~$K$}.
\end{ex}

\begin{rem}[Not only Kac--Moody groups]\label{rem:notKM}
    The family of groups with root group data is strictly larger than Kac--Moody groups, even if we ask that the type is $2$-spherical and the root groups are finite (a typical assumption that we will assume in our results). We refer to~\cite[\S 1.1.III]{caprem:KMsup} and the references therein for details. 
\end{rem}

\begin{notat}\label{notat:RGD}
Throughout, we adopt the following notation:
    \begin{enumerate}
        \item For every~$s\in S$, let~$U_s=U_s^+:=U_{\alpha_s}$ and~$U_s^-:=U_{-\alpha_s}$.
        \item For every~$J\subseteq S$, let~$U_J=U_J^+:=\langle U_s\mid s\in J\rangle$ and~$U_J^-:=\langle U_s^{-}\mid s\in J\rangle$. Conventionally, $U_\emptyset^\pm=\{1\}$.
        \item For every~$s\in S$, let~$q_s:=|U_s|=|U_s^-|$ (see \Cref{rem:RGDrem} for the latter equality). If~$q_s<\infty$ for every~$s\in S$, define also
         \begin{equation*}
    \qmin:=\min\{q_s\mid s\in S\}\quad\text{and}\quad \qmax:=\max\{q_s\mid s\in S\},
\end{equation*}
and for every~$J\subseteq S$ let
 \begin{equation*}
    \qminJ:=\min\{q_s\mid s\in J\}\quad\text{and}\quad \qmaxJ:=\max\{q_s\mid s\in J\}.
\end{equation*}
    \end{enumerate}
\end{notat}

From now on, we set a group~$G$ with a root group datum~$(U_\alpha)_{\alpha\in \Phi}$ of type~$(W,S)$. 

\smallskip

\noindent Recall that, for every~$J\subseteq S$, the set~$\Phi_J$ as in~\eqref{eq:PhiJ} can be identified with the root system~$\Phi(W_J,J)$ (see \Cref{rem:PhiJ}). The following lemma is a specialisation of~\cite[Lemma~2.2]{CapRem:root} to the only case of interest for this paper. See also~\cite[Exercise~7.126(a)]{ab:build}.
\begin{lem}[\protect{\cite[Lemma~2.2, Remark~2.3]{CapRem:root}}]\label{lem:RGDsubsyst}
    Given~$J\subseteq S$, $(U_\alpha)_{\alpha\in \Phi_J}$ is a root group datum of type~$(W_J,J)$ for the subgroup~$L_J:=\langle U_\alpha\mid \alpha\in \Phi_J\rangle$.
\end{lem}
\noindent The following subgroup is called the \emph{Levi subgroup of type~$J$} of~$G$ (with respect to~$(U_\alpha)_{\alpha\in \Phi}$):
\begin{equation*}
  L_J:=\langle U_\alpha\mid \alpha\in \Phi_J\rangle.  
\end{equation*}

\noindent One of the main reasons why we are interested in groups with root group data is summarised in the following theorem.

\begin{thm}[\protect{\cite[\S 4.1.2, Corollary~5.4, Proposition~5.5]{CapRem:root}}]\label{thm:twinBNassoc}
   The subgroups
    \begin{equation}\label{eq:TBN}
 T:=\langle\mu(u)^{-1}\mu(v)\mid u,v\in U_{s}\setminus\{1\},\,s\in S\rangle, \quad B_\pm:=TU_\pm\quad\text{and}\quad N:=T\langle \mu(u)\mid U_s\setminus\{1\},\,s\in S\rangle
   \end{equation}
  determine a twin BN-pair~$(B_+,B_-,N)$ of~$G$ of type~$(W,S)$ and thickness vector~$(q_s)_{s\in S}$, that is~$q_s=|B_+sB_+/B_+|=|B_-sB_-/B_-|$ for every~$s\in S$.
Moreover,
\begin{enumerate}
    \item \label{thm:twinBN1} $B_\pm$ splits as an inner semi-direct product~$B_\pm=U_\pm\rtimes T$;
    \item \label{thm:twinBN2} $U_+\cap B_-=\{1\}$;
    \item \label{thm:twinBN3} $B_+\cap B_-=B_+\cap N=B_-\cap N=T$. In particular, $T\unlhd N$ and there is a group isomorphism $\nu\colon N/T\to W$ such that~$\nu(uT)=\mu(u)=s$ for all~$u\in U_s\setminus\{1\}$ and~$s\in S$;
    \item \label{thm:twinBN4} $T=\bigcap_{\alpha\in \Phi}N_G(U_\alpha)$.
\end{enumerate}
\end{thm}
\begin{rem}\label{rem:Tnorm}
    By \Cref{thm:twinBNassoc}\eqref{thm:twinBN4},~$T\langle U_\alpha\mid \alpha\in \Psi\rangle$ is a subgroup of~$G$ for every subset~$\Psi\subseteq \Phi$.
\end{rem}
\begin{rem}\label{rem:symmBN}
    After \Cref{thm:twinBNassoc}, we can refine \Cref{rem:defnRGDsym} as follows. Given a root group datum~$(U_\alpha)_{\alpha\in \Phi}$ of type~$(W,S)$ for~$G$ with associated twin BN-pair~$(B_+,B_-,N)$, the datum~$(\wtil{U}_\alpha)_{\alpha\in \Phi}$ given by~$\wtil{U}_\alpha:=U_{-\alpha}$ for every~$\alpha\in \Phi$ is a root group datum for~$G$ of type~$(W,S)$, and its associated twin BN-pair~$(\wtil{B}_+,\wtil{B}_-,\wtil{N})$ satisfies
    \begin{equation*}
      \wtil{B}_\pm=B_\mp\qquad \text{and}\qquad \wtil{N}=N. 
    \end{equation*}
    Indeed, by \Cref{thm:twinBNassoc}\eqref{thm:twinBN4} we have~$\wtil{B}_+\cap \wtil{B}_-=B_+\cap B_-$ and hence~$\wtil{N}=W(\wtil{B}_+\cap \wtil{B}_-)=W(B_+\cap B_-)=N$. Moreover, using~\eqref{eq:TBN} we also deduce that~$\wtil{B}_\pm=B_\mp$.
\end{rem}

For every~$J\subseteq S$ denote by
\begin{equation*}
    P_J^+:=B_+W_JB_+\quad\text{and}\quad P_J^-:=B_-W_JB_-
\end{equation*}
the \emph{positive} and \emph{negative standard parabolic subgroups} of type~$J$, respectively. For simplicity, for every~$s\in S$ we typically write~$P_s^\pm$ instead of~$P_{\{s\}}^\pm$.

Since~$B_\pm=U_\pm T$ and~$W\cong N/T$, one easily deduces the following refinement of~\eqref{eq:birk}:
\begin{cor}[Refined Birkhoff decomposition]\label{cor:refBirk}
   Given~$(G,(U_\alpha)_{\alpha\in \Phi})$ as before, we have
    \begin{equation*}
        G=\bigsqcup_{w\in W}U_+wB_-=\bigsqcup_{w\in W}U_-wB_+.
    \end{equation*}
\end{cor}

\begin{cor}[\protect{\cite[\S 8.6.2(2)]{ab:build}}]\label{cor:wUconj}
   For all~$w\in W$ and~$\alpha\in \Phi$ we have~$wU_\alpha w^{-1}=U_{w\alpha}$. In particular, $sU_ss^{-1}=U_s^-$.
\end{cor}

For every~$w\in W$, following~\cite[\S 1.1.5 and p.~44]{CapRem:root}, consider the subgroup
\begin{equation}\label{eq:Uw}
    U_w:=\langle U_\alpha\mid \alpha\in \Phi^+,\,w^{-1}\alpha\in \Phi^-\rangle.
\end{equation}
If~$w=s\in S$, the corresponding group in~\eqref{eq:Uw} coincides with~$U_{\alpha_s}$~\cite[Lemma~1.5(ii)]{CapRem:root}. Thus the notation~$U_w$ is compatible with \Cref{notat:RGD}. Moreover, notice that~$U_1=\{1\}$.

\begin{lem}[\protect{\cite[Proposition~7.104 and~p.~472]{ab:build}, \cite[Corollary~5.4(i)]{CapRem:root}}]\label{lem:BwBUw}
For every~$w\in W$ we have~$B_+wB_+=U_wwB_+$ and~$U_+\cap wB_-w^{-1}=U_w$.
\end{lem}
\begin{rem}\label{rem:Ps-U}
   By \Cref{lem:BwBUw} and \Cref{cor:wUconj}, for every~$s\in S$ we have
   $$P_s^-=sP_s^-=sB_-\sqcup sU_s^-sB_-=sB_-\sqcup U_sB_-.$$ Since~$U_+\cap B_-=\{1\}$ (see \Cref{thm:twinBNassoc}) and~$U_+\cap sB_-=\emptyset$ (see \Cref{cor:refBirk}), we have~$U_+\cap P_s^-=U_s$. The latter formula is an instance of the more general \Cref{lem:UJUPinters}.
\end{rem}

From \Cref{lem:BwBUw} we deduce the following decomposition of~$P_J^+$. 
\begin{cor}\label{cor:weakLevi}
    For every~$J\subseteq S$ we have~$P_J^+=L_JB_+:=\{xb\mid x\in L_J,\,b\in B_+\}$.
\end{cor}
\begin{proof}
  We firstly check that~$L_JB_+\subseteq P_J^+$. By \Cref{cor:wUconj}, for every~$w\in W_J$ and~$s\in J$ we have~$U_{w\alpha_s}=wU_sw^{-1}\subseteq wB_+w^{-1}\subseteq P_J^+$. Hence~$L_J\subseteq P_J^+$ and~$L_JB_+\subseteq P_J^+$.

   It remains to show that~$B_+wB_+\subseteq L_JB_+$ for every~$w\in W_J$. By \Cref{thm:twinBNassoc} applied to~$(L_J,(U_\alpha)_{\alpha\in \Phi_J})$ we have~$w\in L_J$. By \Cref{lem:BwBUw}, we also have~$B_+wB_+=U_wwB_+$. We now claim that~$U_w\subseteq L_J$. In detail, $U_w=\langle U_\alpha\mid \alpha\in \Phi_w\rangle$ where~$\Phi_w:=\{\alpha\in \Phi^+\mid w^{-1}\alpha\in \Phi^-\}$. By~\cite[Lemma~1.5(iii)]{CapRem:root}, given a reduced expression~$w=s_1\dots s_n$ in~$(W,S)$, we have $\Phi_w=\{\alpha_{s_n},s_n\alpha_{s_{n-1}},\dots, s_n\cdots s_2\alpha_{s_1}\}$. Since~$w\in W_J$, we have~$s_1,\dots, s_n\in J$ (see~\cite[Theorem~5.5(b)]{hum:ref}) and hence~$\Phi_w\subseteq \Phi_J$. We have now proved that~$U_w\subseteq L_J$, and from this we conclude that~$B_+wB_+=U_wwB_+\subseteq L_JB_+$.
\end{proof}

\begin{rem}
    The decomposition in \Cref{cor:weakLevi} can be seen as a "rough version" of the so-called \emph{Levi decomposition} of~$P_J^+$, namely the inner semidirect product decomposition~$P_J^+:=TL_J\rtimes \caU_{S\setminus J}$, where~$\caU_{S\setminus J}:=\LL U_\alpha\mid \alpha\in \Phi^+\setminus \Phi_J\GG$. The latter decomposition holds for all~$J\subseteq S$ only if the root group datum satisfies a refined version of~\hyperref[RGD3]{(RGD3)}; see~\cite[\S 6.2.2]{rem02} for details or~\cite[Theorem~B.39]{marq:km} for an English version. For an arbitrary root group datum, the Levi decomposition for~$P_J^+$ given before holds only if~$J\subseteq S$ is assumed to be spherical~\cite[\S 6.2]{CapRem:root}.
\end{rem}

\smallskip

Let~$(B_+,B_-,N)$ be the twin BN-pair as in \Cref{thm:twinBNassoc}. Following \Cref{sus:build}, let
\begin{equation}\label{eq:buildRGD}
    (\caC_+,\caC_-,\delta)
\end{equation}
be the associated twin building and denote the simplicial complex corresponding to~$\caC_\pm$ by
\begin{equation}\label{eq:XpmRGD}
    X_\pm:=\caC\caC\Big(G, \big(P_{S\setminus\{s\}}^\pm\big)_{s\in S}\Big).
\end{equation}
Let also
\begin{equation}\label{eq:c0Delta0}
    c_0^\pm:=\big\{1P_{S\setminus \{s\}}^\pm\big\}_{s\in S}
\end{equation}
be the \emph{fundamental chamber} in~$X_\pm$.

\smallskip

A well-known chamber subsystem of~$\caC_-$ is the \emph{negative opposite system~$\,\caC_-^{\op}$}, i.e., the chamber system over~$S$ whose underlying set of chambers is 
\begin{equation}\label{eq:Copppm}
    \caC_-^{\op}:=\{c\in \caC_-\mid \delta(c,1B_+)=1\}=B_+\cdot 1B_-=U_+\cdot 1B_-.
\end{equation}
The last two passages in~\eqref{eq:Copppm} are due to~\eqref{eq:defndelta} and~$B_+=U_+T=U_+(B_+\cap B_-)$ (see \Cref{thm:twinBNassoc}). For every~$s\in S$, the $s$-adjacency relation in~$\caC_-^{\op}$ is given by restricting the respective $s$-adjacency relation in~$\caC_-$: i.e., for all~$u,u'\in U_+$,
\begin{equation}\label{eq:adjCop}
    uB_-\sim_s u'B_-\,\Longleftrightarrow\, u'u^{-1}\in U_s=U_+\cap P_s^-.
\end{equation}
see \Cref{rem:Ps-U} for the latter equality.

As recalled in \Cref{ex:GrpChambSys}, $\caC_-$ is $G$-equivariant isomorphic to the set~$\Ch(X_-)$ of chambers of~$X_-$ as chamber systems over~$S$. We denote by~$Y_-$ the subcomplex of~$X_-$ spanned by the image of~$\caC_-^{\op}$ in~$\Ch(X_-)$ via the previous isomorphism. In other words, 
\begin{equation}\label{eq:Ypm}
   Y_-=\Big\{\big\{uP_{S\setminus \{s\}}^-\big\}_{s\in J}\,\Big\vert\, u\in U_+,\,J\subseteq S\Big\}.
\end{equation}
We call~$Y_-$ the \emph{negative opposite complex}.

We will see later in \Cref{cor:Y-Uchar} that, under mild conditions on the root group datum, $Y_-$ can be described as a coset complex of~$U_+$ and the subgroups~$(U_{S\setminus\{s\}})_{s\in S}$.
\begin{rem}
    If~$(W,S)$ is spherical, $\caC_-^{\op}$ (exceptionally) admits a complete description in terms of~$\caC_-$ and its Weyl-distance~$\delta_-$. In detail, let~$w_0$ be the element of~$W$ of maximal Coxeter length. By~\cite[Corollary~5.6]{CapRem:root} we have~$B_-=w_0B_+w_0^{-1}$ and hence
    \begin{equation*}
        \caC_-^{\op}=\{c\in \caC_-\mid \delta_-(c,1B_-)=w_0\}.
    \end{equation*}
\end{rem}

The following theorem is a weaker form of the Curtis--Tits theorem by Abramenko--M\"uhlherr~\cite{AM97}, and will be used repeatedly in this paper.
\begin{thm}[\protect{\cite{AM97}, \cite[Theorem~8.2]{CapRem:root}}]\label{thm:OpConn}
  Let~$(W,S)$ be $2$-spherical and assume that, for all~$s,t\in S$ with~$s\neq t$, the Levi subgroup~$L_{\{s,t\}}$ satisfies
  \begin{equation*}\label{eq:co}
      \hspace{2cm} L_{\{s,t\}}/Z(L_{\{s,t\}})\not\in\{B_2(2)=\mathrm{Sp}_4(\F_2), G_2(\F_2), G_2(\F_3), {}^2F_4(\F_2)\}. \hspace{1.5cm} \mathrm{(co)}
  \end{equation*}
 Then~$U_+=\langle U_s\mid s\in S\rangle$. In other terms, $\caC_-^{\op}$ and hence~$Y_-$ are gallery-connected.
\end{thm}
\begin{rem}\label{rem:co}
    \begin{enumerate}
        \item \label{rem:co1} 
        Note that~\hyperref[eq:co]{\co} is automatically satisfied if~$\qmin\geq 4$.
        \item \label{rem:co2} Without the hypothesis of $2$-sphericity or without~\hyperref[eq:co]{\co}, the conclusion of \Cref{thm:OpConn} is false; see~\cite[Proposition~6.2]{MuhRon95} and~\cite[\S I.6, Corollary~8]{ab:twin}. See also~\cite[\S 5.11]{ab:build} (in particular, \cite[Remark~5.212]{ab:build}) for an overview.
    \end{enumerate}
\end{rem}

Combining \Cref{lem:RGDsubsyst} and \Cref{thm:OpConn} we deduce the following well-known result; see just after~\cite[Lemma~3.6]{GdP-VP:KMS} for the corresponding statement in the Kac--Moody case.
\begin{cor}\label{cor:UJUgroup}
Let~$(W,S)$ be $2$-spherical, assume~\hyperref[eq:co]{\co} and let~$J\subseteq S$. Then
\begin{equation*}
   U_J=\langle U_{\alpha}\mid \alpha\in \Phi_J^+\rangle\qquad\text{and}\qquad L_J=\langle U_s,U_s^{-}\mid s\in J\rangle.
\end{equation*}
In particular, $U_J$ is the "$U_+$-group" of the root group datum~$(U_\alpha)_{\alpha\in \Phi_J}$ of type~$(W_J,J)$ for~$L_J$.
\end{cor}
The following lemma is a generalisation of~\cite[Lemma~3.7]{GdP-VP:KMS}. 
\begin{lem}\label{lem:UJUPinters}
Let~$(W,S)$ be $2$-spherical and assume~\hyperref[eq:co]{\co}. Then, for all~$J\subseteq S$ and~$w\in W$ such that~$\ell(ws)>\ell(w)$ for every~$s\in J$ we have
\begin{equation*}
    U_+\cap wP_J^-w^{-1}=wU_Jw^{-1}U_w.
\end{equation*}
In particular, $U_+\cap P_J^-=U_J$ for every~$J\subseteq S$.
\end{lem}
\begin{proof}
    We first prove that~$wU_Jw^{-1}U_w\subseteq U_+\cap wP_Jw^{-1}$.
Since the Coxeter length~$\ell(\argu)$ on~$(W,S)$ is invariant under taking inverses, we have~$w\alpha_s\in \Phi^+$ for every~$s\in J$. By \Cref{cor:wUconj} and \Cref{cor:weakLevi} we deduce that
$$wU_Jw^{-1}=\langle wU_sw^{-1}\mid s\in J\rangle=\langle U_{w\alpha_s}\mid s\in J\rangle\subseteq U_+\cap wL_Jw^{-1}\subseteq U_+\cap wP_J^-w^{-1}.$$ 
Moreover, clearly~$U_w\subseteq U_+$, and by~\cite[Lemma~8.65]{ab:build} we also have~$w^{-1}U_ww=\langle U_\alpha\mid \alpha\in \Phi^-,\,w^{-1}\alpha\in \Phi^+\rangle$. Thus~$U_w\subseteq wU_-w^{-1}\subseteq wP_J^-w^{-1}$.

It remains to prove that~$U_+\cap wP_J^-w^{-1}\subseteq wU_Jw^{-1}U_w$. Let
\begin{equation*}
   \caD:=\{c\in \Res_J(wB_-)\mid \delta(c,1B_+)=w\},
\end{equation*}
where the residue is taken in~$\caC_-$ and~$\delta(\argu,\argu)$ is the twinning of~$(\caC_+,\caC_-)$.
We first show that
\begin{equation}\label{eq:Ddefn}
    \caD=wU_J\cdot 1B_-,
\end{equation}
Recall that~$P_J^-=L_JB_-$ by \Cref{cor:weakLevi}. Moreover, by \Cref{lem:RGDsubsyst} and \Cref{thm:twinBNassoc}, $L_J$ admits a twin BN-pair of type~$(W_J,J)$. We denote by~$B_-^J$ the associated "$B_-$-group". By \Cref{cor:UJUgroup}, the "$U_+$-group" of~$(L_J,(U_\alpha)_{\alpha\in \Phi_J})$ is~$U_J$. By \Cref{cor:refBirk} applied to~$L_J$, we have~$L_J=U_JW_JB_-^J$ and then~$P_J^-=L_JB_-=U_JW_JB_-$. Hence,
\begin{equation*}
    \Res_J(wB_-)=w\cdot\Res_J(1B_-)=wP_J^-\cdot 1B_-=wU_JW_J\cdot 1B_-=(wU_Jw^{-1})wW_J\cdot 1B_-.
\end{equation*}
From the first part of the proof we have~$wU_Jw^{-1}\subseteq U_+\subseteq B_+$. Hence, for all~$u\in wU_Jw^{-1}$ and~$x\in W_J$ we note that~$\delta(uwxB_-,1B_+)=\delta(wxB_-,1B_+)=wx$. Hence~$\caD=(wU_Jw^{-1})\cdot wB_-$ and~\eqref{eq:Ddefn} is proved.

Let~$u\in U_+\cap wP_J^-w^{-1}$. Notice that~$uwB_-\in \caD$: indeed,~$u\in wP_J^-w^{-1}$ implies that~$uwB_-\in w\cdot \Res_J(1B_-)=\Res_J(wB_-)$, and~$u\in U_+$ gives~$\delta(uwB_-,1B_+)=\delta(wB_-,1B_+)=w$.
By~\eqref{eq:Ddefn}, there is~$v\in wU_Jw^{-1}$ such that~$uwB_-=vwB_-$, i.e., $v^{-1}u\in wB_-w^{-1}$. As~$wU_Jw^{-1}\subseteq U_+$, we conclude that~$v^{-1}u\in U_+\cap wB_-w^{-1}=U_w$ (see \Cref{lem:BwBUw}). This proves that~$U_+\cap wP_J^-w^{-1}\subseteq wU_Jw^{-1}U_w$.
\end{proof}

\begin{cor}\label{cor:UJfinite}
Let~$G$ be a group with a locally finite root group datum~$(U_\alpha)_{\alpha\in \Phi}$ of  $2$-spherical type~$(W,S)$ and that satisfies~\hyperref[eq:co]{\co}. Then $U_J$ is finite for every spherical subset~$J\subseteq S$.
    If in particular~$(W,S)$ is spherical, then~$U_+$ is finite.
\end{cor}
\begin{proof}
   By~\Cref{rem:GBspher}, the coset space~$P_J^-/B_-$ is finite.
  Moreover, by \Cref{lem:UJUPinters} we have~$U_J=U_+\cap P_J^-$ and~$U_+\cap B_-=\{1\}$. Hence, the canonical projection~$G\twoheadrightarrow G/B_-$ restricts to an injective map~$U_J\hookrightarrow P_J^-/B_-$.
\end{proof}

As already observed in~\cite[\S 3.2, p.~11]{GdP-VP:KMS} for the Kac--Moody case, another consequence of \Cref{lem:UJUPinters} is the following characterisation of~$Y_-$ as a coset complex. 
\begin{cor}\label{cor:Y-Uchar}
     Let~$(W,S)$ be $2$-spherical and assume~\hyperref[eq:co]{\co}. Then the map
        \begin{equation*}
            uP_{S\setminus\{s\}}^-\longmapsto uU_{S\setminus\{s\}},\qquad \forall\,u\in U_+,\,s\in S
         \end{equation*}
    is well-defined and induces a $U_+$-equivariant simplicial isomorphism
    \begin{equation*}
       Y_-\longrightarrow \caC\caC\Big(U_+, \big(U_{S\setminus\{s\}}\big)_{s\in S}\Big). 
    \end{equation*}
\end{cor}

As a final note, we remark that twin buildings that come from a root group datum as discussed before form a broad family of twin buildings. For this remark, we assume familiarity with the definitions of twin buildings, twin roots and and associated root groups (see~\cite[Definitions~5.133, 5.195, and~8.16]{ab:build}). 
\begin{rem}[Root group data and Moufang twin buildings]\label{rem:Mouf}
   Recall that a thick\footnote{I.e., both~$\caC_+$ and~$\caC_-$ are thick buildings.} twin building~$(\caC_+,\caC_-)$ is \emph{Moufang} if, for every twin root~$\alpha$, the associated root group acts transitively on the set of all twin apartments that contain~$\alpha$~\cite[Definition~8.18]{ab:build}. 
    By~\cite[Theorem~8.81]{ab:build} and~\cite[Remarks~8.57(c)]{ab:build}, a twin building~$(\caC_+,\caC_-,\delta)$  arises from a group with a root group datum as described in~\eqref{eq:buildRGD} if and only if it is \emph{Moufang}. 

    Moreover, by~\cite[Theorem~8.27]{ab:build}, every thick irreducible $2$-spherical twin building of rank~$\geq 3$ that satisfies~\hyperref[eq:co]{\co} is Moufang.
\end{rem}
\subsection{Twin building lattices}\label{sus:TBL}
Let~$G$ be a group with a locally finite root group datum~$(U_\alpha)_{\alpha\in \Phi}$ of type~$(W,S)$.
The type-preserving actions of~$G$ on~$X_+$ and on~$X_-$ allow us to construct two metrics~$\mathrm{dist}_+$ and~$\mathrm{dist}_-$ on~$G$ (see~\cite[p.~51]{CapRem:root} for an explicit definition), and each of them makes~$G$ a topological group. We denote by~$\barG_\pm$ the topological closure of~$G$ with respect to~$\mathrm{dist}_\pm$. As recorded in~\cite[Proposition~6.1]{CapRem:root}, 
\begin{enumerate}
    \item \label{geom0a} $\barG_\pm$ is a tdlc group with respect to the topology induced by the extension of~$\mathrm{dist}_\pm$.
    \item \label{geom0b}  $\barG_\pm$ is discrete if and only if~$W$ is finite. 
    \item \label{geom1} Denote by $\overline{B}_\pm$ the closure of~$B_\pm$ in~$\barG_\pm$ (see \Cref{thm:twinBNassoc} for~$B_\pm$). Then the pair~$(\overline{B}_\pm,N)$ is a BN-pair of type~$(W,S)$ for~$\barG_\pm$, and there is a $G$-equivariant isomorphism of buildings
    \begin{equation*}
        \Delta(\barG_\pm,\overline{B}_\pm)\cong X_\pm=\Delta(G, B_\pm).
    \end{equation*}
    Moreover, the kernel of the $\barG_\pm$-action on~$X_\pm$ is discrete and is equal to~$Z(\barG_\pm)=Z(G)$.
    \item \label{geom2} Denoting by~$\varphi_\pm\colon G\hookrightarrow \barG_\pm$ the inclusion map, the image
    \begin{equation*}
        \Lambda:=(\varphi_+\times \varphi_-)(G)\leq \barG_+\times \barG_-
    \end{equation*}
    is a discrete subgroup of~$G_+\times G_-$. Moreover, $U_+$ is a discrete subgroup of~$G_-$.
\end{enumerate}

Following~\cite[Theorem~6.8]{CapRem:root}, if additionally the subgroup~$T$ is finite (see~\eqref{eq:TBN}) and~$\sum_{w\in W}(q_{\min})^{\ell(w)}$ converges (e.g.,~$\qmin>|S|$), then~$U_+$ is a lattice in~$G_-$ and~$\Lambda$ is an irreducible\footnote{I.e., the projection of~$\Lambda$ in each factor of~$\barG_+\times \barG_-$ is dense.} lattice in~$\barG_+\times \barG_-$.  
In this latter case, we say that~$\Lambda$ is the \emph{twin building lattice} of type~$(W,S)$ associated to the root group datum~$(U_\alpha)_{\alpha\in \Phi}$.

\begin{rem}\label{rem:TBL}
    By part~\eqref{geom1}, $\Lambda\cong G/Z(G)$. Moreover, if~$\Lambda\cong G/Z(G)$ is a twin building lattice, the centre~$Z(G)$ of~$G$ (which is contained in~$T$) is finite. 
\end{rem}

The following theorem records some notable properties of twin building lattices of non-Euclidean type proved in~\cite{caprem:KMsup}. 
\begin{thm}[\protect{\cite[Theorem~7.8]{CapRem:root}, \cite[Corollary~16, Theorem~19]{caprem:KMsup}}]\label{thm:TBLsimp}
    Let~$\Lambda$ be a twin building lattice of irreducible non-Euclidean type~$(W,S)$ and whose associated root groups~$U_s$, $s\in S$, are all (finite) nilpotent. Then~$\Lambda$ is an infinite, finitely generated, virtually simple group. If moreover~$(W,S)$ is $2$-spherical and~$\qmin\geq 3$, then~$\Lambda$ is infinite simple, hence it is not residually finite and non-linear.
\end{thm}
\begin{rem}
    The first part of \Cref{thm:TBLsimp} is~\cite[Theorem~7.8]{CapRem:root}. In the second part, the simplicity of~$\Lambda$  follows by combining~\cite[Corollary~16]{caprem:KMsup} and~\cite[Theorem~18]{caprem:KMsup}. Notice that, in the latter references, the group called by~$\Lambda^\dagger$ is our group~$G$, and our~$\Lambda$ is isomorphic to~$G/Z(G)$.
    
    Finally, an infinite simple group is clearly not residually finite. Given this, the non-linearity of~$\Lambda$ follows from a well-known theorem by Mal'cev: every finitely generated linear group is residually finite. These observations were already recorded within~\cite[\S 5]{caprem:KMsup}.
\end{rem}

\begin{rem}\label{rem:KMsimp}
  \Cref{thm:TBLsimp} applies to the twin building lattice of Kac--Moody type~$\caG_A(\F_q)^\dagger/Z(\caG_A(\F_q)^\dagger)$ as in \Cref{ex:KMgroups} provided~$q\geq 3$ and~$(W_A,S)$ is irreducible, non-Euclidean and $2$-spherical. By~\cite[\S 5.1]{marq:km} and~\cite[Proposition~4.22]{marq:km} and recalling \Cref{ex:genCart}, the previous properties on~$(W_A,S)$ are equivalent to asking, respectively, that~$A$ is of indefinite type, the support matrix of~$2I-A$ is the adjacency matrix of a connected graph, and~$A_{st}A_{ts}\leq 3$ for all~$s,t\in S$ with~$s\neq t$.
\end{rem}

\subsection{Connectivity and acyclicity properties of the complex~$Y_-$}\label{sus:homotY-}
For the subsection, we set once and for all a root group datum~$(U_\alpha)_{\alpha\in \Phi}$ of type~$(W,S)$ for a given group~$G$. Throughout, $\ell$ denotes the word length function of~$(W,S)$ by~$\ell(\argu)$, and~$X_-=\Delta(G,B_-)$ is the simplicial building associated to the negative BN-pair~$(B_-,N)$ of~$G$, i.e.,
\begin{equation}\label{eq:notaconn}
    X_-=\big\{\sigma(g,I)\,\big|\, g\in G,\,I\subseteq S\big\}, \qquad\sigma(g,I):=\big\{gP_{S\setminus \{s\}}^-\,\big|\, s\in I\}.
\end{equation}
By \Cref{cor:refBirk}, note that 
\begin{equation}\label{eq:X_-Uact}
    X_-=\{\sigma(uw,I)\mid u\in U_+,\,w\in W,\,I\subseteq S\}.
\end{equation}
Moreover,~$c_0^-=\sigma(1,S)=\big\{1P_{S\setminus\{s\}}^-\,\big\vert\, s\in S\big\}$ denotes the fundamental chamber of~$X_-$ (see~\eqref{eq:c0Delta0}). 

Since the set of chambers~$\Ch(X_-)$ of~$X_-$ is naturally identified with~$\caC_-$ via the map~$\sigma(g,S)\mapsto gB_-$ (see \Cref{ex:GrpChambSys}), in this section we sometimes assume that the Weyl distance~$\delta_-$ on~$\caC_-$ is defined on pairs of chambers in~$X_-$. We follow an analogous convention for the twinning~$\delta$ of the twin building~$(\caC_+,\caC_-)$ associated to the twin BN-pair~$(B_+,B_-,N)$.

\smallskip

By Solomon--Tits theorem, $X_-$ is $(|S|-2)$-connected and, when~$(W,S)$ is non-spherical, it is even contractible. In particular, $X_-$ is always $(|S|-2)$-acyclic over~$\C$. The goal of this section is to establish that~$Y_-$ is $(k-1)$-connected or, at least, $(k-1)$-acyclic over~$\Q$, whenever~$(W,S)$ is $(k+1)$-spherical (plus some other mild hypotheses on the root group datum); see \Cref{thm:Yacyc} and \Cref{thm:Y-homot}. 

A common strategy, also adopted here, to transfer connectivity and acyclicity properties from a CW~complex~$X$ to a given subcomplex~$Y$ is the following: Construct a (countable) \emph{filtration}~$(X_d)_{d\in \Z_{\geq 0}}$ of~$X$ (i.e., a sequence of subcomplexes~$X_d$ of~$X$ such that~$X_d\subseteq X_{d+1}$ for every~$d\geq 0$ and~$X=\bigcup_{d\geq 0} X_d$) such that~$X_0=Y$ and that satisfies the hypotheses of the following well-known fact:
\begin{lem}\label{prop:filtHom}
    Let~$X$ be a CW~complex with a filtration~$X=\bigcup_{d\in \Z_{\geq 0}}X_d$ by subcomplexes with~$X_0\neq\emptyset$. For every~$d\geq 1$, let~$(Y_{id})_{i\in I_d}$ be a family of contractible subcomplexes of~$X$ such that, for all~$d\geq1$,
    \begin{equation}\label{eq:filtr}
        X_d=X_{d-1}\cup \bigcup_{i\in I_d}Y_{id}\qquad \text{and}\qquad Y_{id}\cap Y_{jd}\subseteq X_{d-1},\quad\forall\,i,j\in I_d,\,i\neq j.
    \end{equation}
   Let~$m\geq 0$ be an integer such that, for all~$d\geq 1$ and~$i\in I_d$,~$X_{d-1}\cap Y_{id}$ is $m$-connected. Then if~$X$ is $m$-connected, so is~$X_d$ for every~$d\geq 0$.

   The same statement holds by replacing "contractible" with "acyclic over~$R$" and "$m$-connected" with "$m$-acyclic over~$R$", where~$R$ is any unitary ring.
\end{lem}
\begin{proof}
  Let~$\ell\geq 0$ and~$x_0\in X_0$. The inclusion~$X_d\hookrightarrow X$,~$d\geq 0$, induces an isomorphism of groupoids
   \begin{equation*}
       \varinjlim_{d\geq 0}\pi_\ell(X_d,x_0)\cong\pi_\ell(X,x_0),
   \end{equation*}
   where~$(\pi_\ell(X_d,x_0))_{d\geq 0}$ is the direct system with transition maps~$\iota_{d,e}^\ell\colon\pi_\ell(X_d,x_0)\to \pi_\ell(X_{d+e},x_0)$ induced by the inclusions~$X_d\hookrightarrow X_{d+e}$ for all~$d,e\geq 0$.
   By~\cite[Lemma~1.7]{sches}, for every~$d\geq 0$ the pair~$(X_{d+1},X_d)$ is $(m+1)$-connected. The inclusion~$X_d\hookrightarrow X_{d+1}$ induces the following long exact sequence of groupoids~\cite[\S 4.4, Item~4.4.9]{GeoBook}: 
   \begin{equation*}
       \xymatrix{\dots\ar[r] & \pi_{\ell+1}(X_{d+1},X_d,x_0) \ar[r] & \pi_\ell(X_d,x_0)\ar[r]^-{\iota_{d,1}^\ell} & \pi_\ell(X_{d+1},x_0) \ar[r] & \pi_\ell(X_{d+1},X_d,x_0)\ar[r] & \dots }
   \end{equation*}
   Hence, $\iota_{d,1}^\ell$ is an isomorphism for all~$d\geq 0$ and~$0\leq \ell\leq m$. Since~$\iota_{d,e}^\ell=\iota_{d+e-1,1}^\ell\circ\dots \circ \iota_{d,1}^\ell$ for every~$e\geq 1$,~$\iota_{d,e}^\ell$ is an isomorphism for all~$d\geq 0$, $e\geq 1$ and~$0\leq \ell\leq m$. Thus~$\pi_\ell(X,x_0)\cong \pi_\ell(X_d,x_0)$ for all~$d\geq 0$ and~$0\leq \ell\leq m$.

  For the statement in homology, one proceeds analogously using the long exact sequence in relative homology~\cite[Proposition~2.8.1]{GeoBook} instead of the long exact sequence in homotopy.
\end{proof}

Our strategy is to apply \Cref{prop:filtHom} to a filtration of~$X_-$ introduced by Abels~\cite[\S 2.3]{abels}, whose construction is now recalled.

\begin{defn}[Abels' filtration for the $U_+$-action on~$X_-$]\label{defn:Abfiltr}
    Denote by~$d_-$ the gallery-distance on~$X_-$, i.e., $d_-=\ell\circ\delta_-$ where~$\delta_-$ is the Weyl-distance on~$X_-$ and~$\ell$ is the word length on~$(W,S)$. For every simplex~$\sigma\in X_-$, let
    \begin{equation*}
        h(\sigma):=d_-(\sigma,Y_-)=\min\{d_-(c,d)\mid c\in \Ch(X_-),\,c\supseteq \sigma,\,d\in \Ch(Y_-)\}.
    \end{equation*}
    For every~$d\in \Z_{\geq 0}$, define the following subcomplex of~$X_-$:
    \begin{equation*}
        \Delta_d:=\{\sigma\in X_-\mid h(\sigma)\leq d\}.
    \end{equation*}
\end{defn}

\begin{rem}\label{rem:filtr}
\begin{enumerate}
    \item Clearly, $(\Delta_d)_{d\geq 0}$ is a filtration of~$X_-$ satisfying~$\Delta_0=Y_-$.
    \item \Cref{defn:Abfiltr} is the specialisation of the general filtration given in~\cite[\S 4.1]{abels}. Following the notation therein, in our case we have~$\Gamma=U_+$, $\Delta=X_-$ and the simplicial strict fundamental domain of the $U_+$-action on~$X_-$ is the fundamental apartment~$\Sigma_-=\{\sigma(w,I)\mid w\in W,\,I\subseteq S\}$ of~$X_-$; see~\eqref{eq:X_-Uact}. Since~$\Ch(Y_-)=U_+\cdot c_0^-$, we have~$h(\sigma)=d_-(\sigma, \Gamma\cdot c_0^-)$ for all~$\sigma\in X_-$. 
\end{enumerate}
\end{rem}
\begin{lem}\label{lem:hl}
   For every chamber~$c=uw\cdot c_0^-$ in~$X_-$ with~$u\in U_+,\,w\in W$, we have~$h(c)=\ell(w)$.
\end{lem}
\begin{proof}
   Recall from~\eqref{eq:Copppm} that~$\Ch(Y_-)=\{c\in \Ch(X_-)\mid \delta(c,c_0^+)=1\}$, where~$\delta$ denotes the twinning of~$(X_+,X_-)$. By~\cite[Lemma~5.143]{ab:build},~$h(c)=d_-(c,Y_-)=\ell(\delta(c,c_0^+))$. Since~$U_+\subseteq B_+$, we also have~$\delta(c,c_0^+)=\delta(wc_0^-,c_0^+)=w$, and the claim is proved.
\end{proof}

For every~$w\in W$ let
\begin{equation}
    D_R(w):=\{s\in S\mid \ell(ws)<\ell(w)\}.
\end{equation}
The set~$D_R(w)$ is usually called the \emph{right descend set} of~$w$. Since~$\ell(ws)\in\{\ell(w)\pm 1\}$, we have
\begin{equation*}
    I(w):=S\setminus D_R(w)=\{s\in S\mid \ell(ws)>\ell(w)\}.
\end{equation*}
It is easily checked that
\begin{equation}\label{eq:I(w)prop}
        I(w)=S \Longleftrightarrow w=1.
\end{equation}
Moreover,~$I(w)=\emptyset$ if and only if~$w$ is an element of maximal length in~$W$. In particular, the case~$I(w)=\emptyset$ can happen if and only if $W$~is finite~\cite[Proposition~2.3.1]{bb:cox}.

\begin{lem}[\protect{\cite[\S 2.1-2.2, Lemma~2.4]{abels}}]\label{lem:restr}
    The following hold for every~$d\geq 0$.
    \begin{enumerate}
        \item \label{lem:restr1} Let~$c$ be a chamber in~$\Delta_{d+1}\setminus \Delta_d$ and let~$u\in U_+$, $w\in W$ be such that~$c=uw\cdot c_0^-$. Then there is a unique smallest face~$\caR^{U_+}(c)$ of~$c$ that is not contained in~$\Delta_d$, namely
        \begin{equation*}
            \caR^{U_+}(c)=\big\{uwP_{S\setminus\{s\}}^-\,\big\vert\, s\in D_R(w)\big\}.
        \end{equation*}
       In particular, $h(\caR^{U_+}(c))=h(c)=\ell(w)$.
        \item \label{lem:restr2} For all~$c,c'\in \Ch(\Delta_{d+1})$, either~$c\cap c'\in \Ch(\Delta_d)$ or~$\caR^{U_+}(c)=\caR^{U_+}(c')$ and~$c'\in U_+\cdot c$.
    \end{enumerate}
\end{lem}

For every~$d\geq 0$, let
\begin{equation*}
\begin{split}
    R_{d+1}& :=\{\caR^{U_+}(c)\mid c\in \Ch(\Delta_{d+1})\setminus \Ch(\Delta_d)\}\\
    & =\{u\cdot \sigma(w,D_R(w))\mid u\in U_+,\,w\in W,\,\ell(w)=d+1\}.
\end{split}
\end{equation*}
The latter equality is due to \Cref{lem:restr}\eqref{lem:restr1}.
For every~$\rho\in R_{d+1}$, let also
\begin{equation}\label{eq:STrho}
S(\rho):=\{\sigma\in \Delta_{d+1}\mid \sigma\cup\rho\in \Delta_{d+1}\}\quad\text{and}\quad T(\rho):=\Delta_d\cap S(\rho).
\end{equation}
Notice that~$S(\rho)$ is precisely the star of~$\rho$ in~$\Delta_{d+1}$.

Before stating the following lemma, we recall the notion of \emph{projection}~$\proj_\sigma(c)$ of a chamber~$c\in \Ch(X_-)$ on any simplex~$\sigma\in X_-$. This definition is given for arbitrary buildings; see~\cite[\S 4.9]{ab:build} for background. In detail, given~$c$ and~$\sigma$ as before, $\proj_\sigma(c)$ is the unique chamber~$c'$ of~$X_-$ containing~$\sigma$ such that~$d_-(c,c')=d_-(c,\sigma)$, i.e., the unique chamber~$c'\supseteq \sigma$ at minimal distance from~$c$. Recall that~$\proj_\sigma(c)$ is contained in every apartment of~$X_-$ that contains both~$c$ and~$\sigma$~\cite[Proposition~4.95]{ab:build}. 

\begin{lem}[\protect{\cite[Lemma~4.2]{abels}, \cite[Lemma~7]{ab94:FPchev}}]\label{lem:Deltadprop}
    The following hold for every~$d\geq 0$:
    \begin{enumerate}
        \item \label{lem:Deltaprop1} $\Delta_{d+1}=\Delta_d\cup \bigcup_{\rho\in R_{d+1}}S(\rho)$ and~$S(\rho)$ is contractible for every~$\rho\in R_{d+1}$.
        \item \label{lem:Deltaprop2} For all pairwise distinct~$\rho,\rho'\in R_{d+1}$ we have~$S(\rho)\cap S(\rho')\subseteq \Delta_d$.
        \item \label{lem:Deltaprop3} Let~$\rho=\sigma(w,D_R(w))\in R_{d+1}$ for some~$w\in W$. Then~$T(u\cdot \rho)=u\cdot T(\rho)$ for all~$u\in U_+$.
        
\noindent Moreover, $T(\rho)$ splits as a join of simplicial complexes\footnote{Given two simplicial complexes~$\Sigma_1$ and~$\Sigma_2$, the join~$\Sigma_1\ast\Sigma_2$ is the simplicial complex whose simplices are all disjoint unions~$\sigma_1\sqcup \sigma_2$, where~$\sigma_i$ is a (possibly empty) simplex of~$\Sigma_i$ and~$i=1,2$.}
\begin{equation*}
    T(\rho)\cong\partial \rho\ast T'(\rho),
\end{equation*}
where
\begin{equation*}
\begin{split}
   \partial\rho & =\{\sigma\mid \sigma\subsetneq \rho\}\simeq \mathbb{S}^{\dim(\rho)-1},\\ 
   T'(\rho) &=(U_+)_\rho\cdot \{\tau\subseteq \proj_\rho(c_0^-)\mid \tau\cap\rho=\emptyset\}=(U_+)_\rho\cdot \{\tau\mid \tau\subseteq \proj_\rho(c_0^-)\setminus \rho\}.
\end{split} 
\end{equation*}
\end{enumerate}
\end{lem}

The next lemma is a possible generalisation of~\cite[\S I.5, Corollary~6]{ab:twin} to an arbitrary~$\rho\in \bigcup_{d\geq 0}R_{d+1}$. In other words, we do not require that the cotype (i.e., the complement of the type) of~$\rho$ is a spherical subset of~$S$.
\begin{prop}\label{prop:descLinks}
  Let~$(W,S)$ be $2$-spherical and assume~\hyperref[eq:co]{\co}. Then, for every~$\rho=\sigma(w,D_R(w))\in \bigcup_{d\geq 0}R_{d+1}$, there is a simplicial isomorphism
    \begin{equation*}
        T'(\rho)\cong \caC\caC\Big(U_{I(w)}, \big(U_{I(w)\setminus\{s\}}\big)_{s\in I(w)}\Big).
    \end{equation*}
    That is, $T'(\rho)$ is simplicially isomorphic to the "$Y_-$-complex" induced by the root group datum~$(U_\alpha)_{\alpha\in \Phi_{I(w)}}$ for~$L_{I(w)}$ (see \Cref{lem:RGDsubsyst}).
\end{prop}
\begin{proof}
Recall that the special coset~$wW_{I(w)}$ has a unique element~$w_1$ of minimal Coxeter length~$\ell$~\cite[Proposition~2.20]{ab:build}. Since~$w_1^{-1}w\in W_{I(w)}\subseteq P_{I(w)}^-$ and~$P_{I(w)}^-$ is the pointwise stabiliser of the simplex~$\sigma(1,D_R(w))$, we deduce that~$\rho=w\cdot\sigma(1,D_R(w))=w_1\cdot\sigma(1,D_R(w))=\sigma(w_1,D_R(w))$. 
Recall that~$c_0^-=\{1P_{S\setminus\{s\}}\mid s\in S\}$.
Given~$c:=\proj_\rho(c_0^-)$, we now claim that
\begin{equation}\label{eq:proj}
    c=w_1\cdot c_0^-.
\end{equation}
Both~$\rho$ and~$c_0^-$ belong to the fundamental apartment~$\Sigma_-=\{\sigma(w',I)\mid w'\in W,\,I\subseteq S\}$. By~\cite[Theorem~4.97]{ab:build}, the chamber~$c$ belongs to~$\Sigma_-$ and hence~$c=w'\cdot c_0^-$ for some~$w'\in W$. 
Since both~$c=w'\cdot c_0^-$ and~$w_1\cdot c_0^-$ are chambers containing~$\rho=\sigma(w_1,D_R(w))$, we have~$w_1^{-1}w'\in P_{I(w)}^-\cap W=W_{I(w)}$. The latter equality is standard and follows from the Bruhat decompositions of~$G$ and~$P_{I(w)}^-$. Since~$c=\proj_\rho(c_0^-)$, we notice that
\begin{equation*}
    \ell(w')=d_-(c,c_0^-)\leq d_-(w_1c_0^-,c_0^-)=\ell(w_1).
\end{equation*}
Since~$w'\in w_1W_{I(w)}$ and~$w_1$ is the unique element of smallest Coxeter length in that coset, this yields~$w'=w_1$ and hence~\eqref{eq:proj}.

Notice also that the simplex~$c\setminus \rho$ is the only face of~$c$ of type~$I(w)$, i.e.,~$c\setminus \rho=\sigma(w_1,I(w))$.
By \Cref{lem:Deltadprop}\eqref{lem:Deltaprop3},
    \begin{equation*}
 T'(\rho)=\{\sigma(uw_1,J)\mid u\in (U_+)_\rho,\,J\subseteq I(w)\}=(U_+)_{\rho}w_1\cdot\Delta_{I(w)},
    \end{equation*}
 where~$\Delta_{I(w)}=\{\sigma(1,J)\mid J\subseteq I(w)\}$.
Clearly,~$\ell(w_1s)>\ell(w_1)$ for every~$s\in I(w)$. Hence, by \Cref{lem:UJUPinters} we have
 \begin{equation*}
(U_+)_{\rho}=U_+\cap w_1P_{I(w)}^-w_1^{-1}=w_1U_{I(w)}w_1^{-1}U_{w_1}.
 \end{equation*}
Moreover, by \Cref{lem:BwBUw}, $U_{w_1}=U_+\cap w_1B_-w_1^{-1}$ and hence~$U_{w_1}$ is the pointwise stabiliser of the chamber~$w_1c_0^-$ in~$U_+$. This implies that~$U_{w_1}$ acts trivially on~$w_1\cdot\Delta_{I(w)}$ and
\begin{equation*}
    T'(\rho)=w_1U_{I(w)}w_1^{-1}U_{w_1}w_1\cdot \Delta_{I(w)}=w_1U_{I(w)}w_1^{-1}(w_1\cdot\Delta_{I(w)})=w_1U_{I(w)}\cdot \Delta_{I(w)}.
\end{equation*}
 Hence, the left action by~$w_1^{-1}$ induces a simplicial isomorphism
 \begin{equation*}
     T'(\rho)=w_1U_{I(w)}\cdot \Delta_{I(w)}\stackrel{\cong}{\longrightarrow}U_{I(w)}\cdot \Delta_{I(w)}.
 \end{equation*}

 We now focus on~$U_{I(w)}\cdot\Delta_{I(w)}$.
 For every~$J\subseteq I(w)$, from~\eqref{eq:intersParab} and \Cref{lem:UJUPinters} we obtain
 \begin{equation*}
 \begin{split}
     \stab_{U_{I(w)}}\hspace{-0.1cm}\Big(\big\{1P^-_{S\setminus \{s\}}\mid s\in J\big\}\Big) &=U_{I(w)}\cap P_{S\setminus J}^-=U_+\cap P_{I(w)}^-\cap P_{S\setminus J}^-=U_+\cap P_{I(w)\cap S\setminus J}^-\\
     &=U_+\cap P_{I(w)\setminus J}^-=\bigcap_{s\in J}\big(U_+\cap P_{I(w)\setminus\{s\}}^-\big)=\stab_{U_+}\hspace{-0.1cm}\Big(\{1U_{I(w)\setminus\{s\}}\mid s\in J\}\Big). 
 \end{split}
 \end{equation*}
 Therefore, the assignment
 \begin{equation*}
   \big\{uP_{S\setminus \{s\}}^-\,\big\vert\,s\in J\}\longmapsto \big\{uU_{I(w)\setminus \{s\}}\,\big\vert\, s\in J\big\}\qquad \forall\,u\in U_{I(w)},\,J\subseteq S
 \end{equation*}
 is a well-defined simplicial isomorphism
 \begin{equation*}
     U_{I(w)}\cdot \Delta_{I(w)}\longrightarrow \caC\caC\Big(U_{I(w)}, \big(U_{I(w)\setminus\{s\}}\big)_{s\in I(w)}\Big).\qedhere
 \end{equation*}
\end{proof}

\begin{prop}\label{prop:YconnSph}
   Let~$G$ be a group with a locally finite root group datum~$(U_\alpha)_{\alpha\in \Phi}$ of spherical type~$(W,S)$. Let~$|S|=n+1\geq 1$ and assume~\hyperref[eq:co]{\co}. Suppose also that, for all~$s,t\in S$ with~$s\neq t$, 
    \begin{equation}\label{eq:thickSpec1}
     \sqrt{\frac{q_sq_t}{\max\{q_s,q_t\}}}> 2n\cdot\cos\Big(\frac{\pi}{m_{st}}\Big),
    \end{equation}
    where~$M=[m_{st}]_{s,t\in S}$ is the Coxeter matrix of~$(W,S)$. Then~$Y_-$ is $(n-1)$-acyclic over~$\Q$.
\end{prop} 
Before proving \Cref{prop:YconnSph} we observe the following general fact, which will be useful anyway for some arguments in \Cref{sus:TkRGD}. 
Recall that, given a finite graph~$X$, $\lambda_2(X)$ denotes the second largest eigenvalue of the normalised weighted adjacency matrix of~$X$ (see~\eqref{eq:GGamma}), and that~$TU_{\{s,t\}}$, $TU_s$ and~$TU_t$ are subgroups of~$G$ for all~$s,t\in S$ (see \Cref{rem:Tnorm}).

\begin{lem}\label{lem:spec1link}
   Let~$G$ be a group with a locally finite root group datum~$(U_\alpha)_{\alpha\in \Phi}$ of type~$(W,S)$. Assume~$|S|\geq 2$ and let~$M=[m_{st}]_{s,t\in S}$ be the Coxeter matrix of~$(W,S)$. Then, for all~$s,t\in S$ with~$s\neq t$,   given~$\Lambda_{st}:=\caC\caC(U_{\{s,t\}},(U_s,U_t))$ and~$\Gamma_{st}:=\caC\caC(TU_{\{s,t\}},(TU_s,TU_t))$ we have
   \begin{equation*}
       \lambda_2(\Lambda_{st})\leq \lambda_2(\Gamma_{st})\leq 2\cdot\cos\Big(\frac{\pi}{m_{st}}\Big)\cdot
       \sqrt{\frac{\max\{q_s,q_t\}}{q_sq_t}}.
   \end{equation*}
\end{lem}

\begin{proof}[Proof of \Cref{lem:spec1link}]
By design, $\Lambda_{st}$ and~$\Gamma_{st}$ are bipartite graphs with respect to the defining vertex bipartitions:
\begin{equation*}
    \Lambda_{st}(0)=U_{\{s,t\}}/U_s\sqcup U_{\{s,t\}}/U_t\qquad\text{and}\qquad \Gamma_{st}(0)=TU_{\{s,t\}}/TU_s\sqcup TU_{\{s,t\}}/TU_t.
\end{equation*}
Notice that they are both $(q_s,q_t)$-regular. Indeed, for every~$x\in\{s,t\}$ we have
$$|U_x:U_s\cap U_t|=|U_x|=q_x$$ because~$U_s\cap U_t=U_+\cap P_s^{-}\cap P_t^-=U_+\cap B_-=\{1\}$;
 see~\Cref{rem:Ps-U}, \eqref{eq:intersParab} and \Cref{thm:twinBNassoc}\eqref{thm:twinBN2}. 
Moreover,~$T\cap U_+\subseteq B_-\cap U_+=\{1\}$ by \Cref{thm:twinBNassoc}\eqref{thm:twinBN2}. Hence,~$T\cap U_s=T\cap U_t=\{1\}$,~$TU_s\cap TU_t=T$ and, for every~$x\in\{s,t\}$,
$$|TU_x:TU_s\cap TU_t|=|TU_x:T|=|U_x|=q_x.$$ 
From~$T\cap U_+=\{1\}$, for every~$x\in\{s,t\}$ we also have
\begin{equation*}
    U_{\{s,t\}}\cap TU_x=U_x.
\end{equation*}
Hence,
\begin{equation*}
    \Lambda_{st}\longrightarrow \Gamma_{st}\qquad uU_x\longmapsto uTU_x,\,\forall\,u\in U_{\{s,t\}},\,x\in \{s,t\}
\end{equation*}
induces an injective simplicial map, and $\Lambda_{st}$ embeds in~$\Gamma_{st}$ as a (necessarily induced) subgraph. By  \Cref{prop:specSubgraph},
\begin{equation*}
    \lambda_2(\Lambda_{st})\leq \lambda_2(\Gamma_{st}).
\end{equation*}

In turn, for every~$x\in\{s,t\}$ we notice that
\begin{equation*}
    P_x^-\cap TU_{\{s,t\}}=TU_x.
\end{equation*}
Indeed, by \Cref{cor:weakLevi} we have~$P_x^-=B_-L_x=TU_-\langle U_x,U_x^-\rangle$. Hence~$TU_x\subseteq P_x^-\cap TU_{\{s,t\}}$. Since~$T\subseteq P_x^-$ and~$P_x^-\cap U_+=U_x$ (see \Cref{lem:UJUPinters}), we also have~$P_x^-\cap TU_{\{s,t\}}\subseteq P_x^-\cap TU_+=TU_x$. 
Given
  \begin{equation*}
      \Xi_{st}:=\Delta(P_{\{s,t\}}, B_+)=\caC\caC\big(P_{\{s,t\}}^-, (P_s^-, P_t^-)\big),
  \end{equation*}
we deduce that 
\begin{equation*}
     \Gamma_{st}\longrightarrow \Xi_{st}\qquad huTU_x\longmapsto huP_x^{-},\,\forall\,h\in T,\,u\in U_{\{s,t\}},\,x\in \{s,t\}
\end{equation*}
is an injective simplicial map.

Finally,~$\Xi_{st}$ is a simplicial building of spherical type~$(W_{\{s,t\}}, \{s,t\})$ and thickness vector~{$(q_s, q_t)$} (see \Cref{sus:build}). In particular,~$\Xi_{st}$ is a finite $(q_s+1,q_t+1)$-regular $(P_{\{s,t\}}^-/P_s^-,P_{\{s,t\}}^-/P_t^-)$-partite graph.
By \Cref{prop:specSubgraph},
  \begin{equation*}
      \lambda_2(\Lambda_{st})\leq \sqrt{\frac{(q_s+1)(q_t+1)}{q_sq_t}}\cdot\lambda_2(\Gamma_{st}).
  \end{equation*}
and now~\eqref{eq:AXineq} applies.
\end{proof}

\begin{proof}[Proof of \Cref{prop:YconnSph}]
    By \Cref{cor:Y-Uchar}, $Y_-$ is isomorphic to~$Z:=\caC\caC(U_+,(U_{S\setminus\{s\}})_{s\in S})$. Hence, by \Cref{lem:CCprop}\eqref{lem:CCprop1}-\eqref{lem:CCprop2}, \Cref{thm:OpConn} and \Cref{cor:UJfinite}, $Z$ is finite pure $n$-dimensional simplicial complex whose positive dimensional links are connected and whose $1$-dimensional links are isomorphic to the coset graph~$\caC\caC(U_{\{s,t\}}, (U_s,U_t))$ for some~$s,t\in S$ with~$s\neq t$. By \Cref{lem:spec1link} and~\eqref{eq:thickSpec1}, 
    \begin{equation*}
        \sup_{\sigma\in Z(n-2)}\lambda_2(Z_\sigma)=\max_{s,t\in S,\,s\neq t}\lambda_2\big(\caC\caC(U_{\{s,t\}}, (U_s,U_t))\big)<\frac{1}{n}.
    \end{equation*}
   By \Cref{thm:OppVan}\eqref{thm:OppVan2}, $Y_-$ is $(n-1)$-acyclic over~$\Q$.
\end{proof}

\begin{rem}
\cite[\S II.2]{ab:twin} collects some examples of locally finite root group data of spherical type~$(W,S)$ in which~$Y_-$ is not $(n-1)$-acyclic over~$\Q$, where~$|S|=n+1$. See in particular~\cite[\S II.2, Proposition~10(ii), Examples~8 and~9]{ab:twin}. In all these cases, \emph{at least one} of the following conditions holds:
\begin{itemize}
    \item $|S|\geq 3$ and there are~$s,t\in S$ with~$m_{st}=3$ and~$2\leq q_s^\pm=q_t^\pm\leq 3$; or,
    \item $|S|\geq 4$ and there are~$s,t\in S$ with~$m_{st}=3$ and~$2\leq q_s^\pm=q_t^\pm\leq 4$.
\end{itemize}
In each situation,~\eqref{eq:thickSpec1} is not satisfied. This means that, without assuming~\eqref{eq:thickSpec1}, the conclusion of \Cref{prop:YconnSph} might not hold.
\end{rem}

For the following theorem, it might be useful to recall the definitions of~$\sigma(W,S)$ and~$\kappa(W,S)$ for a given Coxeter system~$(W,S)$ given in~\eqref{eq:sk(W,S)}.
\begin{thm}\label{thm:Yacyc}
    Let~$G$ be a group with a locally finite root group datum~$(U_\alpha)_{\alpha\in \Phi}$ of $2$-spherical type~$(W,S)$. Suppose also that~\hyperref[eq:co]{\co} holds and that, for all~$s,t\in S$ with~$s\neq t$, 
     \begin{equation}\label{eq:thickSpec2}
    \sqrt{\frac{q_sq_t}{\max\{q_s,q_t\}}} > 2(\sigma(W,S)-1)\cos\Big(\frac{\pi}{m_{st}}\Big).
    \end{equation}
    where~$M=[m_{st}]_{s,t\in S}$ is the Coxeter matrix of~$(W,S)$. 
    Then~$Y_-$ is $(\kappa(W,S)-2)$-acyclic over~$\Q$.
\end{thm}
\begin{proof}
We proceed by induction on~$f(W,S):=|S|-\kappa(W,S)\geq 0$.
Notice that~$\kappa(W,S)=|S|$ if and only if~$(W,S)$ is spherical. The base case~$f(W,S)=0$ is due to \Cref{prop:YconnSph}. 

Let now~$m\geq 1$ and assume the claim true whenever the type~$(W',S')$ satisfies~$f(W',S')<m$. Suppose that~$(U_\alpha)_{\alpha\in \Phi}$ is as in the hypotheses and its type~$(W,S)$ satisfies~$f(W,S)=m$. For every~$J\subsetneq S$, we claim that
\begin{equation}\label{eq:hmonot}
    f(W_J,J)<f(W,S).
\end{equation}
Since~$\kappa(W,S)=m\geq 1$, \eqref{eq:hmonot} is clear if~$J$ is spherical. Moreover, if~$(W,S)$ is $k$-spherical then~$(W_J,J)$ is $k$-spherical. Hence,
\begin{equation}\label{eq:kmonot}
    \kappa(W,S) \leq \max\{0\leq k\leq |S|\colon (W_J,J)\text{ is $k$-spherical}\}.
\end{equation}
If~$J$ is non-spherical, the term on the right-hand side of~\eqref{eq:kmonot} at most~$|J|-1$ and hence necessarily equals~$\kappa(W_J,J)$. This yields~$\kappa(W,S)\leq \kappa(W_J,J)$. Having assumed~$J\subsetneq S$, we now deduce that~$|J|-\kappa(W_J,J)<|S|-\kappa(W,S)$. This completes the proof of~\eqref{eq:hmonot}.

Consider now Abels' filtration~$X_-=\bigcup_{d\geq 0}\Delta_d$ for~$X_-$ as in \Cref{defn:Abfiltr}, and recall that~$\Delta_0=Y_-$. Let~$d\geq 0$. By \Cref{lem:Deltadprop}, $\Delta_{d+1}=\Delta_d\cup\bigcup_{\rho\in R_{d+1}}S(\rho)$ and, for every~$\rho\in R_{d+1}$, 
\begin{equation*}
   S(\rho)\cap \Delta_d=T(\rho)=\partial\rho\ast T'(\rho).
\end{equation*}
By \Cref{lem:restr}\eqref{lem:restr1} and~\eqref{eq:notaconn}, $\rho=u\cdot\sigma(w,D_R(w))$ for some~$u\in U_+$ and~$w\in W$. Since~$\rho\not\in \Delta_0$, we have~$w\neq 1$. Recall also that~$(W,S)$ is non-spherical.
Then, by~\eqref{eq:I(w)prop} and the lines just afterwards,~$I(w):=S\setminus D_R(w)$ is a proper non-empty subset of~$S$. Being homotopy equivalent to~$\mathbb{S}^{\dim(\rho)-1}$, $\partial \rho$ is $(|D_R(w)|-3)$-connected and hence $(|D_R(w)|-3)$-acyclic over~$\Q$. Moreover, by \Cref{prop:descLinks}, $T'(\rho)$ is the "$Y_-$-complex" associated to the root group datum~$(U_\alpha)_{\alpha\in \Phi_{I(w)}}$ of the Levi subgroup~$L_{I(w)}$. Since~$f(W_{I(w)},I(w))<f(W,S)$,~$T'(\rho)$ is $(\kappa(W_{I(w)},I(w))-2)$-acyclic over~$\Q$ by induction. 

It is well-known that the join~$\Sigma_1\ast\Sigma_2$ of an $m_1$-acyclic simplicial complex~$\Sigma_1$ and an $m_2$-acyclic simplicial complex~$\Sigma_2$ is $(m_1+m_2+2)$-acyclic, where the previous acyclicity conditions are all taken over~$\Q$; see~\cite[\S 9.20]{bjorner}. Hence~$T(\rho)=\partial\rho\ast T'(\rho)$ is~$(|D_R(w)|+\kappa(W_{I(w)},I(w))-3)$-acyclic over~$\Q$. 
However, by~\eqref{eq:hmonot}, 
\begin{equation*}
\begin{split}
    |D_R(w)|+\kappa(W_{I(w)},I(w)) & =|S|-f(W_{I(w)},I(w))\\
    & \geq |S|-f(W,S)+1=\kappa(W,S)+1.
\end{split}
\end{equation*}
Therefore, $T(\rho)$ is $(\kappa(W,S)-2)$-acyclic over~$\Q$ for every~$\rho\in R_{d+1}$. 
By \Cref{lem:Deltadprop}, the filtration~$X_-=\bigcup_{d\geq 0}\Delta_d$ satisfies the hypotheses of \Cref{prop:filtHom} with~$m=\kappa(W,S)-2$. As~$X_-$ is contractible by Solomon--Tits theorem,~$\Delta_0$ (in fact, every~$\Delta_d$) is $(\kappa(W,S)-2)$-acyclic over~$\Q$. 
\end{proof}

\begin{rem}
    Notice that~\eqref{eq:thickSpec2} holds if
    \begin{equation}\label{eq:thickSpec3}
        \frac{\qmin}{\sqrt{\qmax}}>2(|S|-1)\cdot\max_{s,t\in S,\,s\neq t}\cos\Big(\frac{\pi}{m_{st}}\Big).
    \end{equation}
    Moreover, let~$(U_\alpha)_{\alpha\in \Phi}$ be the root group datum of Kac--Moody type over~$\F_q$ for~$G=\caG_A(\F_q)^\dagger$ (see \Cref{ex:KMgroups}) and assume that the associated Weyl group~$(W_A,S)$ is $2$-spherical. In this case, if~$q>3(|S|-1)^2$ then~\eqref{eq:thickSpec2} holds. Indeed, by~\eqref{eq:GCMcos} we have
    \begin{equation*}
        3(|S|-1)^2\geq(|S|-1)^2\cdot\max_{s,t\in S,\,s\neq t}A_{st}A_{ts}=\Bigg(2(|S|-1)\cdot\max_{s,t\in S,\,s\neq t}\cos\Big(\frac{\pi}{m_{st}}\Big)\Bigg)^2.
    \end{equation*} 
\end{rem}

\smallskip

The argument to prove \Cref{thm:Yacyc} can be transferred verbatim to the homotopic case once one has a homotopical analogue of \Cref{prop:YconnSph}. \cite[\S II.8, Corollary~19]{ab:twin} provides such an analogue if~$(W,S)$ is classical and for arbitrary large minimal thickness.
\begin{thm}\label{thm:Y-homot}
    Let~$G$ be a group with a locally finite root group datum~$(U_\alpha)_{\alpha\in \Phi}$ of classical $2$-spherical type~$(W,S)$. Assume~\hyperref[eq:co]{\co} and
    \begin{equation}\label{eq:Yhomot}
        \qmin\geq \max\{2^{\tau(W_J,J)\cdot|J|-1}\mid J\subseteq S\text{ spherical}\,\},
    \end{equation}
    where~$\tau(W_J,J)=1$ if~$(W_J,J)$ is of type~$A_n$ and~$\tau(W_J,J)=2$ otherwise. Then~$Y_-$ is $(\kappa(W,S)-2)$-connected.
\end{thm}
\begin{proof}
  The proof of \Cref{thm:Yacyc} transfers verbatim with the following adjustments:
  \begin{itemize}
      \item The base case of the induction is given by~\cite[\S II.8, Corollary~19]{ab:twin}.
      \item In the last paragraph,~\cite[\S 9.20]{bjorner} also covers that~$\Sigma_1\ast \Sigma_2$ is~$(m_1+m_2+2)$-connected if~$\Sigma_1$ is $m_1$-connected and~$\Sigma_2$ is $m_2$-connected.\qedhere
  \end{itemize}
\end{proof}

\subsection{Higher properties~T and groups with root group data}\label{sus:TkRGD}
Let~$G$ be a group with a root group datum~$(U_\alpha)_{\alpha\in \Phi}$ of type~$(W,S)$. For every~$J\subseteq S$, we set the following notation for the subsection:
\begin{itemize}
    \item[$\bullet$] Let~$(B_+,B_-,N)$ be the twin BN-pair of type~$(W,S)$ for~$G$ given in \Cref{thm:twinBNassoc}. Recall that the associated thickness vector~$(q_s,q_s)_{s\in S}$ is given by~$q_s^\pm=|U_s^\pm|$ for every~$s\in S$;
    \item[$\bullet$] $X_\pm=\Delta(G,B_\pm)$ is the simplicial building associated to the BN-pair~$(B_\pm,N)$ for~$G$;
    \item[$\bullet$] $Y_\pm\subseteq X_\pm$ is the negative opposite complex (see~\eqref{eq:Ypm}).
\end{itemize}
We also freely use \Cref{notat:RGD} and the notation for Coxeter systems introduced in \Cref{sus:cox}. Some results involve the cosine matrix~$A(V_0,\dots, V_n)$ introduced in \Cref{defn:CosMatRho}.

\begin{lem}\label{lem:qJineq}
    Let~$G$ be a group with a locally finite root group datum~$(U_\alpha)_{\alpha\in \Phi}$ of type~$(W,S)$. If
    \begin{equation*}
       \frac{\qmin}{\sqrt{\qmax}}>2(1-\veps_{(W,S)})
    \end{equation*}
    then, for every~$J\subsetneq S$, 
    \begin{equation*}
       \frac{\qminJ}{\sqrt{\qmaxJ}} >2(1-\veps_{(W_J,J)}).
    \end{equation*}
\end{lem}
\begin{proof}
    Let~$J\subsetneq S$ and denote by~$\frC$ and~$\frC_J$ the cosine matrices of~$(W,S)$ and~$(W_J,J)$, respectively.
   Note that there is an ordering of~$S$ such that~$\frC$ admits the following block-decomposition:
   \begin{equation*}
       \frC=\begin{bmatrix}
           \frC_J & \frC_{12}\\
           \frC_{21} & \frC_{22}
       \end{bmatrix}
   \end{equation*}
   for suitable matrices~$\frC_{12}$, $\frC_{21}$ and~$\frC_{22}$. 
  Since~$\frC$ is symmetric, by~\cite[Corollary~2.5.2]{bh:spectra} we also have
  \begin{equation*}
      \veps_{(W_J,J)}\geq \veps_{(W,S)}.
  \end{equation*}
 Finally, since~$\qmin\leq \qminJ\leq \qmaxJ\leq \qmax$, 
  \begin{equation*}
      \frac{\qmin}{\sqrt{\qmax}}\leq \frac{\qminJ}{\sqrt{\qmaxJ}}.\qedhere
  \end{equation*}
\end{proof}

\begin{prop}\label{prop:CMBN1}
   Let~$G$ be a group with a locally finite root group datum~$(U_\alpha)_{\alpha\in \Phi}$ of $2$-spherical type~$(W,S)$ with~$S=\{0,\dots, n\}$ for some~$n\geq 0$. Assume~\hyperref[eq:co]{\co} and
    \begin{equation}\label{eq:CosMatIneq}
        \frac{\qmin}{\sqrt{\qmax}}>2(1-\veps_{(W,S)}).
    \end{equation}
   Then, for every unitary representation~$(\rho,V)$ of~$U_+$, the cosine matrix~$A(V^{U_0},\dots, V^{U_n})$ is positive definite. 
\end{prop}
\begin{proof}
  Without loss of generality, $n\geq 1$. Let~$(\rho,V)$ be any unitary representation of~$U_+$ and let~$s,t\in S$ be such that~$s\neq t$. By \Cref{cor:UJfinite} and since~$(W,S)$ is $2$-spherical, $U_{\{s,t\}}$ is a finite group and clearly~$U_{\{s,t\}}=\langle U_s,U_t\rangle$.
    By \Cref{prop:GRO2}\eqref{prop:GRO2a} applied to the $U_+$-action on~$\caC\caC(U_+,(U_{S\setminus\{s\}})_{s\in S})$,
     \begin{equation}\label{eq:cosUineq}
         \cos\angle(V^{U_s},V^{U_t})\leq  \lambda_2\Big(\caC\caC\big(U_{\{s,t\}}, (U_s, U_t)\big)\Big)
     \end{equation}
  and \Cref{lem:spec1link} then implies
  \begin{equation}\label{eq:cosMatTBL1}
  \cos\angle(V^{U_s},V^{U_t})\leq 2\cdot\cos\Big(\frac{\pi}{m_{st}}\Big)\cdot\frac{\sqrt{\qmax}}{\qmin}.    
  \end{equation}
  This determines the following entry-wise inequality:
  \begin{equation*}
      A(V^{U_0},\dots, V^{U_n})\geq B:=\frac{2\sqrt{\qmax}}{\qmin}\frC_{(W,S)}+\Big(1-\frac{2\sqrt{\qmax}}{\qmin}\Big)I,
  \end{equation*}
  where~$\frC_{(W,S)}$ is the cosine matrix of~$(W,S)$ (see~\eqref{eq:cosmat}).
  From~\cite[Proposition~2.8]{gro:curv}, the smallest eigenvalue~$\alpha$ of~$A(V^{U_0},\dots, V^{U_n})$ is at least the smallest eigenvalue of~$B$, i.e., 
    \begin{equation*}
        \alpha\geq \frac{2\sqrt{\qmax}}{\qmin}\veps_{(W,S)}+1-\frac{2\sqrt{\qmax}}{\qmin}.
    \end{equation*}
By~\eqref{eq:CosMatIneq}, the right-hand side of the previous inequality is positive. Therefore, $A(V^{U_0},\dots, V^{U_n})$ is positive definite.
\end{proof}

A similar statement as in \Cref{prop:CMBN1} holds for every parabolic subgroups~$P_J^-$, but requires more work and will be obtained in \Cref{prop:CMBN2}. In this case we cannot proceed as in the proof of \Cref{prop:CMBN1}:~$P_{\{s,t\}}^-$ is not compact in general (it contains~$U_-$) and thus \Cref{prop:GRO2}\eqref{prop:GRO2a} cannot be applied directly. Before \Cref{prop:CMBN2}, we prove the following auxiliary result.
\begin{lem}\label{lem:Uspmt}
    Let~$G$ be a group with a locally finite root group datum~$(U_\alpha)_{\alpha\in \Phi}$ of $2$-spherical type~$(W,S)$ with~$S=\{0,\dots, n\}$ for some~$n\geq 1$, with finite subgroup~$T$ and such that~\hyperref[eq:co]{\co} holds. 
    For all~$s,t\in S$, consider the subgroup~$G_{s\pm t}:=T\langle U_s,U_t^{\pm}\rangle$ and the graph~$\Gamma_{s\pm t}:=\caC\caC(G_{s\pm t}, (TU_s, TU_t^\pm))$. Then
   \begin{equation*}
     \lambda_2(\Gamma_{s\pm t})\leq \mu_{s,t,\pm},\qquad  \mu_{s,t,\veps}:=\left\{
       \begin{array}{cl}
          2\cdot\cos\big(\frac{\pi}{m_{st}}\big)\cdot\frac{\sqrt{\qmax}}{\qmin}  &  \text{if $s\neq t$ and~$\veps=1$},\\
          0  &  \text{if $s\neq t$ and~$\veps=-1$},\\
          q_s^{-1} & \text{if~$s=t$ and~$\veps=-1$.}
       \end{array}
       \right.
   \end{equation*}
\end{lem}
Recall from \Cref{rem:Tnorm} that~$G_{s\pm t}$, $TU_s^\pm$ and~$TU_t^\pm$ are subgroups of~$G$.

\begin{proof}
We first show that~$G_{s\pm t}$ is finite. 
In detail, $G_{s\pm t}\subseteq TL_{\{s,t\}}$ and~$T$ is finite by hypothesis. Also the Levi subgroup~$L_{\{s,t\}}$ is finite: by \Cref{lem:RGDsubsyst} and \Cref{thm:twinBNassoc}, it admits a BN-pair of spherical type~$(W_{\{s,t\}},\{s,t\})$ and thus its "$B$-group", which is the finite group~$TU_{\{s,t\}}$ by \Cref{cor:UJfinite}, has finite index in~$L_{\{s,t\}}$. 

Thus by \Cref{prop:GRO2}\eqref{prop:GRO2a} we have
\begin{equation*}
    \cos\angle (V^{TU_s},V^{TU_t^\pm})\leq \lambda_2(\Gamma_{s\pm t}).
\end{equation*}
It remains to show that
$$\lambda_2(\Gamma_{s\pm t})\leq \mu_{s,t,\pm}.$$ 
Arguing as to prove~\eqref{eq:cosMatTBL1}, we have~$\lambda_2(\Gamma_{s+t})\leq \mu_{s,t,+}$ for all~$s\neq t$.

We now show that, for all~$s\neq t$,
$$\lambda_2(\Gamma_{s-t})=0.$$
By \Cref{ex:openInt} and~\hyperref[RGD3]{(RGD3)} we have~$[U_s,U_t^-]=\{1\}$. In particular, $G_{s-t}=\langle U_s,U_t^{-}\rangle T=U_sU_t^{-}T$. We deduce that~$\Gamma_{s-t}$ is a complete bipartite graph (with its defining vertex bipartition), i.e., for all~$x,x'\in G_{s-t}$ we have~$xTU_s\cap x'TU_t^-\neq \emptyset$. In detail, write~$x=uvh$ and~$x'=u'v'h'$ for some~$u,u'\in U_s$, $v,v'\in U_t^-$ and~$h,h'\in T$. Since~$T$ normalises both~$U_s$ and~$U_t^-$ (see \Cref{thm:twinBNassoc}\eqref{thm:twinBN4}) and~$[u,v]=[u',v']=1$, we get~$xTU_s=uvTU_s=vuU_sT=vU_sT$ and~$x'TU_t^{-}=u'v'TU_t^-=u'v'U_t^-T=u'U_t^-T$ and~$vu'=u'v\in xTU_s\cap x'TU_t^-$. 
This proves that~$\Gamma_{s-t}$ is a complete bipartite graph, say~$\Gamma_{s-t}\cong K_{a,b}$ for some integers~$a,b\geq 1$. 
Moreover, $\Gamma_{s-t}$ is~$(q_s,q_t)$-regular: indeed, since~$B_+\cap U_-=U_+\cap B_-=\{1\}$ and~$T=B_+\cap B_-$ by \Cref{thm:twinBNassoc}, we deduce that~$TU_s\cap TU_t^-=T$, $U_s\cap T=U_t^-\cap T=\{1\}$ and
\begin{equation*}
    |TU_s:TU_s\cap TU_t^-|=|TU_s:T|=q_s\quad\text{and}\quad|TU_t^-:TU_s\cap TU_t^-|=|TU_t^-:T|=q_t.
\end{equation*}
Notice that~$|\Gamma_{s-t}(0)|=a+b\geq q_s+1\geq 3$ because~$U_s$ is non-trivial by~\hyperref[RGD0]{(RGD0)}.
By~\eqref{eq:specCompar} and~\cite[\S 1.4.2]{bh:spectra}, the normalised weighted adjacency matrix~$G_{\Gamma_{s-t}}$ has spectrum
\begin{equation*}
    \Spec(G_{\Gamma_{s-t}})=\textstyle{\frac{1}{\sqrt{q_sq_t}}}\Spec(A_{K_{a,b}})=\Big\{\pm\textstyle{\sqrt{\frac{ab}{q_sq_t}}},0\Big\}.
\end{equation*}
Having~$a+b\geq 3$ is important to guarantee that~$0$ is an eigenvalue of~$A_{K_{a,b}}$ and hence of~$G_{K_{a,b}}$.
This yields~$\lambda_2(\Gamma_{s-t})=0$.

We conclude the argument by showing that, for every~$s\in S$,
$$\lambda_2(\Gamma_{s-s})=q_s^{-1}.$$
To this end, we prove that~$\Gamma_{s-s}$ is a graph of the form given in \Cref{ex:Bip}.
Note that~$G_{s-s}=\langle U_s,U_s^-\rangle T=L_sT$ and, by \Cref{lem:RGDsubsyst} and \Cref{thm:twinBNassoc}, $L_s$ has a twin BN-pair of type~$(W_s,\{s\})$ and~"$B_\pm$-group"~$B_\pm^s:=U_s^\pm T_s$ for some subgroup~$T_s\leq T$. Thus, by \Cref{cor:refBirk} we have~$L_s=B_\pm^s\sqcup U_s sB_\pm^s$.
For simplicity, let~$H_\pm:=U_s^\pm T$. Since~$B_\pm^sT=H_\pm$, we obtain
\begin{equation}\label{eq:Gs-sdec}
   G_{s-s}=L_sT=H_\pm\sqcup U_s^\pm s H_\pm.
\end{equation}
The union before remains disjoint because~$H_\pm\subseteq B_\pm$ and~$U_s^\pm sH_\pm\subseteq B_\pm s B_\pm$ and clearly~$B_\pm\cap B_\pm sB_\pm=\emptyset$.
Let~$q:=q_s$ and~$U_s=\{u_1,\dots, u_q\}$.
As~$s\in L_s$ we have~$G_{s-s}=sG_{s-s}$. Moreover,~$sU_s^-s=U_s^+$ by \Cref{rem:Tnorm}. Thus, by~\eqref{eq:Gs-sdec} and~$G_{s-s}=H_+\sqcup U_ssH_+=sH_-\sqcup U_sH_-$, the two sides of the defining bipartition of~$\Gamma_{s-s}(0)$ have the form:
\begin{equation*}
    G_{s-s}/H_+=\{1H_+\}\sqcup\{u_isH_+\mid 1\leq i\leq q\}\quad\text{and}\quad G_{s-s}/H_-=\{sH_-\}\sqcup \{u_iH_-\mid 1\leq i\leq q\}.
\end{equation*}
Since~$U_s\cap sH_+s^{-1}=U_s\cap H_-\subseteq U_+\cap B_-=\{1\}$ (see \Cref{thm:twinBNassoc}), we get~$|G_{s-s}:H_\pm|=q+1$. 
It is immediate to check that, in~$\Gamma_{s-s}$, the vertex~$1H_+$ is connected to all vertices of~$G_{s-s}/H_-$ except for~$sH_-$. Similarly, for every~$1\leq i\leq q$, the vertex~$u_isH_+$ is connected to all vertices of~$G_{s-s}/H_-$ except for~$u_issH_-=u_iH_-$. Ordering~$\Gamma_{s-s}(0)$ as
\begin{equation*}
    1H_+<u_1sH_+<\dots<u_qsH_+<sH_-<u_1H_-<\dots<u_qH_-,
\end{equation*}
the adjacency matrix of~$\Gamma_{s-s}$ admits the following decomposition in $q$-dimensional square blocks:
\begin{equation*}
    A_{\Gamma_{s-s}}=\begin{bmatrix}
        0 & \bbone-I\\
        \bbone-I & 0
    \end{bmatrix}.
\end{equation*}
Now \Cref{ex:Bip} yields~$\lambda_2(\Gamma_{s-s})=q^{-1}=q_s^{-1}$.
\end{proof}

\begin{prop}\label{prop:CMBN2}
Let~$G$ be a group with a locally finite root group datum~$(U_\alpha)_{\alpha\in \Phi}$ of $2$-spherical type~$(W,S)$ with~$S=\{0,\dots, n\}$, $n\geq 1$, and finite subgroup~$T$. Assume~\hyperref[eq:co]{\co} and 
\begin{equation}\label{eq:CM2Ineq}
    \frac{\qmin-1}{\sqrt{\qmax}}>2(1-\veps_{(W,S)}).
\end{equation}
Then, for every non-empty set~$J=\{j_0,\dots, j_\ell\}\subseteq S$ and every unitary representation~$(\rho,V)$ of~$P_J^-$, the cosine matrix~$A(V^{P_{j_0}^-},\dots, V^{P_{j_\ell}^-})$ is positive definite.
\end{prop}

\begin{proof}
To ease the notation, let~$J=\{0,\dots, \ell\}$ and, without loss of generality, assume~$\ell\geq 1$. Let also~$(\rho,V)$ be a unitary representation of~$P_J^-$. 

By \Cref{cor:weakLevi} and \Cref{cor:UJUgroup}, for every~$s\in S$ we have
\begin{equation*}
    P_s^{-}=B_-L_s=TU_-L_s=T\langle U_x^-\mid x\in S\rangle\langle U_s^+,U_s^-\rangle=\langle TU_s,TU_x^-\mid x\in S\rangle.
\end{equation*}
Setting~$V_x^\pm:=V^{TU_x^\pm}$ for every~$x\in S$, we deduce that
\begin{equation*}
    V^{P_s^-}=V_s^+\cap \bigcap_{x\in S}V_x^-.
\end{equation*}

Recall now~\cite[Corollary~2.12]{gro:curv}: for all closed subspaces~$H_0,\dots, H_n$ of a given Hilbert space, if the cosine matrix~$A(H_0,\dots, H_n)$ is positive definite, so is~$A(H_0\cap H_n,\dots, H_{n-1}\cap H_n)$.
Applying this fact repeatedly, we deduce that if
\begin{equation}\label{eq:caAJ}
    \caA_J:=A(V_0^+,\dots, V_\ell^+,V_0^-,\dots, V_n^-)
\end{equation}
is positive definite, then so is
\begin{equation*}
   A(V^{P_0^-},\dots, V^{P_\ell^-})=A\big(V_0^+\cap \textstyle{\bigcap_{x\in S}}V_x^-,\dots,V_\ell^+\cap \textstyle{\bigcap_{x\in S}}V_x^-\big).
\end{equation*}

It remains to prove that the matrix~$\caA_J$ in~\eqref{eq:caAJ} is positive definite. 
Consider the matrix~$\caB$ given by the following decomposition in square $(n+1)$-dimensional blocks:
\begin{equation*}
    \caB=[B_{ij}]_{0\leq i,j\leq 2n+1}=\begin{bmatrix}
        \caB_0 & -\qmin^{-1}I\\
        -\qmin^{-1}I & \caB_0
    \end{bmatrix},\qquad \caB_0=\frac{2\sqrt{\qmax}}{\qmin}\frC_{(W,S)}+\Big(1-\frac{2\sqrt{\qmax}}{\qmin}\Big)I.
\end{equation*}
By \Cref{lem:Uspmt}, we have an entry-wise inequality~$\caA_J\geq \caB_J$, where~$\caB_J:=[B_{ij}]_{i,j\in\{0,\dots, \ell\}\sqcup\{n+1,\dots, 2n+1\}}$.
If we prove that~$\caB$ is positive definite, also the principal submatrix~$\caB_J$ is positive definite and~\cite[Proposition~2.8]{gro:curv} then yields~$\caA_J$ is positive definite.

It remains to prove that~$\caB$ is positive definite. By~\eqref{eq:CM2Ineq}, the smallest eigenvalue~$\beta$ of~$\caB_0$ satisfies
\begin{equation}\label{eq:betaineq1}
    \beta=\frac{2\sqrt{\qmax}}{\qmin}\veps_{(W,S)}+1-\frac{2\sqrt{\qmax}}{\qmin}>\frac{1}{\qmin}.
\end{equation}
We now claim that
\begin{equation}\label{eq:betaineq2}
    \Spec(\caB)\subseteq \big\{\mu+\textstyle{\frac{1}{\qmin}}, \mu-\textstyle{\frac{1}{\qmin}}, \mu\,\big\vert\, \mu\in \Spec(\caB_0)\big\}.
\end{equation}
By~\eqref{eq:betaineq1}, proving~\eqref{eq:betaineq2} implies that~$\Spec(\caB)\subseteq \R_{>0}$ and then~$\caB$ is positive definite.
It remains to prove~\eqref{eq:betaineq2}. Using Schur's formula~\eqref{eq:schur}, for every~$\lambda\in \R$ we have
\begin{equation*}
\begin{split}
    \det(\caB-\lambda I)=0 & \Longrightarrow\, \det(\caB_0-\lambda I)=0 \text{ or }\det\big((\caB_0-\lambda I)^2-\qmin^{-2}I\big)=0\\
    & \Longleftrightarrow\, \lambda\in\Spec(\caB_0)\text{ or }\det(\caB_0-(\lambda+\qmin^{-1})I)=0\text{ or }\det(\caB_0-(\lambda-\qmin^{-1})I)=0
    \end{split}
\end{equation*}
and~\eqref{eq:betaineq2} follows.
\end{proof}

\begin{prop}\label{thm:U+propTk}
   Let~$G$ be a group with a locally finite root group datum~$(U_\alpha)_{\alpha\in \Phi}$ of type~$(W,S)$ that satisfies~\hyperref[eq:co]{\co}. Assume~$(W,S)$ is $(k+1)$-spherical for some~$1\leq k\leq |S|-1$ and
\begin{equation}\label{eq:Cond}
    \frac{\qmin}{\sqrt{\qmax}}>2\cdot \max\left\{1-\veps_{(W,S)},\,(\sigma(W,S)-1)\cdot\max_{s,t\in S,\, s\neq t}\Big\{\cos\Big(\frac{\pi}{m_{st}}\Big)\Big\}\right\}
\end{equation}
  Then~$U_+$ has property~$\Tk$.
\end{prop}

\begin{rem}\label{rem:ineqPassJ}
   Notice that~\eqref{eq:Cond} is "stable under passage to a subset~$J\subseteq S$". By \Cref{lem:qJineq}, \eqref{eq:Cond} implies the following for every~$J\subseteq S$:
    \begin{equation*}
     \frac{\qminJ}{\sqrt{\qmaxJ}}>2\cdot \max\left\{1-\veps_{(W_J,J)},\,(\sigma(W_J,J)-1)\cdot\max_{s,t\in J,\, s\neq t}\Big\{\cos\Big(\frac{\pi}{m_{st}}\Big)\Big\}\right\}.
\end{equation*}
\end{rem}

\begin{proof}
If~$(W,S)$ is spherical, then~$U_+$ is finite by \Cref{cor:UJfinite} and hence it has property~$\Tk$ by \Cref{thm:Tkperm}\eqref{thm:Tkperm1}.

We can now prove the statement by induction on~$|S|\geq k+1$. The base case~$|S|=k+1$ is already proved because~$(W,S)$ is spherical.
Let now~$n\geq k+1$ and assume the claim true for all root group data as in the hypotheses and whose type~$(W',S')$ satisfies~$k+1\leq |S'|\leq n$. Let~$|S|=n+1$, say~$S=\{0,\dots, n\}$.

Our goal is to apply \Cref{prop:cohoPlainCC} to the $U_+$-action on the coset complex~$Z:=\caC\caC(U_+, (U_{S\setminus\{s\}})_{s\in S})$. 
By \Cref{cor:Y-Uchar} and \Cref{thm:Yacyc}, $Z\cong Y_-$ and hence~$Z$ is $(k-1)$-acyclic over~$\C$. Moreover, the $U_+$-action on~$Z$ has a strict simplicial fundamental domain that corresponds to the face complex of the chamber~$\{1U_{S\setminus\{s\}}\mid s\in S\}$ of~$Z$. Every non-empty simplex of~$Z$ has the form~$\{uU_{S\setminus\{s\}}\mid s\in J\}$ for some non-empty set~$J\subseteq S$ and some~$u\in U_+$, and thus its stabiliser in~$U_+$ is conjugated to~$\bigcap_{s\in J}U_{S\setminus\{s\}}=U_{S\setminus J}$ (see \Cref{eq:intersParab} and \Cref{lem:UJUPinters}). By \Cref{lem:RGDsubsyst} and \Cref{cor:UJUgroup}, $U_{S\setminus J}$ is the "$U_+$-group" of the root group datum~$(U_\alpha)_{\alpha\in \Phi_{S\setminus J}}$ for the Levi subgroup~$L_{S\setminus J}$. From \Cref{rem:ineqPassJ}, the previous root group datum satisfies the hypotheses of the proposition. Hence, by induction, $U_{S\setminus J}$ has property~$\Tk$.
Finally, by \Cref{prop:CMBN1}, the condition~\eqref{eq:Cond} guarantees that the cosine matrix~$A(V^{U_0},\dots, V^{U_n})$ is positive definite for every unitary representation~$(\rho,V)$ of~$U_+$. By \Cref{prop:cohoPlainCC},~$U_+$ has property~$\Tk$. 
\end{proof}

In fact, a refinement of \Cref{thm:U+propTk} (see \Cref{thm:U+nonTk}) will be obtained as a consequence of the following theorem:

\begin{thm}\label{thm:RGDpropTk}
Let~$G$ be a group with a locally finite root group datum~$(U_\alpha)_{\alpha\in \Phi}$ of type~$(W,S)$ that satisfies~\hyperref[eq:co]{\co} and such that the subgroup~$T$ as in~\eqref{eq:TBN} is finite.
Assume~$(W,S)$ is $(k+1)$-spherical for some~$k\geq 1$, and suppose that
\begin{equation}\label{eq:Cond'}
    \frac{\qmin-1}{\sqrt{\qmax}}>2\cdot \max\left\{1-\veps_{(W,S)},\,(\sigma(W,S)-1)\cdot\max_{s,t\in S,\, s\neq t}\Big\{\cos\Big(\frac{\pi}{m_{st}}\Big)\Big\}\right\}
\end{equation}
Then:
    \begin{enumerate}
        \item \label{thm:RGDTk1} For every subset~$J\subseteq S$, the parabolic subgroups~$P_J^\pm$ have property~$\Tk$. In particular, $G$ has property~$\Tk$.
        \item \label{thm:RGDTk2} Assume additionally that~$(W,S)$ is not $(k+2)$-spherical and that~$\qmin>|S|$. Then, for every non-spherical subset~$J\subseteq S$ with~$|J|=k+2$, the parabolic subgroup~$P_J^\pm$ has property~$\Tk$ but not property~$(T_{k+1})$. 
    \end{enumerate}
\end{thm}
\begin{proof}
By \Cref{rem:defnRGDsym}, it suffices to prove the statement for the positive parabolic subgroups.

The BN-pair~$(B_+,N)$ for~$G$ as in \Cref{thm:twinBNassoc} induces, for every~$J\subseteq S$, a BN-pair~$(B_+, N_J)$ for~$P_J^+$ of type~$(W_J,J)$ (see \Cref{sus:BNpairs}). The minimal and maximal thicknesses of~$(B_+,N_J)$ are~$\qminJ$ and~$\qmaxJ$, respectively. Let~$\Delta_+^J$ be the fundamental chamber of the building~$X_+^J=\Delta(P_J^+,B_+)$.
Let~$J=\{j_0,\dots, j_\ell\}$.
By \Cref{prop:CMBN2} applied to~$P_J^+$ (recall again \Cref{rem:defnRGDsym}), for every unitary representation~$(\rho,V)$ of~$P_J^+$ the cosine matrix~$A(V_{j_0},\dots, V_{j_\ell})$, where~$V_{j_i}=V^{P_{j_i}^+}$ for every~$i$, is positive definite. 
    
The claims now follow from \Cref{cor:vanishBN} applied to each BN-pair~$(B_+,N_J)$, once we prove that~$B_+$ has property~$\Tk$. By \Cref{thm:twinBNassoc}\eqref{thm:twinBN1}, $B_+\cong U_+\rtimes T$. Since~$T$ is finite by hypothesis and~$U_+$ has property~$\Tk$ by \Cref{thm:U+propTk}, $B_+$ has property~$\Tk$ by \Cref{thm:Tkperm}\eqref{thm:Tkperm3}. 
\end{proof}

\begin{cor}\label{cor:KMpropTk}
    Let~$G=\caG_A(\F_q)$ be a minimal Kac--Moody group over a finite field~$\F_q$ and generalised Cartan matrix~$A=[A_{st}]_{s,t\in S}$. Assume that the associated root group datum satisfies~\hyperref[eq:co]{\co} and that the Weyl group~$(W_A,S)$ associated to~$A$ is $(k+1)$-spherical for some~$k\geq 1$. Suppose also that
    \begin{equation}\label{eq:suffKM}
     \sqrt{q}>1+\max\left\{2(1-\veps_{(W_A,S)}),\,(\sigma(W_A,S)-1)\cdot\sqrt{\max_{s,t\in S,\,s\neq t} A_{st}A_{ts}}\right\}.
    \end{equation}
    Then~$G$ has property~$\Tk$. If in particular~$|S|=k+2$, $(W_A,S)$ is non-spherical and~$q>|S|$, then~$G$ has property~$\Tk$ but not~$(T_{k+1})$.
\end{cor}
\begin{proof}
    Following the notation in \Cref{ex:KMgroups}, the subgroup~$H:=\caG_A(\F_q)^\dagger=\langle U_\alpha\mid\alpha\in \Phi_A\rangle$ has finite index in~$G=\caG_A(\F_q)$. By \Cref{thm:Tkperm}\eqref{thm:Tkperm3}, it suffices to prove the claim for~$H$. But this is now a straightforward consequence of \Cref{thm:RGDpropTk} and~\eqref{eq:GCMcos}.
\end{proof}

\begin{rem}\label{rem:KMpropTk}
    All the conclusions of \Cref{cor:KMpropTk} hold if, for instance, 
    \begin{equation}\label{eq:KMcondEasy}
    q>3(|S|-1)^2+1.
    \end{equation}
    Indeed, the claim is true if~$(W_A,S)$ is spherical even without assuming~\eqref{eq:suffKM}. Let now~$(W_A,S)$ be non-spherical and $2$-spherical. By~\eqref{eq:GCMcos} we have
    \begin{equation*}
        \sqrt{3}(|S|-1)>(\sigma(W_A,S)-1)\cdot\sqrt{\max_{s,t\in S,\,s\neq t}A_{st}A_{ts}}.
    \end{equation*}
    Moreover, $1-\veps_{(W_A,S)}$ is the largest eigenvalue of the matrix~$I-\frC_{(W,S)}$. The latter matrix has non-negative real entries and vanishing diagonal. By Gershgorin circle theorem~\cite[Proposition~2.10.2]{bh:spectra} and since~$m_{st}\leq 6$ for all~$s,t\in S$ (see~\eqref{eq:GCMcos}), we have
    \begin{equation*}
        1-\veps_{(W_A,S)}\leq \max_{s\in S}\sum_{t\in S,\,s\neq t}\cos\Big(\frac{\pi}{m_{st}}\Big)\leq \frac{\sqrt{3}}{2}(|S|-1)
    \end{equation*}
    Finally, if~$q>3(|S|-1)^2+1$ and~$(W_A,S)$ is non-spherical, then~$|S|\geq 3$ and~$q>3(|S|-1)^2+1>|S|$. In particular~$q\geq 4$, and~\hyperref[eq:co]{\co} holds.
\end{rem}

\begin{cor}\label{cor:twinTk}
    Let~$\Lambda$ be a twin building lattice of $(k+1)$-spherical type~$(W,S)$ whose associated root group datum~$(G,(U_\alpha)_{\alpha\in \Phi})$ satisfies~\hyperref[eq:co]{\co} and~\eqref{eq:Cond'}.
    
    Then $\Lambda$ has property~$\Tk$.
    If in particular~$|S|=k+2$, $(W,S)$ is non-spherical and~$\qmin>|S|$, then~$\Lambda$ has property~$\Tk$ but not property~$(T_{k+1})$.
\end{cor}
\begin{proof}
     By \Cref{rem:TBL},~$\Lambda$ is isomorphic to the quotient~$G/Z(G)$ and~$Z(G)$ is finite. The claim now follows from \Cref{thm:Tkperm}\eqref{thm:Tkperm2} and \Cref{thm:RGDpropTk}.
\end{proof}

\begin{thm}\label{thm:U+nonTk}
 Let~$G$ be a group with a locally finite root group datum~$(U_\alpha)_{\alpha\in \Phi}$ of type~$(W,S)$, with finite subgroup~$T$ (see~\eqref{eq:TBN}) and such that~\hyperref[eq:co]{\co},~\eqref{eq:Cond} and~$\qmin>|S|$ hold. Assume~$(W,S)$ is $(k+1)$-spherical but not~$(k+2)$-spherical for some~$k\geq 1$. 
Then the group~$U_+$ has property~$\Tk$ but not~$[T_{k+1}]$.
\end{thm}
\begin{proof}
By \Cref{thm:U+propTk}, $U_+$ has property~$\Tk$. 
 Let now~$J\subseteq S$ be a non-spherical subset such that~$|J|=k+2$. By \Cref{thm:RGDpropTk}\eqref{thm:RGDTk2}, $P_J^+$ does not have property~$[T_{k+1}]$. On the other hand, since~$B_+=U_+\rtimes T$ and by \Cref{thm:Tkperm}\eqref{thm:Tkperm3}, $U_+$ has property~$[T_{k+1}]$ if and only if so has~$B_+$. 
 
 We assume for a contradiction that~$B_+$ has property~$[T_{k+1}]$. 
 Let~$(B_+,N_J)$ be the BN-pair of type~$(W_J,J)$ for~$P_J^+$ induced by~$(B_+,N)$, and notice that it has finite thickness vector~$(q_s)_{s\in J}$.
 By \Cref{rem:GBspher}, for every non-empty subset~$I\subseteq J$, the group $P_{J\setminus I}$ contains~$B_+$ as a finite-index subgroup and hence has property~$[T_{k+1}]$ by \Cref{thm:Tkperm}\eqref{thm:Tkperm3}. By \Cref{cor:cohoPlainBN} and \Cref{rem:ineqPassJ}, if~$B_+$ has property~$[T_{k+1}]$ then~$P_J^+$ has property~$[T_{k+1}]$, impossible.
\end{proof}
   
\Cref{thm:U+nonTk} implies the following statement on Kac--Moody--Steinberg groups (see \Cref{ex:KMgroups}).
\begin{cor}\label{cor:KMSTk}
    Let~$A=[A_{st}]_{s,t\in S}$ be a generalised Cartan matrix of dimension~$|S|\geq 2$ and associated Weyl group~$(W_A,S)$ of $(k+1)$-spherical type for some~$k\geq 2$. Let also~$q$ be a prime power satisfying
    \begin{equation*}
        \sqrt{q}>1+\max\left\{2,\,2(1-\veps_{(W,S)}),\,(\sigma(W_A,S)-1)\cdot\sqrt{\max_{s,t\in S,\,s\neq t} A_{st}A_{ts}}\right\}.
    \end{equation*}
    Then the Kac--Moody--Steinberg group~$\caU_A(\F_q)$ has property~$\Tk$. If~$(W,S)$ is additionally not $(k+2)$-spherical and~$q>|S|$, then~$\caU_A(\F_q)$ does not have property~$[T_{k+1}]$.
\end{cor}

\subsection{Finiteness properties and groups with root group data}\label{sus:FinPropRGD}

In this section, we prove the results stated in \Cref{sus:introFP}. Before this, we recall some generalities on finiteness properties that will be used repeatedly in the following proofs. 

Let~$G$ be a group and~$R$ be a unital ring.
For every~$n\in \Z_{\geq 0}$, an \emph{$n$-good $G$-CW complex} over~$R$ is a CW complex~$X$ which is $(n-1)$-acyclic over~$R$ and that has a cellular $G$-action such that, for every cell~$\sigma$ with~$0\leq \dim(\sigma)\leq n-1$, $G_\sigma$ is of type~$\FP_{n-\dim(\sigma)}(R)$.

The following characterisation is commonly known as "Brown's criterion".
\begin{thm}[\protect{\cite[Theorem~2.2]{BrownFP}}]\label{thm:BroCrit}
    Let~$X$ be an $n$-good $G$-CW complex over~$R$ with a filtration~$X=\bigcup_{\alpha\in I}X_\alpha$ over a directed set~$(I,\leq)$ such that the $n$-skeleton of each~$X_\alpha$ is finite modulo~$G$. Then~$G$ is of type~$\FP_n(R)$ if and only if, for every~$0\leq j\leq n-1$, the sequence of reduced homology groups~$(\wtil{\hH}_j(X_\alpha,R))_{\alpha\in I}$ is \emph{essentially trivial}, i.e., for all~$\alpha\in I$ there is~$\beta\geq \alpha$ such that the map~$\wtil{\hH}_j(X_\alpha,R)\to \wtil{\hH}_j(X_\beta,R)$ induced by the inclusion~$X_\alpha\hookrightarrow X_\beta$ is the zero map.
\end{thm}
An immediate consequence of \Cref{thm:BroCrit} is the following (see also~\cite[Proposition~1.1]{brown:coho}):
\begin{cor}\label{prop:BrownSuff}
    Let~$X$ be an $n$-good $G$-CW complex over~$R$ such that its $n$-skeleton is finite modulo~$G$. Then~$G$ is of type~$\FP_n(R)$.  
\end{cor}

For type~$\tF_k$, we recall the following well-known criteria.
\begin{prop}\label{prop:FPnFn}
Let~$G$ be a group.
\begin{enumerate}
    \item \label{prop:FPnFn1} Assume that~$G$ has a cellular action on a $1$-connected CW complex~$X$ with finitely presented vertex stabilisers, finitely generated edge stabilisers, and with finitely many orbits on the $2$-skeleton of~$X$. Then~$G$ is finitely presented.
    \item \label{prop:FPnFn2} Let~$G$ be a finitely presented group and let~$n\in \Z_{\geq 2}$. Then~$G$ is of type~$\FP_n(\Z)$ if and only if~$G$ is of type~$\tF_n$.
\end{enumerate}
\end{prop}
\begin{proof}
    Part~\eqref{prop:FPnFn1} is~\cite[Proposition~3.1]{brown:coho}. Regarding part~\eqref{prop:FPnFn2}: for the "if part", consider the projective resolution given by the augmented cellular chain complex of a model for~$EG$ over~$\Z$; the "only if" part is recorded, with an explicit proof, in~\cite[Lemma~2.1]{fffza}.
\end{proof}

\begin{prop}[\protect{\cite[\S 2.2, Propositions~2.5 and~2.7]{bieri}}]\label{prop:permFP}
   Let~$G$ be a group. Then:
    \begin{enumerate}
        \item For every finite-index subgroup~$H\leq G$, $H$ is of type~$\FP_n(R)$ if and only if so is~$G$.
        \item For every finite normal subgroup~$K\unlhd G$, $G/K$ is of type~$\FP_n(R)$ if and only if so is~$G$.
    \end{enumerate}
    The same statements hold also for~$\tF_n$.
\end{prop}

We now focus on the finiteness properties for a group~$G$ with locally finite root group data~$(U_\alpha)_{\alpha\in \Phi}$. Throughout, $(B_+,B_-,N)$ is the twin BN-pair associated to~$(G,(U_\alpha)_{\alpha\in \Phi})$ (see \Cref{thm:twinBNassoc}). We set~$T=B_+\cap B_-$ (see also~\eqref{eq:TBN}) and~$X_\pm=\Delta(G,B_\pm)$, and denote by~$Y_-\subseteq X_-$ the negative opposite complex (see~\eqref{eq:Ypm}).

\begin{thm}\label{thm:UFP}
Let~$G$ be a group with a locally finite root group datum~$(U_\alpha)_{\alpha\in \Phi}$ of $(k+1)$-spherical type~$(W,S)$ for some~$k\geq 1$, and such that~\hyperref[eq:co]{\co} and~\eqref{eq:thickSpec2} hold.
Then:
\begin{enumerate}
    \item \label{thm:UFP1} $U_+$ is of type~$\FP_k(\Q)$.
    \item \label{thm:UFP2} If~$|S|=k+2$ and~$(W,S)$ is not spherical, then~$U_+$ is not of type~$\FP_{k+1}(\Q)$.
    \item \label{thm:UFP3} If~$(W,S)$ is classical and~$\qmin$ satisfies~\eqref{eq:Yhomot}, then~$U_+$ is of type~$\tF_k$.
\end{enumerate}
\end{thm}
\begin{proof}
   Let~$(W,S)$ be $(k+1)$-spherical. If~$(W,S)$ is spherical, then~$U_+$ is finite (see \Cref{cor:UJfinite}) and there is nothing to prove.
   From now on, let~$(W,S)$ be non-spherical and denote by~$Y_-$ the negative opposite complex (see~\eqref{eq:Ypm}).

\eqref{thm:UFP1} By \Cref{thm:Yacyc}, $Y_-$ is $(k-1)$-acyclic over~$\Q$. By \Cref{cor:Y-Uchar}, $Y_-$ is $U_+$-equivariant isomorphic to the coset complex~$\caC\caC(U,(U_{S\setminus\{s\}})_{s\in S})$. Hence~$U_+$ has finitely many orbits on~$Y_-$. Moreover, for every~$\ell\geq 0$, the pointwise stabiliser in~$U_+$ of any simplex~$\sigma\in Y_-(\ell)$ is conjugated to~$U_{S\setminus I}$ for some~$I\subseteq S$ with~$|I|=\ell+1$ (see \Cref{lem:UJUPinters} and~\eqref{eq:intersParab}). 
  By \Cref{lem:RGDsubsyst}, for every~$J\subseteq S$, $U_{J}$ is the "$U_+$-group" of the root group datum~$(U_\alpha)_{\alpha\in \Phi_{J}}$ of type~$(W_{J},J)$, and thus it satisfies the hypotheses of the theorem. By \Cref{cor:UJfinite}, whenever~$J$ is spherical, $U_{J}$ is finite and hence of type~$\FP_k(\Q)$.  Moreover, by~\Cref{prop:BrownSuff}, if~$U_{S\setminus I}$ is of type~$\FP_k(\Q)$ for all non-empty subsets~$I\subseteq S$, then~$U_+$ is of type~$\FP_k(\Q)$.
  The claim now follows by induction on~$|S|\geq k+2$.
  
\eqref{thm:UFP3} As~$\tF_1$ is equivalent to~$\FP_1(\Q)$, we may assume~$k\geq 2$.
By \Cref{thm:Y-homot}, $Y_-$ is $(k-1)$-connected and hence $(k-1)$-acyclic over~$\Z$. By \Cref{prop:FPnFn}\eqref{prop:FPnFn2}, it suffices to prove that~$U_+$ is finitely presented and of type~$\FP_k(\Z)$.
In both cases, the argument used to prove part~\eqref{thm:UFP1} can be transferred verbatim with only one minor change: to prove that~$U_+$ is finitely presented, in the induction step, apply \Cref{prop:FPnFn}\eqref{prop:FPnFn1} instead of \Cref{prop:BrownSuff}.
\end{proof}

\begin{thm}\label{thm:nonFPnU+}
    Let~$G$ be a group with a locally finite root group datum~$(U_\alpha)_{\alpha\in \Phi}$ of non-spherical $(k+1)$-spherical type~$(W,S)$, where~$|S|=k+2\geq 2$, and such that~\hyperref[eq:co]{\co} and~\eqref{eq:thickSpec2} hold. Then~$U_+$ is of type~$\FP_k(\Q)$ but not of type~$\FP_{k+1}(\Q)$.
\end{thm}
\begin{proof}
    The positive part of the statement is due to \Cref{thm:UFP}.
For the negative part, we follow the proof strategy of~\cite[\S I.6, Theorem~A]{ab:twin}. Nevertheless, we provide a self-contained proof for completeness. 

Let~$X_-=\Delta(G,B_-)$ and notice that $U_+$ acts simplicially on it. We aim to apply \Cref{thm:BroCrit} to Abels' filtration~$X_-=\bigcup_{d\geq 0}\Delta_d$ (see \Cref{defn:Abfiltr}).
Firstly, by Solomon--Tits theorem, $X_-$ is contractible and hence acyclic over~$\Q$. 

We now prove that the $U_+$-stabiliser of every non-empty simplex in~$X_-$ is finite.
By \Cref{cor:refBirk}, every simplex of~$X_-$ of type~$I\subseteq S$ is of the form~$\sigma(uw,I)=\{uwP_{S\setminus\{s\}}^-\mid s\in I\}$ for some~$u\in U_+$ and~$w\in W$ and, by~\eqref{eq:intersParab}, its $U_+$-stabiliser is $u$-conjugated to~$U_+\cap wP_{S\setminus I}^-w^{-1}$. Denote by~$w_1$ the element in~$wW_{S\setminus I}$ of minimal Coxeter length and observe that~$U_+\cap wP_{S\setminus I}^-w^{-1}=U_+\cap w_1P_{S\setminus I}^-w_1^{-1}=w_1U_{S\setminus I}w_1^{-1}U_{w_1}$. The latter equality follows from \Cref{lem:UJUPinters}, which is applicable because~$\ell(w_1s)>\ell(w_1)$ for every~$s\in S\setminus I$. By~\cite[Lemma~5.8]{CapRem:root}, $U_{w_1}=U_{\alpha_1}\cdots U_{\alpha_\ell}$ for some~$\alpha_1,\dots, \alpha_\ell\in \Phi$. Hence~$U_{w_1}$ is finite. Moreover, if~$I\neq \emptyset$ then $U_{S\setminus I}$ is finite by \Cref{cor:UJfinite}.
[If~$T$ is finite, the finiteness of~$U_+\cap wP_J^-w^{-1}$ for every~$w\in W$ and every spherical subset~$J\subseteq S$ can be directly deduced from~\cite[Lemma~3.4]{AbGates}.]

We also observe that, for every~$d\geq 0$,~$U_+$ has finitely many orbits on~$\Delta_d$. In detail, by \Cref{cor:refBirk} and \Cref{lem:hl} we have
  \begin{equation*}
      \Delta_d =U_+\cdot\{\sigma(w,I)\mid h(\sigma(w,I))\leq d,\,w\in W,\,I\subseteq S\}=U_+\cdot\{\sigma(w,I)\mid \ell(w)\leq d,\,w\in W,\,I\subseteq S\}.
  \end{equation*}
  Since~$S$ is a finite generating set of~$W$,~$\Delta_d$ is finite modulo~$U_+$.

For simplicity, we divide the remaining part of the argument into steps:
\begin{itemize}
    \item \textbf{Step~0}: \emph{For every~$\rho\in \bigcup_{d\geq 0}R_{d+1}$,~$\dim(T(\rho))=k$ and~$T(\rho)$ is $(k-1)$-acyclic over~$\Q$.}
\end{itemize}
\textit{Proof of Step~0}. Clearly,~$\dim(\Delta_d)=k+1$ for every~$d\geq 0$. Moreover, for every~$\rho\in R_{d+1}$, since~$T(\rho)=\partial \rho\ast T'(\rho)$ and~$T'(\rho)$ is the link of~$\rho$ in~$\Delta_{d+1}$, we deduce that
  \begin{equation*}
\dim(T(\rho))=\dim(\partial\rho)+\dim(T'(\rho))+1=\dim(\Delta_{d+1})-1=k.
  \end{equation*}

Notice that~$\rho$ has spherical type in~$X_-$. By \Cref{prop:YconnSph} and \Cref{prop:descLinks}, $T'(\rho)$ is $(\dim(T'(\rho))-1)$-acyclic over~$\C$. By~\cite[\S 9.20]{bjorner}, we conclude that~$T(\rho)=\partial\rho\ast T'(\rho)$ is $(\dim(\rho)-1)$-acyclic.

\begin{itemize}
    \item \textbf{Step~1}: \textit{There exist infinitely many~$d\geq 0$ such that~$\dim(\rho)=k$ for some~$\rho\in R_{d+1}$. In that case,~$\wtil{\hH}_k(T(\rho),\Q)\neq\{0\}$.}
\end{itemize}
\textit{Proof of Step~1}. The first part of the statement was proved in~\cite[\S I.5, Lemma~13]{ab:twin}, and the second part was claimed therein. For completeness, we provide a proof of~$\wtil{\hH}_k(T(\rho),\Q)\neq\{0\}$. 

Since~$\rho$ has codimension~$1$ in~$\Delta_{d+1}$ and~$\Delta_{d+1}$ is locally finite, the $0$-dimensional link~$T'(\rho)=(\Delta_{d+1})_\rho$ consists of finitely many vertices~$\{v_1,\dots, v_n\}$. Since~$(W,S)$ is non-spherical, $\rho=\sigma(uw,I(w))$ for some~$u\in U_+$, $w\in W$ such that~$\ell(w)=d+1$ (see \Cref{lem:restr}\eqref{lem:restr1}) and such that~$I(w)\subseteq S$ is a proper non-empty subset (see~\eqref{eq:I(w)prop} and the lines below). Hence, every chamber of the form~$uwk\cdot c_0^-$ for some~$k\in P_{I(w)}^-$ contains~$\rho$. Notice that there are exactly~$|P_{I(w)}^-:B_-|\geq 2$ chambers like this. It follows that~$|T'(\rho)|\geq 2$. Thus,
\begin{equation*}
    T(\rho)=\partial\rho\ast \{v_1,\dots, v_n\}=\displaystyle{\bigcup_{i=1}^n}D_i,\qquad D_i:=\partial\rho\ast\{v_i\}=\{\sigma\mid \sigma\subseteq\rho\sqcup\{v_i\}\}.
\end{equation*}
Each~$D_i$ is contractible, being a cone. Moreover, $|\partial\rho|\simeq \mathbb{S}^{k-1}$ because~$\dim(\rho)=k$ by assumption. Using the Mayer--Vietoris sequence in relative homology~\cite[p.~71]{GeoBook}, we deduce that
\begin{equation*}
    \wtil{\hH}_{k}(D_1\cup D_2,\Q)\cong \wtil{\hH}_{k-1}(D_1\cap D_2,\Q)=\wtil{\hH}_{\dim(\partial\rho)}(\partial\rho,\Q)\cong \Q
\end{equation*}
and, for every~$1\leq \ell\leq n-1$, since~$(D_1\cup\dots\cup D_\ell)\cap D_{\ell-1}=\partial\rho$, we get a monomorphism
\begin{equation*}
    \wtil{\hH}_k(D_1\cup \dots\cup D_{\ell},\Q)\cong  \wtil{\hH}_k(D_1\cup \dots\cup D_{\ell},\Q)\oplus \wtil{\hH}_k(D_{\ell+1})\hookrightarrow \wtil{\hH}_k(D_1\cup\dots \cup D_{\ell+1},\Q).
\end{equation*}
This yields~$\wtil{\hH}_k(T(\rho),\Q)=\wtil{\hH}_k(D_1\cup\dots\cup D_n,\Q)\neq\{0\}$.

\begin{itemize}
    \item \textbf{Step~2}: \textit{For every~$d\geq 0$, the inclusion~$\Delta_d\hookrightarrow \Delta_{d+1}$ induces a surjective map~$\wtil{\hH}_k(\Delta_d,\Q)\twoheadrightarrow \wtil{\hH}_k(\Delta_{d+1},\Q)$ whose kernel contains an isomorphic copy of~$\wtil{\hH}_k(T(\rho),\Q)$ for every~$\rho\in R_{d+1}$.}
\end{itemize}
\textit{Proof of Step~2}. Choose~$\rho\in R_{d+1}$ and label the elements of the countable set~$R_{d+1}$ using the positive integers as~$\rho_1=\rho, \rho_2, \dots\,\,$.~(Recall that already~$X_-$ is countable, being locally finite). Consider the filtration~$\Delta_{d+1}=\bigcup_{m\geq 0}Z_m$ given by~$Z_0:=\Delta_{d}$ and~$Z_m:=\Delta_d\cup \bigcup_{\ell=1}^mS(\rho_i)$. 
For all~$\ell\geq 0$ and~$m\geq 1$ we notice that
\begin{equation}\label{eq:Zmrel}
    {\hH}_\ell(Z_{m},Z_{m-1};\Q)\cong \hH_\ell(S(\rho_m),T(\rho_m);\Q)\cong \wtil{\hH}_{\ell-1}(T(\rho_m),\Q).
\end{equation}
The first isomorphism before follows by excision~\cite[Theorem~2.8.7]{GeoBook}, while the second from the long exact sequence in relative reduced homology for~$(S(\rho_m),T(\rho_m))$~\cite[p.~71]{GeoBook} and the contractibility of~$S(\rho_m)$.
Moreover, from the long exact sequence in relative reduced homology for~$(Z_m,Z_{m-1})$ and~\eqref{eq:Zmrel} we get the following long exact sequence:
\begin{equation}\label{eq:niceLES}
     \xymatrix{\wtil{\hH}_{k+1}(Z_{m-1},\Q)\ar@{^{(}->}[r]^-{j_m^{k+1}} & \wtil{\hH}_{k+1}(Z_m,\Q)\ar[r] & \wtil{\hH}_{k}(T(\rho_m),\Q)\ar[r] & \wtil{\hH}_{k}(Z_{m-1},\Q)\ar@{->>}[r]^-{j_m^k} & \wtil{\hH}_k(Z_m,\Q)}
\end{equation}
where the labelled morphisms are induced by inclusions. The monomorphism and the epimorphism are due to Step~0.

By~\eqref{eq:niceLES} and standard properties of direct limits~\cite[Proposition~1.2.4]{ribzal}, the inclusion~$\Delta_d=Z_0\hookrightarrow \Delta_{d+1}=\bigcup_{m\geq 0}Z_m$ induces an epimorphism
\begin{equation*}
\wtil{\hH}_k(\Delta_d,\Q)=\wtil{\hH}_k(Z_0,\Q)\twoheadrightarrow \varinjlim_{m\geq 0}\wtil{\hH}_k(Z_m,\Q)\cong \wtil{\hH}_k(\Delta_{d+1},\Q)
\end{equation*}
with kernel containing the kernel of~$j_1^k\colon \wtil{\hH}_k(Z_0,\Q)\to \wtil{\hH}_k(Z_1,\Q)$, and a monomorphism
\begin{equation*}
    \wtil{\hH}_{k+1}(\Delta_d,\Q)\hookrightarrow\varinjlim_{m\geq 0}\wtil{\hH}_{k+1}(Z_m,\Q)\cong \wtil{\hH}_{k+1}(\Delta_{d+1},\Q).
\end{equation*}
In turn, arguing analogously, the latter monomorphism induces a monomorphism
\begin{equation*}
    \wtil{\hH}_{k+1}(\Delta_d,\Q)\hookrightarrow\varinjlim_{d'\geq 0}\wtil{\hH}_{k+1}(\Delta_{d'},\Q)\cong \wtil{\hH}_{k+1}(X_-,\Q).
\end{equation*}
Since~$X_-$ is contractible, we obtain that~$\wtil{\hH}_{k+1}(\Delta_d,\Q)=\{0\}$ for every~$d\geq 0$. As a consequence, for every~$d\geq 0$, for each term~$Z_m$ of the filtration~$\Delta_{d+1}=\bigcup_{m\geq 0}Z_m$,~$\wtil{\hH}_{k+1}(Z_m,\Q)$ vanishes (as it embeds in~$\wtil{\hH}_{k+1}(\Delta_{d+1},\Q)$). By~\eqref{eq:niceLES}, the kernel of~$j_1^k$ is isomorphic to~$\wtil{\hH}_k(T(\rho_1),\Q)$.

\begin{itemize}
    \item \textbf{Step~3}: \textit{The sequence~$(\wtil{\hH}_k(\Delta_d),\Q)_{d\geq 0}$ is not essentially trivial. In particular, $U_+$ is not of type~$\FP_{k+1}(\Q)$.}
\end{itemize}
\textit{Proof of Step~3}. Assume for a contradiction that~$(\wtil{\hH}_k(\Delta_d),\Q)_{d\geq 0}$ is essentially trivial. Thus there is~$d_0\geq 1$ such that, for all~$d\geq d_0$, the map~$\wtil{\hH}_k(\Delta_0,\Q)\to \wtil{\hH}_k(\Delta_d,\Q)$  induced by the inclusion is zero. By Step~2, the latter map is surjective for every~$d\geq 0$. Hence~$\wtil{\hH}_k(\Delta_d,\Q)=\{0\}$ for every~$d\geq d_0$. However, by Step~1, we can always find~$d\geq d_0$ such that the map~$\wtil{\hH}_k(T(\rho),\Q)\neq\{0\}$ for some~$\rho\in R_{d+1}$. But this is incompatible with what we get from Step~$2$ and the previous vanishing statements. This proves that~$(\wtil{\hH}_k(\Delta_d),\Q)_{d\geq 0}$ is not essentially trivial. By \Cref{thm:BroCrit} applied to the $U_+$-action on~$X_-=\bigcup_{d\geq 0}\Delta_d$, we conclude that~$U_+$ is not of type~$\FP_{k+1}(\Q)$.
\end{proof}

\begin{thm}\label{thm:GFP}
  Let~$G$ be a group with a locally finite root group datum~$(U_\alpha)_{\alpha\in \Phi}$ of type~$(W,S)$ with finite subgroup~$T$ as in~\eqref{eq:TBN}, and assume that~\hyperref[eq:co]{\co} and~\eqref{eq:thickSpec2}. Let~$(W,S)$ be $(k+1)$-spherical for some~$k\geq 1$. Then the following hold for every~$J\subseteq S$:
  \begin{enumerate}
      \item \label{thm:GFP1} $P_J^+$ is of type~$\FP_k(\Q)$;
      \item \label{thm:GFP2} If~$(W,S)$ is not spherical, then~$G$ is not of type~$\FP_{2\cdot\sigma(W,S)-2}(\Q)$;
      \item \label{thm:GFP3} Assume that~$(W,S)$ is classical and that~$\qmin$ satisfies~\eqref{eq:Yhomot}. Then~$P_J^+$ is of type~$\tF_k$.
  \end{enumerate}
\end{thm}
\begin{proof}
 Let~$(W,S)$ be $(k+1)$-spherical and let~$J\subseteq S$. 
 
  \eqref{thm:GFP1} If~$J$ is a spherical set, then~$P_J^+$ contains~$B_+$ and hence~$U_+$ as finite index subgroups; see \Cref{rem:GBspher} and \Cref{thm:twinBNassoc}\eqref{thm:twinBN1}. By \Cref{prop:permFP}, $P_J^+$ is of type~$\FP_k(\Q)$ whenever~$J$ is spherical. In particular, $P_J^+$ is of type~$\FP_k(\Q)$ whenever~$|J|\leq k+1$.

We now prove that~$P_J^+$ is of type~$\FP_k(\Q)$ by induction on~$|J|\geq k+1$. The base case~$|J|=k+1$ is already proved. 
Let~$|J|>k+1$ and assume that~$P_{J'}^+$ is of type~$\FP_k(\Q)$ for all~$J'\subseteq S$ with~$|J'|<|J|$. 
Recall that~$P_J^+$ acts on the building~$\Delta(P_J^+,B_+)$ coming the BN-pair~$(B_+,N_J)$ of type~$(W_J,J)$ induced from~$(B,N)$ (see \Cref{sus:BNpairs}). This action has finitely many orbits and the stabilisers of $\ell$-dimensional simplices are isomorphic to~$P_{J\setminus I}^+$ for some~$I\subseteq J$ with~$|I|=\ell+1$. By Solomon--Tits theorem, $\Delta(P_J^+,B_+)$ is $(|J|-1)$-acyclic over~$\Z$ (and hence over~$\Q$) and it is contractible whenever~$J$ is not spherical. By~\Cref{prop:BrownSuff},~$P_J^+$ is of type~$\FP_k(\Q)$.
 
\eqref{thm:GFP2} The claim is well-known and follows from a particular case of Gandini's theorem~\cite{gandini}, namely~\cite[Proposition~3.2]{AbGates}. In detail, apply the latter proposition to the case when~$Z(\Delta(G,B_\pm))$ is the Davis realisation of the building~$\Delta(G,B_\pm)$, i.e., the flag complex of all the spherical special cosets~$\{gP_J^+\mid g\in G,\,J\subseteq S\text{ spherical}\,\}$ ordered by inclusion. Clearly,~$Z(\Delta(G,B_{\pm}))$ is $(\sigma(W,S)-1)$-dimensional.

\eqref{thm:GFP3} Since being of type~$\tF_1$ is equivalent to being of type~$\FP_1(\Q)$, we may assume~$k\geq 2$. 
By \Cref{prop:FPnFn}\eqref{prop:FPnFn2}, it suffices to prove that~$P_J^+$ is finitely presented and of type~$\FP_k(\Z)$. In both cases, the strategy to prove part~\eqref{thm:GFP1} is transferred verbatim with only one minor change: to prove that~$P_J^+$ is finitely presented, apply \Cref{prop:FPnFn}\eqref{prop:FPnFn1} instead of \Cref{prop:BrownSuff}.
\end{proof}

\begin{cor}\label{cor:TBLFP}
     Let~$k\in \Z_{\geq 2}$ and let~$\Lambda$ be a twin building lattice of $(k+1)$-spherical type~$(W,S)$ for some~$k\geq 1$, whose associated root group datum~$(G,(U_\alpha)_{\alpha\in \Phi})$ satisfies~\hyperref[eq:co]{\co} and~\eqref{eq:thickSpec2}. Then~$\Lambda$ is of type~$\FP_k(\Q)$ and, if $(W,S)$ is non-spherical, not of type~$\FP_{2\cdot\sigma(W,S)-2}(\Q)$.
 If in addition~$(W,S)$ is classical and~$\qmin$ satisfies~\eqref{eq:Yhomot}, then~$\Lambda$ is of type~$\tF_k$.
\end{cor}
\begin{proof}
    The group~$\Lambda$ is isomorphic to~$G/Z(G)$ and the centre~$Z(G)$ is a finite group (see \Cref{rem:TBL}). The claim now follows from~\Cref{prop:permFP} and \Cref{thm:GFP}.
\end{proof}

The negative part of the statement in \Cref{thm:GFP} puts restrictions on the homological acyclicity of the negative opposite complex~$Y_-$. 
\begin{cor}\label{cor:OpXhomot}
    Let~$G$ be a group with a locally finite root group datum~$(U_\alpha)_{\alpha\in \Phi}$ of type~$(W,S)$ with finite subgroup~$T$ and that satisfies~\hyperref[eq:co]{\co} and~\eqref{eq:thickSpec2}. Assume also that~$(W,S)$ is~$(|S|-1)$-spherical but not spherical. Then~$Y_-$ is not~$(|S|-2)$-acyclic over~$\Q$. 
\end{cor}
\begin{proof}
Recall that~$Y$ is $U_+$-equivariant simplicially isomorphic to~$\caC\caC(U_+,(U_{S\setminus\{s\}})_{s\in S})$ and that~$U_J$ is finite whenever~$J\subseteq S$ is spherical; see \Cref{cor:Y-Uchar} and \Cref{cor:UJfinite}.
Hence~$U_+$ acts simplicially, properly, and with finitely many orbits on~$Y_-$. If~$Y_-$ is $(|S|-2)$-acyclic over~$\Q$, by \Cref{prop:BrownSuff} we would have that~$U_+$ is of type~$\FP_{|S|-1}(\Q)$, contradicting \Cref{thm:UFP}\eqref{thm:UFP2}.
\end{proof}

Below we record another straightforward consequence of the negative part of \Cref{thm:GFP}.
\begin{cor}\label{cor:homotRes}
    Let~$G$ be a group with a locally finite root group datum~$(U_\alpha)_{\alpha\in \Phi}$ of non-spherical $2$-spherical type~$(W,S)$, with finite subgroup~$T$ and such that~\hyperref[eq:co]{\co} and~\eqref{eq:thickSpec2} hold. Consider the flag complex of the following poset ordered by entry-wise inclusion:
   \begin{equation*}
       \big\{(gP_I^+,gP_J^-)\mid g\in G,\,I,J\subseteq S\text{ spherical}\,\big\}.
   \end{equation*}
   Then~$Z$ is not $(2\cdot\sigma(W,S)-3)$-acyclic over~$\Q$.
\end{cor}
\begin{proof}
    The group~$G$ acts simplicially, properly (see~\cite[Lemma~3.4]{AbGates}) and with finitely many orbits on~$Z$. Let~$n=2\cdot\sigma(W,S)-2$. If~$Z$ is $(n-1)$-acyclic over~$\Q$, then $Z$ is an $n$-good $G$-CW complex. By \Cref{prop:BrownSuff},~$G$ is of type~$\FP_n(\Q)$, which contradicts \Cref{thm:GFP}\eqref{thm:GFP2}. 
\end{proof}

\printbibliography

\end{document}